\documentclass[landscape]{article}
\usepackage{amsmath,amssymb,latexsym,amsthm}

\usepackage{longtable}

\usepackage[english]{babel}

\newtheorem{defn}{Definition}[section]
\newtheorem{lemma}[defn]{Lemma}

\newtheorem{prop}[defn]{Proposition}

\theoremstyle{definition}
\newtheorem{ex}{Example}[subsection]
\newtheorem{rem}[ex]{Remark}

\newcommand{\h}{{\cal H}}

\newcommand{\mn}{\mathbb N}

\newcommand{\neweq}{\overset{\scriptscriptstyle{\nabla}}{=}} 

\def\range{{\cal R}}

\def\h{{\cal H}}

\def\bp{\noindent{\bf Proof: \ }}
\def\ep{\noindent{$\Box$}}

\def\<{\langle}
\def\>{\rangle}

\newcommand{\seq[1]}{(#1_n)}

\def\newin {\,\kern-0.4em\in\kern-0.15em}
\def\newsubset {\kern-0.2em\subset\kern-0.2em}

\def\normsn{$\|\!\cdot\!\|$-$SN$ }

\usepackage{times}
\usepackage[dvips,colorlinks=true,urlcolor=blue,bookmarks=false,pdfpagemode=None]{hyperref}
\usepackage{color}

\definecolor{xxl}{rgb}{0.2,0.08,0.4}

\title{Detailed characterization of unconditional convergence \\ and invertibility of multipliers}
\author{
 D.\,T. Stoeva$^{a)\,b)}$ and P. Balazs$^{b)}$ \\ 
$^{a)}$ Department of Mathematics, \\ University of Architecture, Civil Engineering and Geodesy,\\
Blvd Christo Smirnenski 1, 1046 Sofia, Bulgaria\\
$^{b)}$ Acoustics Research Institute, \\
Wohllebengasse 12-14, Vienna A-1040, Austria \\
}

\begin{document}
\maketitle \pagestyle{myheadings} 
{\footnotetext[1]{This work was supported by  the WWTF project MULAC ('Frame Multipliers: Theory and Application in Acoustics; MA07-025)} }

\begin{abstract} In this paper we investigate the possibility of unconditional convergence and invertibility of multipliers $M_{m,\Phi,\Psi}$ depending on the properties of the sequences $\Psi$,$\Phi$ and $m$. 
 We characterize a complete set of conditions for the invertibility and the unconditional convergence of multipliers, and collect those results in tables. 
We either prove that unconditional convergence and invertibility is not possible, that one or both of these conditions are always the case for the given parameters, or we give examples for the feasible combinations. 
We give a full list of examples for all conditions.
\end{abstract}



\section{Introduction}

In \cite{schatt1}, R. Schatten provided a detailed study of ideals of compact operators using their singular decomposition. He
investigated the operators of the form $\sum_{k}\lambda_{k}\varphi_{k}\otimes\overline{\psi_{k}}$ where
$(\phi_{k})$ and $(\psi_{k})$ are orthonormal families. Later such operators were investigated for Gabor frames \cite{feinow1}. In \cite{xxlmult1} abstract Bessel
and frame sequences were used to define Bessel and frame multipliers, where several basic properties of frame multipliers were investigated. Recently this concept was extended to $p$-frames in Banach spaces\cite{rahxxlXX} and generalized frames \cite{rahgenmul10}.

While multipliers are interesting from a theoretical point of view \cite{benepfand07,xxlframehs07,doetor09} there are also interesting for applications, in particular in the fields of audio and acoustics.
Time-invariant filters are often used for audio applications. These systems can be described by the multiplication of the frequency spectrum of the signal by a fixed function. Using the Fourier transform to calculate the spectrum, such an operator can be called a {\em Fourier multiplier} \cite{fena06,bgrgrok05}. 
As a particular way to implement time-frequency filters or time-variant filtering, Gabor (frame) multipliers \cite{feinow1} can be used.
In signal processing they are used 
 under the name
'Gabor filters' \cite{hlawatgabfilt1}. 
\\

In the present table we investigate the unconditional convergence and the invertibility of multipliers $M_{m,\Phi,\Psi}$. 
In \cite{balsto09} the focus was on existence results and formulas for the inversion. Here we
 give a complete set of examples varying the type of the sequences $\Phi=\seq[\phi]$, $\Psi=\seq[\psi]$ 
(non-Bessel, Bessel non-frames, frames non-Riesz bases, Riesz bases; norm-semi-normalized, non-norm-semi-normalized with all possible combinations)
and varying the symbol $m$ (semi-normalized, $\in\ell^\infty$ but non-semi-normalized, $\notin\ell^\infty$). 
We list all possible combinations.
We give a full classification, if multipliers under those conditions can be ('POSSIBLE'), have to be ('ALWAYS') or never can be ('NOT POSSIBLE') 
unconditionally convergent and invertible (resp. non-invertible) on the given Hilbert space.
Please note that examples for multipliers that are identical to the identity give examples for those cases, where sequences can be dual to each other.
We only consider sequences with non-zero elements, as in this case, for example, the invertible identity operator and the zero 
 operator can be described as multiplier, if zeros are put at appropriate places (see \cite{balsto09} for details.)\footnote{For future work such sequences including zeros could be linked to the excess of frames \cite{bacahela03,bacahela03-1}.}

The paper is organized as follows. 
In Section \ref{sec:prelnot0} we fix the notation in this paper, summarize known and prove some new results needed for the rest of the paper.
In \mbox{Section \ref{sec:tables0}} we determine if unconditional convergence and invertibility is always, sometimes or never possible for the complete set of possibilities for the sequences $\Phi=\seq[\phi]$ and  $\Psi=\seq[\psi]$: non-Bessel, Bessel non-frames, frames non-Riesz bases, Riesz bases, combined with all the possibilities for norm-boundedness;
and varying the symbol $m$ to be semi-normalized, bounded above or bounded below or non-bounded. 
These results are collected in tables, which are linked to Section \ref{sec:examples1}.
There  we give concrete examples for all the combinations in the tables.

Please note that in the electronic version of this paper, the examples and the tables are connected with hyperlinks.

\section{Preliminaries and Notations}\label{sec:prelnot0}

Throughout the whole paper, $\h$ denotes an 
 (infinite-dimensional)$^2$ Hilbert spaces. \footnotetext[2]{
All results will be valid for finite-dimensional spaces also, but some classifications might not be useful, as e.g. every finite sequence in a finite-dimensional vector space is a frame on its span.
}
The range of an operator $G$ is denoted by $\range(G)$. 
The identity operator on $\h$ is denoted by $I_\h$. The operator $G:\h\to \h$ is called {\it invertible on $\h$} if 
 there exists bounded operator $G^{-1}:\h\to \h$  such that $GG^{-1}=G^{-1}G=I_\h$. 

The notation $\Phi$ (resp. $\Psi$) is used to denote the sequence $\seq[\phi]$ (resp. $\seq[\psi]$) with elements from $\h$; 
$m$ denotes a 
 complex scalar sequence $\seq[m]$
and $\overline{m}$ denotes the sequence of the complex conjugates of $m_n$.
The sequence $\seq[e]$ denotes an orthonormal basis of $\h$. 
The index set of sequences will be often omitted, in such cases the set $\mn=\{1,2,3,\ldots\}$ is assumed implicitly.
   
Recall that $\Phi$ is called a {\it Bessel sequence} (in short, {\it Bessel}) {\it for $\h$ with bound $B_\Phi$} if $B_\Phi>0$ 
and $\sum |\<h,\phi_n\>|^2 \leq B_\Phi\|h\|^2$ for every $h\in\h$. A Bessel sequence $\Phi$ with bound $B_\Phi$  is called a {\it frame for $\h$ with bounds $A_\Phi, B_\Phi$}, if $A_\Phi>0$ and $A_\Phi\|h\|^2\leq  \sum |\<h,\phi_n\>|^2 $ for every $h\in\h$. 
The sequence $\Phi$ is called a Riesz basis for $\h$ with bounds $A_\Phi,B_\Phi$,  
if $A_\Phi>0$ and 
$A_\Phi \sum |c_n|^2\leq \|\sum c_n \phi_n\|^2\leq B_\Phi \sum |c_n|^2$,
$\forall \seq[c]\in\ell^2$.
Every Riesz basis for $\h$ with bounds $A,B$ is a frame for $\h$ with bounds $A,B$. 
For standard references for frame theory and related topics see \cite{Casaz1,ole1,he98-1}.

For a given Bessel sequence $\Phi$, the mapping $U_\Phi:\h\to \ell^2$ given by $U_\Phi f =(\<f,\phi_n\>)$ is called {\it the analysis operator for $\Phi$} and the mapping $T_\Phi$ given by $T_\Phi \seq[c]= \sum c_n \phi_n$   is called {\it the synthesis operator for $\Phi$}. 

Given sequences $m, \Phi,$ and $\Psi$, the operator $M_{m,\Phi,\Psi}$ given by $$M_{m,\Phi,\Psi}f=\sum_{n=1}^\infty m_n\<f,\psi_n\>\phi_n, f\in\h,$$ is called {\it multiplier}. 
The multiplier $M_{m,\Phi,\Psi}$ is called {\it unconditionally convergent on $\h$} (resp. {\it surjective}) if $\sum m_n\<f,\psi_n\>\phi_n$ converges unconditionally for every $f\in\h$ (resp. $\range(M_{m,\Phi,\Psi})=\h$).
The sequence $m$ is called {\it semi-normalized} if $0<\inf_n |m_n|\leq  \sup_n |m_n| <\infty$. The sequence $\Phi$ is called norm-bounded below (resp. norm-bounded above) if $\inf_n \|\phi_n\| >0$ (resp. $\sup_n \|\phi_n\| <\infty$) and $\Phi$ is called norm-semi-normalized if  $0<\inf_n \|\phi_n\|\leq  \sup_n \|\phi_n\| <\infty$.

To shorten the file we use the following abbreviations: R.b. - Riesz basis,
fr. - frame, 
B. - Bessel sequence, 
$SN$ - semi-normalized, 
\normsn - norm-semi-normalized, 
$NBB$ - norm-bounded below, 
$NBA$ - norm-bounded above, 
unc.\,conv. - unconditionally convergent on $\h$, 
INV. - invertible on $\h$.
Recall that a Riesz basis is always \normsn and a Bessel sequence is always $NBA$. That is why the tables do not include the impossible cases \lq\lq non-$NBA$ Bessel sequence\rq\rq\ and \lq\lq non-\normsn Riesz basis\rq\rq.

\subsection{Concerning Well-definedness and Unconditional Convergence}

First note that the sequence $m\Phi$ is Bessel, frame, Riesz basis or satisfies the lower frame condition for $\h$, if and only if $\overline{m}\,\Phi$ is Bessel, frame, Riesz basis or satisfies the lower frame condition for $\h$, respectively \cite{balsto09}. Furthermore, if $m$ is $SN$, then the sequence $\Phi$ is Bessel (resp. frame, Riesz basis) for $\h$ if and only if $m\Phi$ is Bessel (resp. frame, Riesz basis) for $\h$.

\vspace{.1in} We will often use the following result:
\begin{prop} \label{bmh} {\rm \cite[Theorem 6.1 (1)]{xxlmult1}} 
Let $\Phi$ and $\Psi$ be Bessel sequences for $\h$. If $m\newin\ell^\infty$, then $M_{m,\Phi,\Psi}$ 
is well defined from $\h$ into $\h$, bounded and unconditionally convergent on $\h$.
\end{prop}

 To shorten notation, for $\nu=(\nu_n), \Theta=(\theta_n), \Xi=(\xi_n)$, we will write  
$M_{m,\Phi,\Psi} \neweq  M_{\nu,\Xi,\Theta}$ if there exist scalar sequences $(c_n)$, $(d_n)$ so that $\xi_n=c_n \phi_n$, $\theta_n=d_n\psi_n$ and $m_n=\nu_n c_n d_n$ for every $n$. This means that in the series the summands are the same element-wise.
We need this equality for conclusions of unconditional convergence. If $M_{\nu,\Xi,\Theta}$ is unconditionally convergent on $\h$ and $M_{m,\Phi,\Psi} \neweq  M_{\nu,\Xi,\Theta}$, then $M_{m,\Phi,\Psi}$ is clearly also unconditionally convergent on $\h$. Note that if $M_{\nu,\Xi,\Theta}$ is unconditionally convergent on $\h$ and $M_{m,\Phi,\Psi} =  M_{\nu,\Xi,\Theta}$, then $M_{m,\Phi,\Psi}$ might not be unconditionally convergent on $\h$. Consider for example the sequences $\Phi$ and $\Psi$ in Remark \ref{r1}(a) - the multiplier $M_{m,\Phi,\Psi}$ is equal to the identity operator and thus it can be written as $M_{(1),(e_n),(e_n)}$, which is unconditionally convergent on $\h$; however, $M_{m,\Phi,\Psi}$ is not unconditionally convergent on $\h$.

\begin{prop} \label{lem31} {\rm \cite[Prop. 3.3 and Lemma 3.5]{balsto09}} Let $M_{m,\Phi,\Psi}$ or $M_{m,\Psi,\Phi}$ be unconditionally convergent on $\h$. 
\begin{itemize}
\item[{\rm (i)}] If $\Phi$ (resp. $m\Phi$) is $NBB$, then $m\Psi$ (resp. $\Psi$) is Bessel for $\h$.
\item[{\rm (ii)}] 
If both $\Phi$ and $\Psi$ are $NBB$, then $m\in\ell^\infty$.
\end{itemize}
\end{prop}

\begin{prop} \label{lem36} {\rm \cite[Prop. 3.9]{balsto09}}
Let $\Phi$ be a Riesz basis and let $M_{m,\Phi,\Psi}$ (resp. $M_{m,\Psi,\Phi}$) be well defined on $\h$. Then the following holds.
\begin{itemize}
\item[{\rm (i)}] The sequence $m\Psi$ is Bessel for $\h$.
\item[{\rm (ii)}] If $\Psi$ is $NBB$, then $m\in\ell^\infty$.
 \end{itemize}
\end{prop}

\subsection{Concerning invertibility}

\begin{prop} \label{c1}  {\rm \cite[Theorem 4.3]{balsto09}}
Let $M_{m,\Phi,\Psi}$ be invertible on $\h$. 
If $\Psi$ and $\Phi$ are Bessel sequences for $\h$ and $m\in\ell^\infty$, then $\Psi$ and $\Phi$ are frames for $\h$.
\end{prop}

 \begin{prop} \label{rbinv} {\rm \cite[Prop. 7.7]{xxlmult1}} 
Let $\Phi$ and $\Psi$ be Riesz bases for $\h$. If $m$ is $SN$, then $M_{m,\Phi,\Psi}$ is invertible on $\h$.
\end{prop}

\begin{prop} \label{rc} {\rm \cite[Corollary 4.12]{balsto09}}
Let $\Phi$ be a Riesz basis for $\h$. Then $M_{m,\Phi,\Psi}$ (resp. $M_{m,\Psi,\Phi}$) 
is invertible on $\h$ if and only if $m\Psi$ is a Riesz basis for $\h$. This may happen only in the following cases:

\indent $\bullet$ $\Psi$ is a Riesz basis for $\h$ and $m$ is $SN$;

$\bullet$ $\Psi$ is non-$NBB$ and Bessel for $\h$, which is not a frame for $\h$, and $m$ is $NBB$, but not in $\ell^\infty$;

$\bullet$ $\Psi$ is non-$NBA$, $NBB$, and non-Bessel for $\h$, $m$ is non-$NBB$ and $m\in\ell^\infty$;

$\bullet$ $\Psi$ is non-$NBA$, non-$NBB$, and non-Bessel for $\h$, $m$ is non-$NBB$ and $m\notin\ell^\infty$.

\end{prop}

\begin{prop} \label{rbis}
Let $\Phi$ be a Riesz basis for $\h$, $\Psi$ be an overcomplete frame for $\h$ and $m$ be $SN$. Then 

\vspace{.05in}
(a) $M_{m,\Phi,\Psi}$ is injective, but not surjective.

\vspace{.05in}
(b) $M_{m,\Psi,\Phi}$ is surjective, but not injective.
\end{prop}

\bp First note that $m\Psi$ and $\overline{m}\,\Psi$ are also overcomplete frames for $\h$. 

(a) Let $M_{m,\Phi,\Psi}f=0$ for some $f\in\h$. 
Since $\Phi$ is a Riesz basis for $\h$, it follows that $\<f,\overline{m}_n\psi_n\>=0$, $\forall n$, see \cite[Theorem 6.1.1]{ole1}. Since $\overline{m}\,\Psi$ is complete in $\h$, 
we obtain that $f=0$. Therefore, $M_{m,\Phi,\Psi}$ is injective.
It is proved in \cite[Theorem 4.11]{balsto09} that $M_{m,\Phi,\Psi}$ is not surjective. 
  
(b) 
Since $\Phi$ is a Riesz basis for $\h$, $U_{\Phi}$ is a bijection of $\h$ onto $\ell^2$.  Since $m\Psi$ is a frame for $\h$, $T_{m\Psi}$ is surjective from $\ell^2$ onto $\h$. Therefore, $M_{m,\Psi,\Phi}=T_{m\Psi} U_{\Phi}$ is surjective. 
It is proved in \cite[Theorem 4.11]{balsto09} that $M_{m,\Psi,\Phi}$ is not injective. 
\ep

\begin{prop} \label{propnon}
Let $\Phi$ be a $NBB$ Bessel for $\h$, which is not a frame for $\h$. Then, for any $\Psi$ and any $m$, the multiplier $M_{m,\Phi,\Psi}$ (resp. $M_{m,\Psi,\Phi}$) can not be both unconditionally convergent on $\h$ and invertible on $\h$. 
\end{prop}
\bp
Assume that $M_{m,\Phi,\Psi}$ (resp. $M_{m,\Psi,\Phi}$) is unconditionally convergent on $\h$. Write $M_{m,\Phi,\Psi}=M_{(1),\Phi,\overline{m}\,\Psi}$ (resp. $M_{m,\Psi,\Phi}=M_{(1),m\Psi,\Phi}$). By Proposition \ref{lem31}(i), the sequence $\overline{m}\,\Psi$ (resp. $m\Psi$) is Bessel for $\h$. Now Proposition \ref{c1} implies that  $M_{(1),\Phi,\overline{m}\,\Psi}$ (resp. $M_{(1),m\Psi,\Phi}$) is not invertible on $\h$.
\ep

\begin{lemma} \label{lemg}
Let $G_k$ denote the multiplier  $M_{(\frac{1}{n^k}),(e_n),(e_n)} $, $k\in\mn$. 
Then $G_k$ is unconditionally convergent on $\h$ and not invertible on $\h$.
\end{lemma}
\bp
By Proposition \ref{bmh}, $G_k$ is well defined from $\h$ into $\h$ and unconditionally convergent on $\h$.  The multiplier $G_k$ is injective, but not surjective - for example, the element $\sum \frac{1}{n^k}e_n\in\h$ does not belong to the range of $G_k$.
\ep

\newpage
\section{Classification Tables}\label{sec:tables0} 

\vspace{.1in}

{\small 
\centerline{Table 1: two non-Bessel sequences} 
\vspace{.1in}

\begin{longtable}{|l|l|l|l|l|}

\hline   & & & & \\
\label{table1}

\!$\phi$\! -\! not B. \!\!\!\!& 
\!$\psi$\! -\! not B.  \!\!\!\!& 
$m$ - SN \hspace{.1in}&
$m\in\ell^\infty$, but non-$SN$ \hspace{.1in} &
$m\notin\ell^\infty $   \\
\hline   & & & &   \\

 &
 &
$M_{m,\Phi,\Psi}$,\,$M_{m,\Psi,\Phi}$ \hspace{.08in} $M_{m,\Phi,\Psi}$,\,$M_{m,\Psi,\Phi}$\! &
 $M_{m,\Phi,\Psi}$,\,$M_{m,\Psi,\Phi}$ \hspace{.08in} $M_{m,\Phi,\Psi}$,\,$M_{m,\Psi,\Phi}$\! & 
$M_{m,\Phi,\Psi}$,\,$M_{m,\Psi,\Phi}$ \hspace{.08in} $M_{m,\Phi,\Psi}$,\,$M_{m,\Psi,\Phi}$ \! \\

 &
 &
unc.\,conv.\,\& \hspace{.4in} unc.\,conv.\,\& &
unc.\,conv.\,\&\hspace{.43in} unc.\,conv.\,\& & 
unc.\,conv.\,\& \hspace{.4in} unc.\,conv.\,\& \\

 &
 &
INV. \hspace{.78in} NON-INV. &
INV. \hspace{.77in} NON-INV. & 
INV. \hspace{.78in} NON-INV. \\

 & &  &  & \\ 
\hline   

\endfirsthead

\hline   & & & & \\
{\footnotesize continued from}  & & & & \\
 {\footnotesize the previous page}  & & & & \\ 
\hline & & & & \\
\!$\phi$\! -\! not B. \!\!\!\!& 
\!$\psi$\! -\! not B.  \!\!\!\!& 
$m$ - SN \hspace{.1in}&
$m\in\ell^\infty$, but non-$SN$ \hspace{.1in} &
$m\notin\ell^\infty $   \\
\hline   & & & &   \\

 &
 &
$M_{m,\Phi,\Psi}$,\,$M_{m,\Psi,\Phi}$ \hspace{.08in} $M_{m,\Phi,\Psi}$,\,$M_{m,\Psi,\Phi}$\! &
 $M_{m,\Phi,\Psi}$,\,$M_{m,\Psi,\Phi}$ \hspace{.08in} $M_{m,\Phi,\Psi}$,\,$M_{m,\Psi,\Phi}$\! & 
$M_{m,\Phi,\Psi}$,\,$M_{m,\Psi,\Phi}$ \hspace{.08in} $M_{m,\Phi,\Psi}$,\,$M_{m,\Psi,\Phi}$ \! \\

 &
 &
unc.\,conv.\,\& \hspace{.4in} unc.\,conv.\,\& &
unc.\,conv.\,\&\hspace{.43in} unc.\,conv.\,\& & 
unc.\,conv.\,\& \hspace{.4in} unc.\,conv.\,\& \\

 &
 &
INV. \hspace{.78in} NON-INV. &
INV. \hspace{.77in} NON-INV. & 
INV. \hspace{.78in} NON-INV. \\

 & &  &  & \\
\hline  & &  &  & \\ 
\endhead

& & & & \\
& & & & {\footnotesize continued on the next page} \\
\hline
\endfoot

\endlastfoot

\hline
& &  &  & \\

$\|\!\cdot\!\|$-SN &
$\|\!\cdot\!\|$-SN & 
NOT POSSIBLE  \hspace{.11in} NOT POSSIBLE &
POSSIBLE \hspace{.42in}  POSSIBLE & 
NOT POSSIBLE  \hspace{.12in} NOT POSSIBLE  \\

&
&
\hspace{.07in}  $M_{m,\Phi,\Psi},M_{m,\Psi,\Phi}$ - not unc.\,conv., &
Example \ref{nbnb1}(i) \hspace{.16in} Example \ref{nbnb1}(ii)  &
\hspace{.07in}  $M_{m,\Phi,\Psi},M_{m,\Psi,\Phi}$ - not unc.\,conv., \\

  &
 &
\hspace{.3in} see Proposition \ref{lem31}(i) &
 &
\hspace{.3in} see Proposition \ref{lem31}(ii) \\
 
 &
 &
 Remark \ref{r1}(a) \hspace{.15in} Remark \ref{r1}(b) &
 &
 \\

& &  &  & \\
\hline
& &  &  & \\

$\|\!\cdot\!\|$-SN &
$NBA$ \& & 
NOT POSSIBLE  \hspace{.115in} NOT POSSIBLE &
POSSIBLE \hspace{.42in} POSSIBLE  & 
POSSIBLE \hspace{.41in} POSSIBLE  \\

&
non-$NBB$ &
\hspace{.07in}  $M_{m,\Phi,\Psi},M_{m,\Psi,\Phi}$ - not unc.\,conv., &
Example \ref{nbnb9}(i) \hspace{.16in} Example \ref{nbnb10} &
Example \ref{nbnb9}(ii) \hspace{.11in} Example \ref{nbnb9}(iii) \\

  &
 &
\hspace{.3in} see Proposition \ref{lem31}(i) &
 &
 \\
 
 &
 &
Remark \ref{r11}(a) \hspace{.15in} Remark \ref{r11}(b) &
 &
 \\
 
  & & & & \\
\hline
 & & & & \\

$\|\!\cdot\!\|$-SN &
non-$NBA$  $\&$ & 
 NOT POSSIBLE  \hspace{.1in} NOT POSSIBLE &
POSSIBLE \hspace{.43in} POSSIBLE  & 
NOT POSSIBLE \hspace{.11in} NOT POSSIBLE  \\

&
$NBB$ &
\hspace{.07in}  $M_{m,\Phi,\Psi},M_{m,\Psi,\Phi}$ - not unc.\,conv.,   &
Example \ref{nbnb21}(i) \hspace{.17in} Example \ref{nbnb21}(ii) &
\hspace{.07in}  $M_{m,\Phi,\Psi},M_{m,\Psi,\Phi}$ - not unc.\,conv., \\

 &
 &
\hspace{.3in} see Proposition \ref{lem31}(i) &
 &
\hspace{.3in} see Proposition \ref{lem31}(ii) \\

  & & & & \\
\hline
 & & & & \\

$\|\!\cdot\!\|$-SN &
non-$NBA$  $\&$ & 
 NOT POSSIBLE  \hspace{.1in} NOT POSSIBLE &
POSSIBLE \hspace{.43in} POSSIBLE  & 
POSSIBLE \hspace{.4in} POSSIBLE  \\

&
non-$NBB$ &
\hspace{.07in}  $M_{m,\Phi,\Psi},M_{m,\Psi,\Phi}$ - not unc.\,conv.,   &
Example \ref{nbnbnew1}(i) \hspace{.17in} Example \ref{nbnbnew2} &
Example \ref{nbnbnew1}(ii) \hspace{.1in} Example \ref{nbnbnew1}(iii) \\

 &
 &
\hspace{.3in} see Proposition \ref{lem31}(i) &
 &
\\

\hline
 & & & & \\

$NBA$ \& &
$NBA$ \&  & 
 POSSIBLE  \hspace{.4in}  POSSIBLE &
POSSIBLE \hspace{.43in} POSSIBLE  & 
POSSIBLE \hspace{.4in} POSSIBLE   \\

non-$NBB$ &
non-$NBB$  &
Example \ref{nbnbnew5}(i) \hspace{.12in} Example \ref{nbnbnew5}(ii) &
Example \ref{nbnbnew6} \hspace{.21in} Example \ref{nbnbnew7}(i) &
Example \ref{nbnbnew7}(ii) \hspace{.05in} Example \ref{nbnbnew7}(iii) \\

& & & & \\
\hline
 & & & & \\

$NBA$ \& &
non-$NBA$ \& & 
 NOT POSSIBLE  \hspace{.14in} NOT POSSIBLE&
POSSIBLE \hspace{.43in} POSSIBLE  & 
POSSIBLE \hspace{.45in} POSSIBLE   \\

non-$NBB$ &
$NBB$ &
\hspace{.07in}  $M_{m,\Phi,\Psi},M_{m,\Psi,\Phi}$ - not unc.\,conv., &
Example \ref{nbnbnew8}(i) \hspace{.101in} Example \ref{nbnbnew8}(ii) &
Example \ref{nbnb1718}(i) \hspace{.12in} Example \ref{nbnb1718}(ii) \\

 &
 &
\hspace{.3in} see Proposition \ref{lem31}(i) &
 &
 \\

 & & & & \\
\hline
 & & & & \\

$NBA$ \& &
non-$NBA$ \& & 
POSSIBLE  \hspace{.435in}  POSSIBLE &
POSSIBLE \hspace{.43in} POSSIBLE  & 
POSSIBLE \hspace{.45in} POSSIBLE   \\

non-$NBB$ &
non-$NBB$ &
Example \ref{nbnb13}(i) \hspace{.1in} Example \ref{nbnb14}(i) &
Example \ref{nbnbnew9} \hspace{.21in} Example \ref{nbnb13}(ii) &
Example \ref{nbnb14}(ii) \hspace{.09in} Example \ref{nbnbnew9new} \\

 & & & & \\
\hline
 & & & & \\

non-$NBA$ \& &
non-$NBA$ \& & 
NOT POSSIBLE  \hspace{.15in} NOT POSSIBLE &
POSSIBLE \hspace{.42in}  POSSIBLE & 
NOT POSSIBLE  \hspace{.16in} NOT POSSIBLE  \\

$NBB$ &
$NBB$ &
\hspace{.07in}  $M_{m,\Phi,\Psi},M_{m,\Psi,\Phi}$ - not unc.\,conv., &
Example  \ref{nbnb15i} \hspace{.21in} Example  \ref{nbnb15ii}  &
\hspace{.07in}  $M_{m,\Phi,\Psi},M_{m,\Psi,\Phi}$ - not unc.\,conv., \\

  &
 &
\hspace{.3in} see Proposition \ref{lem31}(i) &
 &
\hspace{.3in} see Proposition \ref{lem31}(ii) \\

 & & & & \\
\hline
 & & & & \\

non-$NBA$ \& &
non-$NBA$ \& & 
NOT POSSIBLE  \hspace{.15in} NOT POSSIBLE &
POSSIBLE \hspace{.43in} POSSIBLE  & 
POSSIBLE \hspace{.45in} POSSIBLE   \\

$NBB$ &
non-$NBB$ &
\hspace{.07in}  $M_{m,\Phi,\Psi},M_{m,\Psi,\Phi}$ - not unc.\,conv., &
Example \ref{nbnbnew10}(i) \hspace{.1in} Example \ref{nbnbnew10}(ii) &
Example \ref{nbnbnew11}(i)  \hspace{.12in} Example \ref{nbnbnew11}(ii) \\
 
  &
 &
\hspace{.3in} see Proposition \ref{lem31}(i) &
 &
 \\

 & & & & \\
\hline
 & & & & \\

 non-$NBA$ \& &
non-$NBA$ \&  & 
 POSSIBLE  \hspace{.435in}  POSSIBLE &
POSSIBLE \hspace{.43in} POSSIBLE  & 
POSSIBLE \hspace{.45in} POSSIBLE   \\

non-$NBB$  &
non-$NBB$  &
Example \ref{nbnbnew12} \hspace{.22in} Example \ref{nbnbnew13} &
Example \ref{nbnbnew14}(i) \hspace{.1in} Example \ref{nbnbnew14}(ii) &
Example \ref{nbnbnew15}(i) \hspace{.12in} Example \ref{nbnbnew15}(ii) \\
 
  & & & & \\
\hline
\end{longtable}
}

\newpage
{ \small
\centerline{Table 2: one Bessel not a frame, one non-Bessel sequence}
\vspace{.1in}

\begin{longtable}{|l|l|l|l|l|}

\hline   & & & & \\
\label{table2}

\!$\phi$\! -\! B. \!\!\!\!& 
\!$\psi$\! -\! not B. \!\!\!\!& 
$m$ - SN \hspace{.1in}&
$m\in\ell^\infty$, but non-$SN$ \hspace{.1in} &
$m\notin\ell^\infty $   \\

not fr. \!\!\!\!& 
& 
&
&
\\

\hline   & & & &   \\

 &
 &
$M_{m,\Phi,\Psi}$,\,$M_{m,\Psi,\Phi}$ \hspace{.11in} $M_{m,\Phi,\Psi}$,\,$M_{m,\Psi,\Phi}$\! &
 $M_{m,\Phi,\Psi}$,\,$M_{m,\Psi,\Phi}$ \hspace{.16in} $M_{m,\Phi,\Psi}$,\,$M_{m,\Psi,\Phi}$\! & 
$M_{m,\Phi,\Psi}$,\,$M_{m,\Psi,\Phi}$ \hspace{.12in} $M_{m,\Phi,\Psi}$,\,$M_{m,\Psi,\Phi}$ \! \\

 &
 &
unc.\,conv.\,\& \hspace{.42in} unc.\,conv.\,\& &
unc.\,conv.\,\&\hspace{.51in} unc.\,conv.\,\& & 
unc.\,conv.\,\& \hspace{.43in} unc.\,conv.\,\& \\

 &
 &
INV. \hspace{.79in} NON-INV. &
INV. \hspace{.85in} NON-INV. & 
INV. \hspace{.81in} NON-INV. \\

& &  &  & \\

\endfirsthead

\hline   & & & & \\
{\footnotesize continued from}  & & & & \\
 {\footnotesize the previous page}  & & & & \\

\hline   & & & & \\

\!$\phi$\! -\! B. \!\!\!\!& 
\!$\psi$\! -\! not B. \!\!\!\!& 
$m$ - SN \hspace{.1in}&
$m\in\ell^\infty$, but non-$SN$ \hspace{.1in} &
$m\notin\ell^\infty $   \\

not fr. \!\!\!\!& 
& 
&
&
\\

\hline   & & & &   \\

 &
 &
$M_{m,\Phi,\Psi}$,\,$M_{m,\Psi,\Phi}$ \hspace{.12in} $M_{m,\Phi,\Psi}$,\,$M_{m,\Psi,\Phi}$\! &
 $M_{m,\Phi,\Psi}$,\,$M_{m,\Psi,\Phi}$ \hspace{.17in} $M_{m,\Phi,\Psi}$,\,$M_{m,\Psi,\Phi}$\! & 
$M_{m,\Phi,\Psi}$,\,$M_{m,\Psi,\Phi}$ \hspace{.12in} $M_{m,\Phi,\Psi}$,\,$M_{m,\Psi,\Phi}$ \! \\

 &
 &
unc.\,conv.\,\& \hspace{.42in} unc.\,conv.\,\& &
unc.\,conv.\,\&\hspace{.51in} unc.\,conv.\,\& & 
unc.\,conv.\,\& \hspace{.43in} unc.\,conv.\,\& \\

 &
 &
INV. \hspace{.79in} NON-INV. &
INV. \hspace{.85in} NON-INV. & 
INV. \hspace{.81in} NON-INV. \\

& &  &  & \\

\endhead

\hline  
& & & & \\
& & & & {\footnotesize continued on the next page} \\
\hline
\endfoot

\endlastfoot 

\hline
& &  &  & \\

$\|\!\cdot\!\|$-SN & 
$\|\!\cdot\!\|$-SN & 
NOT POSSIBLE  \hspace{.16in} NOT POSSIBLE &
NOT POSSIBLE \hspace{.23in} POSSIBLE  & 
NOT POSSIBLE  \hspace{.18in} NOT POSSIBLE  \\

&
& 
\hspace{.07in}  $M_{m,\Phi,\Psi},M_{m,\Psi,\Phi}$ - not unc.\,conv.,&
see Prop. \ref{propnon}  \hspace{.43in} Example \ref{bnb4}(i)  &
 \ \ \  $M_{m,\Phi,\Psi}$, $M_{m,\Psi,\Phi}$ - not unc.\,conv.,\\

&
& 
\hspace{.2in} \ \ \  see Proposition \ref{lem31}(i)  &
&
\hspace{.36in} see Proposition \ref{lem31}(ii) \\

& &  &  & \\
\hline
& &  &  & \\

$\|\!\cdot\!\|$-SN & 
$NBA$ $\&$ &
NOT POSSIBLE  \hspace{.19in} NOT POSSIBLE  &
NOT POSSIBLE \hspace{.22in} POSSIBLE   & 
NOT POSSIBLE \hspace{.19in} POSSIBLE \\

&
non-$NBB$ & 
\ \   $M_{m,\Phi,\Psi}, M_{m,\Psi,\Phi}$ - not unc.\,conv., & 
see Prop. \ref{propnon}   \hspace{.42in} Example \ref{bnb4}(ii) &
see Prop. \ref{propnon} \hspace{.41in} Example \ref{exadd2} \\

&
& 
\hspace{.36in} see Proposition \ref{lem31}(i) &
 &
 \\

& &  &  & \\
\hline
& &  &  & \\

$\|\!\cdot\!\|$-SN & 
non-$NBA$  \& &
NOT POSSIBLE  \hspace{.17in} NOT POSSIBLE  &
NOT POSSIBLE \hspace{.21in} POSSIBLE   & 
NOT POSSIBLE  \hspace{.17in} NOT POSSIBLE  \\

&
$NBB$ & 
\ \   $M_{m,\Phi,\Psi}, M_{m,\Psi,\Phi}$ - not unc.\,conv., & 
see Prop. \ref{propnon}    \hspace{.41in} Example \ref{bnb10} &
 \ \ \  $M_{m,\Phi,\Psi}$, $M_{m,\Psi,\Phi}$ - not unc.\,conv., \\

&
& 
\hspace{.36in} see Proposition \ref{lem31}(i) &
  &
\hspace{.36in} see Proposition \ref{lem31}(ii)   \\

& &  &  & \\
\hline
& &  &  & \\

$\|\!\cdot\!\|$-SN & 
non-$NBA$ \& &
NOT POSSIBLE  \hspace{.17in} NOT POSSIBLE  &
NOT POSSIBLE \hspace{.21in} POSSIBLE   & 
NOT POSSIBLE \hspace{.19in} POSSIBLE \\

&
non-$NBB$ & 
\ \   $M_{m,\Phi,\Psi}, M_{m,\Psi,\Phi}$ - not unc.\,conv., & 
see Prop. \ref{propnon}    \hspace{.42in} Example \ref{bnbnew1} &
see Prop. \ref{propnon} \hspace{.41in} Example \ref{bnb12} \\

&
& 
\hspace{.36in} see Proposition \ref{lem31}(i) &
  &
  \\

& &  &  & \\
\hline
& &  &  & \\

non-$NBB$  & 
$\|\!\cdot\!\|$-SN &
POSSIBLE \hspace{.47in} POSSIBLE  &
POSSIBLE \hspace{.53in} POSSIBLE  & 
POSSIBLE \hspace{.49in} POSSIBLE  \\

&
& 
Example \ref{bnb13}(i) \hspace{.2in} Example \ref{bnb14}  & 
Example \ref{bnb15} \hspace{.38in} Example \ref{bnb13}(ii)  &
Example \ref{bnb17}(i) \hspace{.23in} Example \ref{bnb17}(ii) \\

& &  &  & \\
\hline
& &  &  & \\

non-$NBB$  & 
$NBA$ \& &
POSSIBLE \hspace{.477in} POSSIBLE   &
POSSIBLE \hspace{.52in} POSSIBLE  & 
POSSIBLE \hspace{.495in} POSSIBLE  \\

&
non-$NBB$ & 
Example \ref{bnbnew2} \hspace{.27in} Example \ref{bnbnew3}  & 
Example \ref{bnbnew2new} \hspace{.31in} Example \ref{bnbnew4}(i)  &
Example \ref{bnbnew4}(ii) \hspace{.13in} Example \ref{bnbnew5} \\

   & &  &  & \\
\hline
& &  &  & \\

non-$NBB$  & 
non-$NBA$   \& &
POSSIBLE \hspace{.477in} POSSIBLE   &
POSSIBLE \hspace{.53in} POSSIBLE  & 
POSSIBLE \hspace{.495in} POSSIBLE  \\

&
$NBB$ & 
Example \ref{bnb19} \hspace{.27in} Example \ref{bnb20}  & 
Example \ref{bnb21i} \hspace{.32in} Example \ref{bnb21ii}  &
Example \ref{bnb23i} \hspace{.28in} Example \ref{bnb23ii} \\

     & &  &  & \\
\hline
& &  &  & \\

non-$NBB$  & 
non-$NBA$ \& &
POSSIBLE \hspace{.477in} POSSIBLE   &
POSSIBLE \hspace{.53in} POSSIBLE  & 
POSSIBLE \hspace{.495in} POSSIBLE  \\

&
non-$NBB$ & 
Example \ref{bnbnew6} \hspace{.27in} Example \ref{bnbnew7}  & 
Example \ref{bnbnew8}(i) \hspace{.2in} Example \ref{bnbnew8}(ii)  &
Example \ref{bnbnew10}(i) \hspace{.16in} Example \ref{bnbnew10}(ii) \\
   
& & & & \\

\hline
\end{longtable}
}

\newpage
{ \small
\centerline{Table 3: two Bessel sequences which are not frames}
\vspace{.1in}

\begin{tabular}{|l|l|l|l|l|}
\hline   & & & & \\
\label{table3}

\!$\phi$\! -\! B. \!\!\!\!& 
\!$\psi$\! -\! B. \!\!\!\!& 
$m$ - SN \hspace{.1in}&
$m\in\ell^\infty$, but non-$SN$ \hspace{.1in} &
$m\notin\ell^\infty $   \\

not fr. \!\!\!\!& 
not fr. \!\!\!\!& 
&
&
\\

\hline   & & & &   \\

 &
 &
$M_{m,\Phi,\Psi}$,\,$M_{m,\Psi,\Phi}$ \hspace{.15in} $M_{m,\Phi,\Psi}$,\,$M_{m,\Psi,\Phi}$\! &
 $M_{m,\Phi,\Psi}$,\,$M_{m,\Psi,\Phi}$ \hspace{.1in} $M_{m,\Phi,\Psi}$,\,$M_{m,\Psi,\Phi}$\! & 
$M_{m,\Phi,\Psi}$,\,$M_{m,\Psi,\Phi}$ \hspace{.09in} $M_{m,\Phi,\Psi}$,\,$M_{m,\Psi,\Phi}$ \! \\

 &
 &
unc.\,conv.\,\& \hspace{.47in} unc.\,conv.\,\& &
unc.\,conv.\,\&\hspace{.45in} unc.\,conv.\,\& & 
unc.\,conv.\,\& \hspace{.4in} unc.\,conv.\,\& \\

 &
 &
INV. \hspace{.845in} NON-INV. &
INV. \hspace{.79in} NON-INV. & 
INV. \hspace{.78in} NON-INV. \\

& &  &  & \\
\hline
& &  &  & \\

$\|\!\cdot\!\|$-SN & 
$\|\!\cdot\!\|$-SN & 
NOT POSSIBLE \hspace{.21in} ALWAYS  &
NOT POSSIBLE \hspace{.18in} ALWAYS  & 
NOT POSSIBLE  \hspace{.16in} NOT POSSIBLE  \\

&
& 
see Prop. \ref{c1} \hspace{.43in} apply Prop. \ref{bmh}, \ref{c1} &
see Prop. \ref{c1} \hspace{.39in} apply Prop. \ref{bmh}, \ref{c1} &
\ \ \ $M_{m,\Phi,\Psi}$, $M_{m,\Psi,\Phi}$ - not unc.\,conv.,\\

&
& 
\hspace{1.1in} Example \ref{bb12}(i) &
\hspace{1.07in} Example \ref{bb12}(ii)  &
\hspace{.36in} see Proposition \ref{lem31}(ii) \\

& &  &  & \\
\hline
& &  &  & \\

$\|\!\cdot\!\|$-SN & 
non-$NBB$ &
NOT POSSIBLE \hspace{.19in} ALWAYS  &
NOT POSSIBLE \hspace{.18in} ALWAYS  & 
NOT POSSIBLE \hspace{.16in} POSSIBLE \\

&
& 
see Prop. \ref{c1} \hspace{.4in} apply Prop. \ref{bmh}, \ref{c1}  & 
see Prop. \ref{c1} \hspace{.38in} apply Prop. \ref{bmh}, \ref{c1}  &
see Prop. \ref{propnon}  \hspace{.37in} Example \ref{bb34}(iii) \\

&
& 
\hspace{1.08in} Example \ref{bb34}(i) &
\hspace{1.06in} Example \ref{bb34}(ii)  &
 \\

& &  &  & \\
\hline
& &  &  & \\

non-$NBB$  & 
non-$NBB$ &
NOT POSSIBLE \hspace{.17in} ALWAYS  &
NOT POSSIBLE \hspace{.16in} ALWAYS  & 
POSSIBLE \hspace{.45in} POSSIBLE \\

&
& 
see Prop. \ref{c1} \hspace{.38in} apply Prop. \ref{bmh}, \ref{c1}  & 
see Prop. \ref{c1} \hspace{.37in} apply Prop. \ref{bmh}, \ref{c1}  &
Example \ref{bb7}(iii)  \hspace{.12in} Example \ref{bb7}(iv) \\

&
& 
\hspace{1.05in} Example \ref{bb7}(i) &
\hspace{1.05in} Example \ref{bb7}(ii)  &
 \\ 
   
& & & & \\

\hline
\end{tabular}
}

\newpage
{\small
\centerline{Table 4: one overcomplete frame, one non-Bessel sequence}

\vspace{.1in}

\begin{longtable}{|l|l|l|l|l|} 
\hline   & & & & \\
\label{table4}
\!$\phi$\! -\! fr. \!\!\!\!& 
\!$\psi$\! -\! not B. \!\!\!\!& 
$m$ - SN \hspace{.1in}&
$m\in\ell^\infty$, but non-$SN$ \hspace{.1in} &
$m\notin\ell^\infty $   \\

not R.b. \!\!\!\!& 
& 
&
&
\\
\hline   & & & & \\
 &
 &
$M_{m,\Phi,\Psi}$,\,$M_{m,\Psi,\Phi}$ \hspace{.21in} $M_{m,\Phi,\Psi}$,\,$M_{m,\Psi,\Phi}$\! &
 $M_{m,\Phi,\Psi}$,\,$M_{m,\Psi,\Phi}$ \hspace{.04in} $M_{m,\Phi,\Psi}$,\,$M_{m,\Psi,\Phi}$\! & 
$M_{m,\Phi,\Psi}$,\,$M_{m,\Psi,\Phi}$ \hspace{.12in} $M_{m,\Phi,\Psi}$,\,$M_{m,\Psi,\Phi}$ \! \\

 &
 &
unc.\,conv.  \& \hspace{.505in} unc.\,conv. \& &
unc.\,conv.  \& \hspace{.335in} unc.\,conv. \& & 
unc.\,conv.  \& \hspace{.41in} unc.\,conv. \& \\

 &
 &
INV. \hspace{.91in} NON-INV. &
INV. \hspace{.73in} NON-INV. & 
INV. \hspace{.815in} NON-INV. \\

& &  &  & \\
\endfirsthead
\hline   & & & & \\
{\footnotesize continued from}  & & & & \\
 {\footnotesize the previous page}  & & & & \\
\hline   & & & & \\

\!$\phi$\! -\! fr. \!\!\!\!& 
\!$\psi$\! -\! not B. \!\!\!\!& 
$m$ - SN \hspace{.1in}&
$m\in\ell^\infty$, but non-$SN$ \hspace{.1in} &
$m\notin\ell^\infty $   \\

not R.b. \!\!\!\!& 
& 
&
&
\\

\hline   & & & &  \\

 &
 &
$M_{m,\Phi,\Psi}$,\,$M_{m,\Psi,\Phi}$ \hspace{.21in} $M_{m,\Phi,\Psi}$,\,$M_{m,\Psi,\Phi}$\! &
 $M_{m,\Phi,\Psi}$,\,$M_{m,\Psi,\Phi}$ \hspace{.04in} $M_{m,\Phi,\Psi}$,\,$M_{m,\Psi,\Phi}$\! & 
$M_{m,\Phi,\Psi}$,\,$M_{m,\Psi,\Phi}$ \hspace{.12in} $M_{m,\Phi,\Psi}$,\,$M_{m,\Psi,\Phi}$ \! \\

 &
 &
unc.\,conv.  \& \hspace{.505in} unc.\,conv. \& &
unc.\,conv.  \& \hspace{.335in} unc.\,conv. \& & 
unc.\,conv.  \& \hspace{.41in} unc.\,conv. \& \\

 &
 &
INV. \hspace{.91in} NON-INV. &
INV. \hspace{.73in} NON-INV. & 
INV. \hspace{.815in} NON-INV. \\

& &  &  & \\
\hline   & & & & \\
\endhead

& & & & {\footnotesize continued  on the next page} \\
\hline
\endfoot

\hline
\endlastfoot

\hline
& &  &  & \\

$\|\!\cdot\!\|$-SN & 
$\|\!\cdot\!\|$-SN & 
NOT POSSIBLE  \hspace{.26in} NOT POSSIBLE &
POSSIBLE  \hspace{.37in} POSSIBLE  & 
NOT POSSIBLE  \hspace{.17in} NOT POSSIBLE  \\

&
& 
\ \ $M_{m,\Phi,\Psi}, M_{m,\Psi,\Phi}$ - not unc.\,conv., &
Example \ref{ex11table4}  \hspace{.22in} Example \ref{fnb2}  &
  $M_{m,\Phi,\Psi}$, $M_{m,\Psi,\Phi}$ - not unc.\,conv.,\\

&
& 
\hspace{.36in} see Proposition \ref{lem31}(i) &
  &
\hspace{.36in} see Proposition \ref{lem31}(ii) \\

& &  &  & \\
\hline
& &  &  & \\

$\|\!\cdot\!\|$-SN & 
$NBA$ $\&$ &
NOT POSSIBLE  \hspace{.26in} NOT POSSIBLE  &
POSSIBLE \hspace{.37in} POSSIBLE   & 
POSSIBLE \hspace{.46in} POSSIBLE \\

&
 non-$NBB$  & 
\ \ $M_{m,\Phi,\Psi}, M_{m,\Psi,\Phi}$ - not unc.\,conv.,  & 
Example \ref{fnb51}  \hspace{.22in} Example \ref{fnb52}(i) &
Example \ref{fnb5b} \hspace{.31in} Example \ref{fnb52}(ii) \\

&
& 
\hspace{.36in} see Proposition \ref{lem31}(i) &
 &
  \\

& &  &  & \\
\hline
& &  &  & \\

$\|\!\cdot\!\|$-SN & 
non-$NBA$  \& &
NOT POSSIBLE  \hspace{.26in} NOT POSSIBLE  &
POSSIBLE \hspace{.38in} POSSIBLE   & 
NOT POSSIBLE \hspace{.17in} NOT POSSIBLE \\

&
 $NBB$ & 
\ \ $M_{m,\Phi,\Psi}, M_{m,\Psi,\Phi}$ - not unc.\,conv.,    & 
Example \ref{fnb3}(i)  \hspace{.11in} Example \ref{fnb3}(ii) &
\ \ $M_{m,\Phi,\Psi}, M_{m,\Psi,\Phi}$ - not unc.\,conv.,  \\

&
& 
\hspace{.36in} see Proposition \ref{lem31}(i) &
&
\hspace{.36in} see Proposition \ref{lem31}(ii) \\

& &  &  & \\
\hline
& &  &  & \\

$\|\!\cdot\!\|$-SN & 
non-$NBA$ \& &
NOT POSSIBLE  \hspace{.26in} NOT POSSIBLE  &
POSSIBLE \hspace{.38in} POSSIBLE   & 
POSSIBLE \hspace{.465in} POSSIBLE \\

&
non-$NBB$  & 
\ \ $M_{m,\Phi,\Psi}, M_{m,\Psi,\Phi}$ - not unc.\,conv.,    & 
Example \ref{fnb31}  \hspace{.23in} Example \ref{fnb32} &
Example \ref{fnb5} \hspace{.32in} Example \ref{fnb6} \\

&
& 
\hspace{.36in} see Proposition \ref{lem31}(i) &
&
\\

& &  &  & \\
\hline
& &  &  & \\

non-$NBB$  & 
$\|\!\cdot\!\|$-SN &
POSSIBLE \hspace{.555in} POSSIBLE  &
POSSIBLE \hspace{.36in} POSSIBLE  & 
POSSIBLE \hspace{.48in} POSSIBLE  \\

&
& 
Example \ref{fnb7}(i) \hspace{.21in} Example \ref{fnb8}  & 
Example \ref{fnb9} \hspace{.15in} Example \ref{fnb7}(ii)  &
Example \ref{fnb11}(i) \hspace{.15in} Example \ref{fnb11}(ii) \\

& &  &  & \\
\hline
& &  &  & \\

non-$NBB$  & 
$NBA$ \& &
POSSIBLE \hspace{.57in} POSSIBLE   &
POSSIBLE \hspace{.35in} POSSIBLE  & 
POSSIBLE \hspace{.48in} POSSIBLE  \\

&
non-$NBB$ & 
Example \ref{fnb13} \hspace{.36in} Example \ref{fnb14}(i)  & 
Example \ref{fnb33} \hspace{.14in} Example \ref{fnb16}  &
Example \ref{fnb14}(ii) \hspace{.12in} Example \ref{fnb18} \\

& &  &  & \\
\hline
& &  &  & \\

non-$NBB$  & 
non-$NBA$  \& &
POSSIBLE \hspace{.57in} POSSIBLE   &
POSSIBLE \hspace{.37in} POSSIBLE  & 
POSSIBLE \hspace{.5in} POSSIBLE  \\

&
$NBB$ & 
Example \ref{fnb34} \hspace{.36in} Example \ref{fnb35}  & 
Example \ref{fnb15}(i) \hspace{.04in} Example \ref{fnb15}(ii)  &
Example \ref{fnb36}(i) \hspace{.17in} Example \ref{fnb36}(ii) \\

& &  &  & \\
\hline
& &  &  & \\

non-$NBB$  & 
non-$NBA$ \& &
POSSIBLE \hspace{.57in} POSSIBLE   &
POSSIBLE \hspace{.37in} POSSIBLE  & 
POSSIBLE \hspace{.5in} POSSIBLE  \\

&
non-$NBB$ & 
Example \ref{fnb37} \hspace{.36in} Example \ref{fnb38}  & 
Example \ref{fnb39}(i) \hspace{.04in}   Example \ref{fnb39}(ii)  &
Example \ref{fnb40} \hspace{.29in} Example \ref{fnb41} \\
   
& & & & \\
   
\end{longtable}

} 

\vspace{.1in}

\newpage
{\small 
\centerline{Table 5: one overcomplete frame, one Bessel sequence which is not a frame}

\vspace{.1in}

\begin{tabular}{|l|l|l|l|l|}
\hline   & & & & \\
\label{table5}

\!$\phi$\! -\! fr. \!\!\!\!& 
\!$\psi$\! -\! B. \!\!\!\!& 
$m$ - SN \hspace{.1in}&
$m\in\ell^\infty$, but non-$SN$ \hspace{.1in} &
$m\notin\ell^\infty $   \\

not R.b. \!\!\!\!& 
not fr. \!\!\!\!& 
&
 &
   \\

\hline   & & & &   \\

 &
 &
$M_{m,\Phi,\Psi}$,\,$M_{m,\Psi,\Phi}$ \hspace{.135in} $M_{m,\Phi,\Psi}$,\,$M_{m,\Psi,\Phi}$\! &
 $M_{m,\Phi,\Psi}$,\,$M_{m,\Psi,\Phi}$ \hspace{.15in} $M_{m,\Phi,\Psi}$,\,$M_{m,\Psi,\Phi}$\! & 
$M_{m,\Phi,\Psi}$,\,$M_{m,\Psi,\Phi}$ \hspace{.09in} $M_{m,\Phi,\Psi}$,\,$M_{m,\Psi,\Phi}$ \! \\

 &
 &
unc.\,conv.\,\& \hspace{.43in} unc.\,conv.\,\& &
unc.\,conv.\,\&\hspace{.495in} unc.\,conv.\,\& & 
unc.\,conv.\,\& \hspace{.405in} unc.\,conv.\,\& \\

 &
 &
INV. \hspace{.81in} NON-INV. &
INV. \hspace{.84in} NON-INV. & 
INV. \hspace{.79in} NON-INV. \\

& &  &  & \\
\hline
& &  &  & \\

$\|\!\cdot\!\|$-SN & 
$\|\!\cdot\!\|$-SN & 
NOT POSSIBLE \hspace{.16in} ALWAYS  &
NOT POSSIBLE \hspace{.2in} ALWAYS  & 
NOT POSSIBLE  \hspace{.14in} NOT POSSIBLE  \\

&
& 
see Prop. \ref{c1} \hspace{.37in} apply Prop. \ref{bmh}, \ref{c1} &
see Prop. \ref{c1} \hspace{.41in} apply Prop. \ref{bmh}, \ref{c1} &
  $M_{m,\Phi,\Psi}$, $M_{m,\Psi,\Phi}$ - not unc.\,conv.,\\

&
& 
\hspace{1.05in} Example \ref{fb1}(i) &
\hspace{1.09in} Example \ref{fb1}(ii)  &
\hspace{.36in} see Proposition \ref{lem31}(ii) \\

& &  &  & \\
\hline
& &  &  & \\

$\|\!\cdot\!\|$-SN & 
non-$NBB$ &
NOT POSSIBLE \hspace{.16in} ALWAYS  &
NOT POSSIBLE \hspace{.2in} ALWAYS  & 
POSSIBLE \hspace{.45in} POSSIBLE \\

&
& 
see Prop. \ref{c1} \hspace{.36in} apply Prop. \ref{bmh}, \ref{c1}  & 
see Prop. \ref{c1} \hspace{.41in} apply Prop. \ref{bmh}, \ref{c1}  &
Example \ref{fb3}(iii) \hspace{.12in} Example \ref{fb6} \\

&
& 
\hspace{1.05in} Example \ref{fb3}(i) &
\hspace{1.1in} Example \ref{fb3}(ii)  &
 \\ 

& &  &  & \\
\hline
& &  &  & \\

non-$NBB$  & 
$\|\!\cdot\!\|$-SN &
NOT POSSIBLE \hspace{.16in} ALWAYS  &
NOT POSSIBLE \hspace{.2in} ALWAYS  & 
NOT POSSIBLE \hspace{.17in} POSSIBLE \\

&
& 
see Prop. \ref{c1} \hspace{.38in} apply Prop. \ref{bmh}, \ref{c1}  & 
see Prop. \ref{c1} \hspace{.4in} apply Prop. \ref{bmh}, \ref{c1}  &
see Prop. \ref{propnon} \hspace{.38in} Example \ref{fb10} \\

&
& 
\hspace{1.06in} Example \ref{fb7}(i) &
\hspace{1.08in} Example \ref{fb7}(ii)  &
  \\

& &  &  & \\
\hline
& &  &  & \\

non-$NBB$  & 
non-$NBB$ &
NOT POSSIBLE \hspace{.16in} ALWAYS  &
NOT POSSIBLE \hspace{.2in} ALWAYS  & 
POSSIBLE \hspace{.45in} POSSIBLE \\

&
& 
see Prop. \ref{c1} \hspace{.36in} apply Prop. \ref{bmh}, \ref{c1}  & 
see Prop. \ref{c1} \hspace{.4in} apply Prop. \ref{bmh}, \ref{c1}  &
Example \ref{fb13}(i)  \hspace{.19in} Example \ref{fb13}(ii) \\

&
& 
\hspace{1.04in} Example \ref{fb11}(i) &
\hspace{1.08in} Example \ref{fb11}(ii)  &
 \\ 
   
& & & & \\

\hline
\end{tabular}
}

\vspace{.1in}

\newpage
{ \small 
\centerline{Table 6: two overcomplete frames}

\vspace{.1in}
\begin{tabular}{|l|l|l|l|l|}
\hline   & & & & \\
\label{table6}

\!$\phi$\! -\! fr. \!\!\!\!& 
\!$\psi$\! -\! fr. \!\!\!\!& 
$m$ - SN \hspace{.1in}&
$m\in\ell^\infty$, but non-$SN$ \hspace{.1in} &
$m\notin\ell^\infty $   \\

not R.b. \!\!\!\!& 
not R.b \!\!\!\!& 
&
 &
   \\

\hline   & & & &   \\

 &
 &
$M_{m,\Phi,\Psi}$,\,$M_{m,\Psi,\Phi}$ \hspace{.11in} $M_{m,\Phi,\Psi}$,\,$M_{m,\Psi,\Phi}$\! &
 $M_{m,\Phi,\Psi}$,\,$M_{m,\Psi,\Phi}$ \hspace{.03in} $M_{m,\Phi,\Psi}$,\,$M_{m,\Psi,\Phi}$\! & 
$M_{m,\Phi,\Psi}$,\,$M_{m,\Psi,\Phi}$ \hspace{.07in} $M_{m,\Phi,\Psi}$,\,$M_{m,\Psi,\Phi}$ \! \\

 &
 &
unc.\,conv.\,\& \hspace{.41in} unc.\,conv.\,\& &
unc.\,conv.\,\&\hspace{.37in} unc.\,conv.\,\& & 
unc.\,conv.\,\& \hspace{.38in} unc.\,conv.\,\& \\

 &
 &
INV. \hspace{.79in} NON-INV. &
INV. \hspace{.71in} NON-INV. & 
INV. \hspace{.76in} NON-INV. \\

& &  &  & \\
\hline
& &  &  & \\

$\|\!\cdot\!\|$-SN & 
$\|\!\cdot\!\|$-SN & 
POSSIBLE \hspace{.43in} POSSIBLE  &
POSSIBLE \hspace{.38in} POSSIBLE  & 
NOT POSSIBLE  \hspace{.13in} NOT POSSIBLE  \\

&
& 
Example \ref{ff1}(i) \hspace{.16in} Example \ref{ff2} &
Example \ref{ff1}(ii) \hspace{.08in} Example \ref{ff1}(iii)    &
\, $M_{m,\Phi,\Psi}$, $M_{m,\Psi,\Phi}$ - not unc. conv.,\\

&
& 
 &
 &
\hspace{.36in} see Proposition \ref{lem31}(ii) \\

& &  &  & \\
\hline
& &  &  & \\

$\|\!\cdot\!\|$-SN & 
non-$NBB$ &
POSSIBLE \hspace{.43in} POSSIBLE &
POSSIBLE  \hspace{.38in} POSSIBLE  &
POSSIBLE \hspace{.44in} POSSIBLE \\

&
& 
Example \ref{ff9}(i) \hspace{.17in} Example \ref{ff6}   & 
Example \ref{ff9}(ii)  \hspace{.08in} Example \ref{ff9}(iii)  &
Example \ref{ff9}(iv) \hspace{.12in} Example \ref{ff9}(v) \\

& &  &  & \\
\hline
& &  &  & \\

non-$NBB$  & 
non-$NBB$ &
POSSIBLE \hspace{.42in} POSSIBLE &
POSSIBLE  \hspace{.4in} POSSIBLE  &
POSSIBLE \hspace{.43in} POSSIBLE \\

&
& 
Example \ref{ff11} \hspace{.28in} Example \ref{ff12}   & 
Example \ref{ff13}(i)  \hspace{.13in} Example \ref{ff13}(ii)  &
Example \ref{ff15}(i) \hspace{.17in} Example \ref{ff15}(ii) \\

& & & & \\

\hline
\end{tabular}
}

\vspace{.1in}

\newpage

{\small 
\centerline{Table 7: one Riesz basis, one non-Bessel sequence}

\vspace{.1in}
\begin{tabular}{|l|l|l|l|l|}
\hline   & & & & \\
\label{table7}

\!$\phi$\! -\! R.b. \!\!\!\!& 
\!$\psi$\! -\! not B. \!\!\!\!& 
$m$ - SN \hspace{.1in}&
$m\in\ell^\infty$, but non-$SN$ \hspace{.1in} &
$m\notin\ell^\infty $   \\

\hline   & & & &   \\

 &
 &
$M_{m,\Phi,\Psi}$,\,$M_{m,\Psi,\Phi}$ \hspace{.16in} $M_{m,\Phi,\Psi}$,\,$M_{m,\Psi,\Phi}$\! &
 $M_{m,\Phi,\Psi}$,\,$M_{m,\Psi,\Phi}$ \hspace{.1in} $M_{m,\Phi,\Psi}$,\,$M_{m,\Psi,\Phi}$\! & 
$M_{m,\Phi,\Psi}$,\,$M_{m,\Psi,\Phi}$ \hspace{.06in} $M_{m,\Phi,\Psi}$,\,$M_{m,\Psi,\Phi}$ \! \\

 &
 &
unc.\,conv.\,\& \hspace{.47in} unc.\,conv.\,\& &
unc.\,conv.\,\&\hspace{.45in} unc.\,conv.\,\& & 
unc.\,conv.\,\& \hspace{.38in} unc.\,conv.\,\& \\

 &
 &
INV. \hspace{.85in} NON-INV. &
INV. \hspace{.79in} NON-INV. & 
INV. \hspace{.76in} NON-INV. \\

& &  &  & \\
\hline
& &  &  & \\

$\|\!\cdot\!\|$-SN & 
$\|\!\cdot\!\|$-SN & 
NOT POSSIBLE \hspace{.19in} NOT POSSIBLE  &
NOT POSSIBLE \hspace{.16in} POSSIBLE  & 
NOT POSSIBLE  \hspace{.14in} NOT POSSIBLE  \\

&
& 
\ \ \ $M_{m,\Phi,\Psi}$, $M_{m,\Psi,\Phi}$ - not well defined, &
see Prop. \ref{rc} \hspace{.38in} Example \ref{rnb1}  &
\ \ $M_{m,\Phi,\Psi}$, $M_{m,\Psi,\Phi}$ - not well defined,\\

&
& 
\hspace{.45in} see Proposition \ref{lem36}(i) 
 &
 &
\hspace{.36in} see  Proposition \ref{lem36}(ii) \\

& &  &  & \\
\hline
& &  &  & \\

$\|\!\cdot\!\|$-SN & 
$NBA$ \& & 
NOT POSSIBLE \hspace{.19in} NOT POSSIBLE  &
NOT POSSIBLE \hspace{.16in} POSSIBLE  & 
NOT POSSIBLE  \hspace{.15in} POSSIBLE  \\

&
non-$NBB$ & 
\ \ \ $M_{m,\Phi,\Psi}$, $M_{m,\Psi,\Phi}$ - not well defined, &
see Prop. \ref{rc} \hspace{.37in} Example \ref{rnb31}  &
see Prop. \ref{rc} \hspace{.35in} Example \ref{rnb32}\\

&
& 
\hspace{.45in} see Proposition \ref{lem36}(i) 
 &
 &
 \\

& &  &  & \\
\hline
& &  &  & \\

$\|\!\cdot\!\|$-SN & 
non-$NBA$ \&  &
NOT POSSIBLE \hspace{.2in} NOT POSSIBLE &
POSSIBLE  \hspace{.48in} POSSIBLE  &
NOT POSSIBLE \hspace{.16in} NOT POSSIBLE  \\

&
$NBB$  & 
\ \ \ $M_{m,\Phi,\Psi}$, $M_{m,\Psi,\Phi}$ - not well defined,   & 
Example \ref{rnb2}  \hspace{.33in} Example \ref{rnb3}  &
\ \ \ $M_{m,\Phi,\Psi}$, $M_{m,\Psi,\Phi}$ - not well defined, \\

&  
& 
\hspace{.45in} see Proposition \ref{lem36}(i) &
  &
\hspace{.45in} see Proposition \ref{lem36}(ii) \\
   
& & & & \\

& &  &  & \\
\hline
& &  &  & \\

$\|\!\cdot\!\|$-SN & 
non-$NBA$ \& &
NOT POSSIBLE \hspace{.2in} NOT POSSIBLE &
NOT POSSIBLE  \hspace{.16in} POSSIBLE  &
POSSIBLE \hspace{.44in} POSSIBLE \\

&
non-$NBB$ & 
\ \ \ $M_{m,\Phi,\Psi}$, $M_{m,\Psi,\Phi}$ - not well defined,   & 
see Prop. \ref{rc} \hspace{.37in} Example \ref{rnb33}  &
Example \ref{rnb4} \hspace{.29in} Example \ref{rnb5} \\

&  
& 
\hspace{.45in} see Proposition \ref{lem36}(i) &
  &
 \\
   
& & & & \\
   
\hline
\end{tabular}
}

\vspace{.1in}

\newpage
{ \small
\centerline{Table 8: one Riesz basis, one Bessel which is not a frame}

\vspace{.1in}
\begin{tabular}{|l|l|l|l|l|}
\hline   & & & & \\
\label{table8}

\!$\phi$\! -\! R.b. \!\!\!\!& 
\!$\psi$\! -\! B.  \!\!\!\!& 
$m$ - SN \hspace{.1in}&
$m\in\ell^\infty$, but non-$SN$ \hspace{.1in} &
$m\notin\ell^\infty $   \\

& 
not fr.& 
&
&
 \\

\hline   & & & &   \\

 &
 &
$M_{m,\Phi,\Psi}$,\,$M_{m,\Psi,\Phi}$ \hspace{.15in} $M_{m,\Phi,\Psi}$,\,$M_{m,\Psi,\Phi}$\! &
 $M_{m,\Phi,\Psi}$,\,$M_{m,\Psi,\Phi}$ \hspace{.09in} $M_{m,\Phi,\Psi}$,\,$M_{m,\Psi,\Phi}$\! & 
$M_{m,\Phi,\Psi}$,\,$M_{m,\Psi,\Phi}$ \hspace{.1in} $M_{m,\Phi,\Psi}$,\,$M_{m,\Psi,\Phi}$ \! \\

 &
 &
unc.\,conv.\,\& \hspace{.46in} unc.\,conv.\,\& &
unc.\,conv.\,\&\hspace{.45in} unc.\,conv.\,\& & 
unc.\,conv.\,\& \hspace{.42in} unc.\,conv.\,\& \\

 &
 &
INV. \hspace{.83in} NON-INV. &
INV. \hspace{.78in} NON-INV. & 
INV. \hspace{.79in} NON-INV. \\

& &  &  & \\
\hline
& &  &  & \\

$\|\!\cdot\!\|$-SN & 
$\|\!\cdot\!\|$-SN & 
NOT POSSIBLE \hspace{.19in} ALWAYS  &
NOT POSSIBLE \hspace{.13in} ALWAYS  & 
NOT POSSIBLE  \hspace{.16in} NOT POSSIBLE  \\

&
& 
see Prop. \ref{rc} \hspace{.39in} apply Prop. \ref{bmh}, \ref{rc} &
see Prop. \ref{rc} \hspace{.33in} apply Prop. \ref{bmh}, \ref{rc}   &
\ \ $M_{m,\Phi,\Psi}$, $M_{m,\Psi,\Phi}$ - not well defined,\\

&
& 
\hspace{1.07in} Example \ref{rb1} &
\hspace{1.02in} Example \ref{rb1b} &
\hspace{.36in} see Proposition \ref{lem36}(ii) \\

& &  &  & \\
\hline
& &  &  & \\

$\|\!\cdot\!\|$-SN & 
non-$NBB$ &
NOT POSSIBLE \hspace{.17in} ALWAYS &
NOT POSSIBLE  \hspace{.116in} ALWAYS  &
POSSIBLE \hspace{.46in} POSSIBLE \\

&
& 
see Prop. \ref{rc}  \hspace{.38in} apply Prop. \ref{bmh}, \ref{rc}  & 
see Prop. \ref{rc} \hspace{.32in} apply Prop. \ref{bmh}, \ref{rc} &
Example \ref{rb4} \hspace{.31in} Example \ref{rb5} \\

&  
& 
 \hspace{1.06in} Example \ref{rb2} &
\hspace{1.01in} Example \ref{rb3}  &
 \\
   
& & & & \\
   
\hline
\end{tabular}
}

\vspace{.1in}

\newpage
{ \small
\centerline{Table 9: one Riesz basis, one overcomplete frame}

\vspace{.1in}
\begin{tabular}{|l|l|l|l|l|}
\hline   & & & & \\
\label{table9}

\!$\phi$\! -\! R.b. \!\!\!\!& 
\!$\psi$\! -\! fr.  \!\!\!\!& 
$m$ - SN \hspace{.1in}&
$m\in\ell^\infty$, but non-$SN$ \hspace{.1in} &
$m\notin\ell^\infty $   \\

& 
not R.b.& 
&
&
 \\

\hline   & & & &   \\

 &
 &
$M_{m,\Phi,\Psi}$,\,$M_{m,\Psi,\Phi}$ \hspace{.24in} $M_{m,\Phi,\Psi}$,\,$M_{m,\Psi,\Phi}$\! &
 $M_{m,\Phi,\Psi}$,\,$M_{m,\Psi,\Phi}$ \hspace{.1in} $M_{m,\Phi,\Psi}$,\,$M_{m,\Psi,\Phi}$\! & 
$M_{m,\Phi,\Psi}$,\,$M_{m,\Psi,\Phi}$ \hspace{.09in} $M_{m,\Phi,\Psi}$,\,$M_{m,\Psi,\Phi}$ \! \\

 &
 &
unc.\,conv.\,\& \hspace{.56in} unc.\,conv.\,\& &
unc.\,conv.\,\&\hspace{.46in} unc.\,conv.\,\& & 
unc.\,conv.\,\& \hspace{.41in} unc.\,conv.\,\& \\

 &
 &
INV. \hspace{.93in} NON-INV. &
INV. \hspace{.78in} NON-INV. & 
INV. \hspace{.78in} NON-INV. \\

& &  &  & \\
\hline
& &  &  & \\

$\|\!\cdot\!\|$-SN & 
$\|\!\cdot\!\|$-SN & 
NOT POSSIBLE \hspace{.3in} ALWAYS  &
NOT POSSIBLE \hspace{.16in} ALWAYS  & 
NOT POSSIBLE  \hspace{.15in} NOT POSSIBLE  \\

&
& 
see Prop. \ref{rbis} \hspace{.52in} apply Prop. \ref{bmh}, \ref{rbis} &
see Prop. \ref{rc} \hspace{.37in} apply Prop. \ref{bmh}, \ref{rc}   &
\ \ $M_{m,\Phi,\Psi}$, $M_{m,\Psi,\Phi}$ - not well defined,\\

&
& 
$M_{m,\Phi,\Psi}$-inj., non-surj. \hspace{.001in} Example \ref{rf1} &
\hspace{1.06in} Example \ref{rf2} &
\hspace{.43in} see Proposition \ref{lem36}(ii) \\ 
   
& 
&  
$M_{m,\Psi,\Phi}$-surj., non-inj. &
 &
\\
 
& &  &  & \\
\hline
& &  &  & \\

$\|\!\cdot\!\|$-SN & 
non-$NBB$ &
NOT POSSIBLE \hspace{.29in} ALWAYS &
NOT POSSIBLE  \hspace{.14in} ALWAYS  &
NOT POSSIBLE \hspace{.14in} POSSIBLE \\

&
& 
see Prop. \ref{rbis}  \hspace{.51in} apply Prop. \ref{bmh}, \ref{rbis}  & 
see Prop. \ref{rc} \hspace{.35in} apply Prop. \ref{bmh}, \ref{rc} &
see Prop. \ref{rc} \hspace{.35in} Example \ref{rf5} \\

&  
& 
$M_{m,\Phi,\Psi}$-inj., non-surj.  \hspace{-.001in} Example \ref{rf3} &
\hspace{1.03in} Example \ref{rf4}  &
 \\
   
& 
& 
$M_{m,\Psi,\Phi}$-surj., non-inj. & 
& \\
   
   & & & & \\
   
\hline
\end{tabular}

}

\vspace{.1in}

{ \small

\centerline{Table 10: two Riesz bases}

\vspace{.1in}
\begin{tabular}{|l|l|l|l|l|}
\hline   & & & & \\
\label{table10}

\!$\phi$\! -\! R.b. \!\!\!\!& 
\!$\psi$\! -\! R.b.  \!\!\!\!& 
$m$ - SN \hspace{.1in}&
$m\in\ell^\infty$, but non-$SN$ \hspace{.1in} &
$m\notin\ell^\infty $   \\
\hline   & & \hspace{1.089in}  & &   \\

 &
 &
$M_{m,\Phi,\Psi}$,\,$M_{m,\Psi,\Phi}$ \hspace{.29in}   $M_{m,\Phi,\Psi}$,\,$M_{m,\Psi,\Phi}$\! &
 $M_{m,\Phi,\Psi}$,\,$M_{m,\Psi,\Phi}$ \hspace{.13in} $M_{m,\Phi,\Psi}$,\,$M_{m,\Psi,\Phi}$\! & 
$M_{m,\Phi,\Psi}$,\,$M_{m,\Psi,\Phi}$ \hspace{.11in} $M_{m,\Phi,\Psi}$,\,$M_{m,\Psi,\Phi}$ \! \\

 &
 &
unc.\,conv.\,\& \hspace{.6in} unc.\,conv.\,\& &
unc.\,conv.\,\&\hspace{.48in} unc.\,conv.\,\& & 
unc.\,conv.\,\& \hspace{.44in} unc.\,conv.\,\& \\

 &
 &
INV. \hspace{.97in} NON-INV. &
INV. \hspace{.82in} NON-INV. & 
INV. \hspace{.81in} NON-INV. \\

 & & & & \\
\hline
 & & & & \\
 
$\|\!\cdot\!\|$-SN &
$\|\!\cdot\!\|$-SN \ \ \ \ \  & 
 ALWAYS  \hspace{.73in} NOT POSSIBLE &
NOT POSSIBLE \hspace{.16in} ALWAYS & 
NOT POSSIBLE  \hspace{.16in} NOT POSSIBLE \\

&
&
apply Prop. \ref{bmh}, \ref{rbinv} \hspace{.23in} apply Prop. \ref{bmh}, \ref{rbinv}  &
see Prop. \ref{rc} \hspace{.35in} apply Prop. \ref{bmh}, \ref{rc} &
\ \ $M_{m,\Phi,\Psi}$, $M_{m,\Psi,\Phi}$ - not well defined, \\

 &
 & 
Example \ref{rry} &
\hspace{1.04in} Example \ref{rrn} &
\hspace{.43in} see Proposition \ref{lem36}(ii) \\

& & & & \\
\hline  

\end{tabular}

}

\newpage

\section{Examples}\label{sec:examples1}

\subsection{Examples for two non-Bessel sequences; Table 1 
on page \pageref{table1}}

\begin{rem}\label{r1}
When $\Phi$ and $\Psi$ are \normsn non-Bessel sequences and $m$ is $SN$, then both $M_{m,\Phi,\Psi}$ and $M_{m,\Psi,\Phi}$ can never be unconditionally convergent on $\h$ due to Proposition \ref{lem31}(i). However, they can be conditionally convergent and invertible (resp. conditionally convergent and non-invertible).
Examples:
\begin{itemize}
\item[{\rm (a)}] Let 
$\Phi=(e_1, \ \ e_2,\ e_2,\, \ \ \, e_2, \ \ e_3, \ e_3,\, \ \ \ e_3, \ e_3, \ \ \ \ e_3, \ \ e_4, \ \, e_4, \ \ \ \, e_4, \ \, e_4, \ \ \ \, e_4, \ e_4, \ \ \ \, e_4, \ldots)$,

\hspace{.18in}
\vspace{.05in}
$\Psi=(e_1, \ \ e_2, \ e_2, \, - e_2, \ \ e_3, \ e_3, \, - e_3, \ e_3,\, - e_3,     \  \ e_4, \ e_4, \, -  e_4, \ e_4, \, - e_4, \ \, e_4,  \, - e_4, \ldots)$.

\noindent  Then $M_{(1),\Phi,\Psi}=M_{(1),\Psi,\Phi}=I$.  
\item[{\rm (b)}] Let 
$\Phi=(e_1, \ \ e_2,\ e_2,\, \ \ \ e_2, \ \ e_3, \ e_3,\, \ \ \ e_3, \ e_3, \, \ \ \ e_3, \ \ e_4, \ \, e_4, \ \ \ e_4, \ \, e_4, \, \ \ \, e_4, \ \, e_4, \ \ \ \, e_4, \ldots)$,

\hspace{.18in}
\vspace{.05in}
$\Psi=(e_1, \ \ e_1, \ e_2, \, - e_2, \ \ e_3, \ e_3, \, - e_3, \ e_3,\, - e_3,     \  \ e_4, \ e_4, \, -  e_4, \ e_4, \, - e_4, \ \, e_4,  \, - e_4, \ldots)$.

\noindent Then $M_{(1),\Phi,\Psi}f = f+ \<f,e_1-e_2\>e_2$ and $M_{(1),\Psi,\Phi}f=f+ \<f,e_2\>(e_1-e_2)$, $f\in\h$. Both $M_{(1),\Phi,\Psi}$ and $M_{(1),\Psi,\Phi}$ are not injective - for example, $M_{m,\Phi,\Psi}e_2=0$ and $M_{m,\Psi,\Phi}e_1=e_1=M_{m,\Psi,\Phi}e_2$.
\end{itemize}

\end{rem}

\begin{rem}\label{r11} 
When $\Phi$ is \normsn non-Bessel for $\h$, $\Psi$ is $NBA$ non-$NBB$ non-Bessel for $\h$ and $m$ is $SN$, then both $M_{m,\Phi,\Psi}$ and $M_{m,\Psi,\Phi}$ can not be unconditionally convergent on $\h$ due to Proposition \ref{lem31}(i), but they can be conditionally convergent and invertible (resp. conditionally convergent and non-invertible). Examples:
\begin{itemize}
\item[{\rm (a)}]
Let 
$\Phi=(e_1, \ \ \ e_2, \ \ \ \ \, e_2, \ \ \ \ \ \ \ e_2, \ \ \ \ e_3, \ \ \  e_3, \ \ \ \ e_3, \ \ \ \ \, e_3, \ \ \ \ \ \ \ e_3, \ \ \  \ e_4, \ \ \  e_4, \ \ \ \ e_4, \ \ \   e_4, \ \ \ \ e_4,  \ \ \ \ \, e_4, \ \ \ \ \ \ \ \, e_4,\ \ldots)$,

\hspace{.18in}
\vspace{.05in}
$\Psi=(e_1, \ \ \ e_2, \ \ \frac{1}{2}  e_2, \, \ -\frac{1}{2}e_2,  \ \ \ \ e_3, \ \ \ e_3, \  - e_3,  \ \ \frac{1}{3}  e_3, \, \ -\frac{1}{3}e_3,  \ \ \  \ e_4, \ \ \  e_4, \ -  e_4, \ \ \,  e_4, \  -  e_4, \ \ \frac{1}{4}  e_4, \, \ -\frac{1}{4}e_4,\ \ldots)$.

\noindent Then 
 $M_{(1),\Phi,\Psi}=M_{(1),\Psi,\Phi}=I$. 
\item[{\rm (b)}]
Let 
$\Phi=(e_1, \ \ \ e_2, \ \ \ \ \, e_2, \ \ \ \ \ \ \, e_2, \ \ \ \ \ e_3, \ \ \ e_3, \ \ \ \ e_3, \ \ \ \ \, e_3, \ \ \ \ \ \, e_3, \ \ \  \ e_4, \ \ \  e_4, \ \ \ \ \, e_4, \ \ \  e_4, \ \ \ \, e_4,  \ \ \ \ e_4, \ \ \ \ \ \ \, e_4,\ \ldots)$,

\hspace{.18in}
\vspace{.05in}
$\Psi=(e_1, \ \frac{1}{2} e_2, \ \ \ \  e_2, \, \ \ \ - e_2,  \ \ \, \frac{1}{3} e_3, \ \ \ e_3, \  - e_3,  \ \ \ \  e_3,  \ \ \ - e_3,  \ \ \frac{1}{4} e_4, \ \ \  e_4, \ -  e_4, \ \ \,  e_4, \  -  e_4, \ \ \ \  e_4, \, \ \ \ - e_4,\ \ldots)$.

\noindent Then 
 $M_{(1),\Phi,\Psi}=M_{(1),\Psi,\Phi}=G_1$ - non-invertible on $\h$ (see Lemma \ref{lemg}). 
\end{itemize}
\end{rem}

There are also other places in the table, where the multipliers can not be unconditionally convergent on $\h$, but they can be conditionally convergent an invertible (resp. conditionally convergent and non-invertible). We will not list such examples any more, as our main aim is to consider the possibilities for the combination of unconditional convergence and invertibility on $\h$.

  \begin{ex} \label{nbnb1}
 Let 
$\Phi=(e_1, \ e_2,\  e_2,\ e_3, \ e_3,\  \ e_3,\  e_4,\  e_4, \ e_4, \ e_4, \ \ldots)$,

\vspace{.1in}
\hspace{.51in} (i) Let  $m=(\  1, \ \ \, \frac{1}{2}, \ \ \, \frac{1}{2},\ \ \frac{1}{3},\ \ \frac{1}{3},\,  \ \ \frac{1}{3},\ \ \, \frac{1}{4},\ \ \frac{1}{4},\ \ \frac{1}{4},\ \ \frac{1}{4}, \ \ldots)$.

\vspace{.05in}
 \noindent
 Then $M_{m,\Phi,\Phi}=I$. The convergence is unconditional on $\h$, because $M_{m,\Phi,\Phi}\neweq M_{(1),(\sqrt{m_n}\,\phi_n),(\sqrt{m_n}\,\phi_n)}$  and the sequence   $(\sqrt{m_n}\,\phi_n)$ is Bessel for $\h$ (apply Prop. \ref{bmh}).

 \vspace{.1in}
\hspace{.49in} (ii) Let
 $m=(\ 1, \  \frac{1}{2^2}, \ \frac{1}{2^2},\, \frac{1}{3^2},\ \frac{1}{3^2},\, \frac{1}{3^2},\, \frac{1}{4^2}, \ \frac{1}{4^2}, \, \frac{1}{4^2},\,  \frac{1}{4^2}, \ \ldots)$.
   
\vspace{.05in} 
\noindent
Then $M_{m,\Phi,\Phi}=G_1$ - non-invertible on $\h$ (see Lemma \ref{lemg}). The unconditional convergence on $\h$ follows as in (i).
\end{ex}

\begin{ex}\label{nbnb9} Let 
  $\Phi=(\,e_1, \ \ \ e_2, \ \ \ \ \ \, e_2, \ \ \ \ \ \ \, \ e_2, \ \ \ \ \ \ \, e_1, \ \ \ e_3,\ \ \ \ \ \, e_3, \ \ \ \ \ \ \ \ \, e_3, \ \ \ \ \ \ \, e_1, \ \ \  e_4, \ \ \ \ \ \, e_4, \ \ \ \ \ \ \ \ e_4,\ldots)$,
  
  \vspace{.05in} \hspace{.91in}
 $\Psi=(\, e_1,\ \ \ e_2, \ \ \ \frac{1}{2}e_2,\ \ -\frac{1}{2}e_2, \ \ \ \ \ \ \,
e_1, \ \ \ e_3, \ \ \ \frac{1}{3}e_3, \ \ \ -\frac{1}{3}e_3, \ \ \ \ \ \ \, e_1, \  \ \ e_4,\ \ \ \frac{1}{4}e_4, \ \ \ -\frac{1}{4}e_4,\ldots)$.

     \vspace{.1in} \hspace{.46in}
(i) Let \  
  $m=(\ \ \frac{1}{2},\ \ \ \ \  1, \ \ \ \ \ \ \ 1, \ \ \ \ \ \ \ \ \ 1, \ \ \ \ \ \ \  \frac{1}{2^2},  \ \ \ \ \, 1, \ \ \ \ \ \ \ 1, \ \ \ \ \ \ \ \ \ \ 1, \ \ \ \ \ \ \   \frac{1}{2^3}, \ \ \ \ \, 1, \ \ \ \ \ \ \, 1, \ \ \ \ \ \ \ \ \ \ 1,  \ldots)$.

  \vspace{.05in}
  \noindent  
 Then  
    $M_{m,\Phi,\Psi}=M_{m,\Psi,\Phi}=I$. The convergence is unconditional on $\h$, because $M_{m,\Phi,\Psi}\neweq M_{(1),(\sqrt{m_n}\,\phi_n),(\sqrt{m_n}\,\psi_n)}$,  $M_{m,\Psi,\Phi}\neweq M_{(1),(\sqrt{m_n}\,\psi_n),(\sqrt{m_n}\,\phi_n)}$ and the sequences $(\sqrt{m_n}\,\phi_n)$, $(\sqrt{m_n}\,\psi_n)$ are Bessel for $\h$ (apply Prop. \ref{bmh}).
    
    \vspace{.1in} \hspace{.43in}
(ii) Let \ $m=(\ \ \frac{1}{2},\ \ \ \ \, 1, \ \ \ \ \ \ \ 2, \ \ \ \ \ \ \ \ \ 2, \ \ \ \ \ \ \ \frac{1}{2^2}, \ \ \ \ \, 1, \ \ \ \ \ \ \ 3, \ \ \ \ \ \ \ \ \ \ 3, \ \ \ \ \ \ \ \frac{1}{2^3}, \ \ \ \ \ 1, \ \ \ \ \ \ \, 4, \ \ \ \ \ \ \ \ \ \, 4,  \ \ldots)$.

   \vspace{.05in}
 \noindent 
 Then  
    $M_{m,\Phi,\Psi}=M_{m,\Psi,\Phi}=I$. 
 The convergence is unconditional on $\h$, because
    $M_{m,\Phi,\Psi}\neweq M_{m,\Psi,\Phi}\neweq M_{\nu,\Theta,\Theta}$, where 
    
      \vspace{.05in}
    $\Theta=(\frac{1}{\sqrt{2}}e_1, \ \ \ \ \, e_2, \ \ \ \ e_2, \ \ \ \ \ e_2,\ \ 
\frac{1}{\sqrt{2^2}}\,e_1, \ \ \ \ \ \ e_3, \ \ \ \ e_3, \ \ \, \ \ \ e_3,\ \ \frac{1}{\sqrt{2^3}}\,e_1, \ \ \ \ \  \ e_4, \ \  \ \ e_4, \ \ \ \ e_4, \ \ldots)$ is Bessel for $\h$ and 

  \vspace{.05in} $\nu\, =(\ \ \ \ \ 1, \ \ \ \ \ \ 1, \ \ \ \ \   1, \ \ \ \ -1, \ \ \ \ \ \ \ \ \ 1, \ \ \ \ \  \ \ \, 1, \ \ \ \ \ 1, \ \ \ \ \, -1, \ \ \ \  \ \ \ \ \, 1, \ \ \ \ \ \ \ \, 1, \ \ \ \ \ 1, \ \ \ \, -1, \ \ldots)$ is $SN$ (apply Prop. \ref{bmh}).

   \vspace{.1in} \hspace{.4in}
(iii) Let \
  $m=(\, \ \frac{1}{2},\ \ \ \ \,  \frac{1}{2},\ \ \ \ \ \ \  2, \ \ \ \ \ \ \ \ \, 2, \ \ \ \ \ \ \ \frac{1}{2^2}, \ \ \ \  \frac{1}{3}, \ \ \ \ \ \ \  3, \ \ \ \ \ \ \ \ \ \ 3, \ \ \ \ \ \ \, \frac{1}{2^3},  \ \ \ \ \, \frac{1}{4}, \ \ \ \ \ \ \ 4, \  \ \ \ \ \ \ \ \ 4, \ \ldots)$.

  \vspace{.05in}
\noindent 
  Then  $M_{m,\Phi,\Psi}=M_{m,\Psi,\Phi}=G_1$ - non-invertible on $\h$ (see Lemma \ref{lemg}).  
The convergence is unconditional on $\h$, because
    $M_{m,\Phi,\Psi}\neweq M_{m,\Psi,\Phi} \neweq M_{\nu,\Theta,\Theta}$, where 
    
      \vspace{.05in}
    $\Theta=(\frac{1}{\sqrt{2}}e_1, \  \frac{1}{\sqrt{2}}e_2, \ \ \ \ \, e_2,\  \ \ \ \ e_2,\ \ 
\frac{1}{\sqrt{2^2}}\,e_1, \ \ \frac{1}{\sqrt{3}}e_3, \ \ \ \, e_3, \ \ \ \ \, \ \ e_3,\ \ \frac{1}{\sqrt{2^3}}\,e_1, \ \ \frac{1}{\sqrt{4}}e_4, \ \  \  e_4, \ \ \  e_4, \ \ldots)$ is Bessel for $\h$ and 

  \vspace{.05in} $\nu \, =(\ \ \ \ \ \, 1, \ \ \ \ \ \  1, \ \ \  \ \ \ 1, \ \ \ \, -1, \  \ \ \ \ \ \ \ \ 1, \ \ \ \ \ \ \ 1, \ \ \ \ \   1, \ \ \ \ \  -1, \ \ \ \ \ \ \ \ \ 1, \ \ \ \ \ \ \  1, \ \ \ \ \, 1, \ \   -1, \ \ldots)$ is $SN$ (apply Prop. \ref{bmh}).

\end{ex}

\begin{ex}\label{nbnb10} Let 
  $\Phi=(\,e_1, \ \ \ \ \ \, e_2,  \ \ \ \ e_1, \ \ \ \ \ \ e_3, \ \ \ \  \, e_1, \ \ \ \ \ e_4, \ \ldots)$,
  
  \vspace{.05in}\hspace{.91in}
 $\Psi=(\,e_1,\ \ \ \frac{1}{2} e_2, \ \ \ \ e_1,\ \ \ \ \frac{1}{3}e_3, \ \ \ \ 
e_1, \ \ \  \frac{1}{4}e_4, \ \ldots)$.
  
  \vspace{.05in}\hspace{.9in}
  $m=(\ \, \frac{1}{2},\ \ \ \ \ \ \ 1, \ \ \ \   \frac{1}{2^2}, \ \ \ \ \ \ \ \ 1,  \ \ \ \  \frac{1}{2^3}, \ \ \ \ \ \ \, 1,\  \ldots)$.
  
 \vspace{.05in}
  \noindent  Then  
    $M_{m,\Phi,\Psi}=M_{m,\Psi,\Phi}=G_1$ - non-invertible on $\h$ (see Lemma \ref{lemg}). The unconditional convergence follows as in Example \ref{nbnb9}(i).
 \end{ex}

\begin{ex}\label{nbnb21}
 Let 
  $\Phi=(\,e_1, \ \ \ \  \, e_2,  \ \ \ \ \ \, e_1, \ \ \ \ \ \, e_3, \ \ \ \  \ e_1, \ \ \ \ \, e_4, \ \ldots)$,
  
  \vspace{.05in} \hspace{.91in}
 $\Psi=(\,e_1,\ \ \ 2 e_2, \ \ \ \ \ e_1,\ \ \ \ 3e_3, \ \ \ \ \ e_1, \ \ \  4e_4, \ \ldots)$.

  \vspace{.1in} \hspace{.47in}
(i) Let \ 
  $m=(\ \, \frac{1}{2},\ \ \ \ \  \  \frac{1}{2},  \ \ \ \  \ \frac{1}{2^2}, \ \  \ \ \ \  \, \frac{1}{3},  \ \ \ \ \ \frac{1}{2^3}, \ \ \ \ \ \ \frac{1}{4},   \  \ldots)$.

  \vspace{.05in}
\noindent  
  Then  
    $M_{m,\Phi,\Psi}=M_{m,\Psi,\Phi}=I$. 
  The convergence is unconditional on $\h$, because
    $M_{m,\Phi,\Psi}\neweq M_{m,\Psi,\Phi} \neweq M_{(1),\Theta,\Theta}$, where 

 \vspace{.05in}
    $\Theta=(\frac{1}{\sqrt{2}}e_1, \  e_2, \ \ 
\frac{1}{\sqrt{2^2}}\,e_1, \  e_3, \ \ \frac{1}{\sqrt{2^3}}\,e_1, \  e_4,\ \ldots)$ is Bessel for $\h$ (apply Prop. \ref{bmh}).

   \vspace{.1in} \hspace{.44in}
(ii) Let \
  $m=(\ \, \frac{1}{2},\ \ \ \ \ \frac{1}{2^2},  \ \ \ \  \, \frac{1}{2^2}, \ \  \ \ \ \, \frac{1}{3^2},  \ \ \ \ \ \frac{1}{2^3}, \ \ \ \   \frac{1}{4^2},    \ \ldots)$.

    \vspace{.05in}
 \noindent  
  Then  $M_{m,\Phi,\Psi}=M_{m,\Psi,\Phi}=G_1$ - non-invertible on $\h$ (see Lemma \ref{lemg}). 
The convergence  is unconditional on $\h$, because
    $M_{m,\Phi,\Psi}\neweq M_{m,\Psi,\Phi} \neweq M_{(1),\Theta,\Theta}$, where 

 \vspace{.05in}
    $\Theta=(\frac{1}{\sqrt{2}}e_1, \  \, \frac{1}{\sqrt{2}}e_2, \ \ 
\frac{1}{\sqrt{2^2}}\,e_1, \ \, \frac{1}{\sqrt{3}} e_3, \ \ \frac{1}{\sqrt{2^3}}\,e_1, \ \, \frac{1}{\sqrt{4}} e_4,\ \ldots)$ is Bessel for $\h$ (apply Prop. \ref{bmh}).  
\end{ex}

\begin{ex}
\label{nbnbnew1} 
Let 
  $\Phi=(\,e_1, \ \ \ e_2, \ \ \ \ \, e_2, \ \ \ \ \ \, \ e_2, \ \ \ \  e_1, \ \ \ e_3,\ \ \ \ \ e_3, \ \ \ \ \ \ \ \ \, e_3, \ \ \ \ e_1, \ \ \  e_4, \ \ \ \ \, e_4, \ \ \ \ \ \ \ \ e_4, \ \ \ \ e_1, \ \ \  e_5, \ \ \ \ \ \, e_5, \ \ \ \ \ \ \ \ e_5,\ \ldots)$,
  
  \vspace{.05in} \hspace{.91in}
 $\Psi=(\,e_1,\ \ \ e_2, \ \ \, 2e_2,\ \, -2e_2, \ \ \ \ e_1, \ \ \ e_3, \ \ \, \frac{1}{3}e_3, \ \ \, -\frac{1}{3}e_3, \ \ \ \ e_1, \  \ \ e_4,\ \ \ 4e_4, \ \ \, -4e_4,\ \ \ \ e_1, \ \ \ e_5, \ \ \ \frac{1}{5}e_5, \ \ \ -\frac{1}{5}e_5, \ \ldots)$.

      \vspace{.1in} \hspace{.47in}
(i) Let \ 
  $m=(\ \, \frac{1}{2},\ \ \ \ \, 1, \ \ \ \ \  \, \frac{1}{2}, \ \ \ \ \ \ \ \ \frac{1}{2}, \ \ \ \,  \frac{1}{2^2},  \ \ \ \ \ 1, \ \ \ \ \ \, \ 1, \ \ \ \ \ \ \ \ \ \ 1, \ \ \ \,  \frac{1}{2^3}, \ \ \ \ \, 1, \ \ \ \ \  \ \frac{1}{4}, \ \ \ \ \  \ \ \ \ \frac{1}{4}, \ \ \ \,  \frac{1}{2^4},  \ \ \ \ \, 1, \ \ \ \ \ \ \ 1, \ \ \ \ \ \ \ \ \ \ 1, \ \ldots)$.

   \vspace{.05in}
  \noindent  
 Then  
    $M_{m,\Phi,\Psi}=M_{m,\Psi,\Phi}=I$. The convergence is unconditional on $\h$, because  $M_{m,\Phi,\Psi}\neweq M_{m,\Psi,\Phi} \neweq M_{\nu,\Theta,\Theta}$, where 

     \vspace{.05in}
    $\Theta=(\frac{1}{\sqrt{2}}e_1, \  \ \ e_2, \ \ \ \ \, e_2,\  \ \ \ \ e_2,\ \ 
\frac{1}{\sqrt{2^2}}\,e_1, \ \ \ \ e_3,  \ \frac{1}{\sqrt{3}} e_3,  \ \  \frac{1}{\sqrt{3}} e_3,\ \ \frac{1}{\sqrt{2^3}}\,e_1, \ \ \ \ e_4, \ \  \  e_4, \ \ \  e_4, \ \ 
\frac{1}{\sqrt{2^4}}\,e_1, \ \ \ \ e_5,  \ \frac{1}{\sqrt{5}} e_5,  \ \  \frac{1}{\sqrt{5}} e_5, \ \ldots)$ is Bessel for $\h$ and 

  \vspace{.05in} $\nu \, =(\ \ \ \ \ \ \  1, \ \ \ \ \  1, \ \ \  \ \ \ 1, \ \ \ \, -1, \ \ \ \ \ \ \ \ \ \ \ 1, \, \ \ \ \ \ 1, \ \ \ \ \ \ \ \ 1, \ \ \ \ \, \  -1, \ \ \ \ \ \ \ \ \ \ \ 1, \  \ \ \ \   1,\ \ \ \ \ \, 1, \ \,  -1, \ \ \ \ \ \ \ \ \ \ \ 1, \, \ \ \ \ \, 1, \ \ \ \ \ \ \ \ 1, \ \ \ \ \ \ -1, \ \ldots)$ is $SN$ (apply Prop. \ref{bmh}).

  \vspace{.1in} \hspace{.44in}
(ii) Let \
   $m=(\ \, \frac{1}{2},\ \ \ \ \, 1, \ \ \ \ \  \, \frac{1}{2},  \ \ \ \ \ \ \ \ \frac{1}{2}, \ \ \ \,  \frac{1}{2^2},  \ \ \ \ \, 1, \ \ \ \ \ \ \, 3, \ \ \ \ \ \ \ \ \, \ 3, \ \ \ \,  \frac{1}{2^3}, \ \ \ \ \, 1, \ \ \ \ \  \ \frac{1}{4}, \ \ \ \ \  \ \ \ \ \, \frac{1}{4}, \ \ \ \,  \frac{1}{2^4},  \ \ \ \ \, 1, \ \ \ \ \ \ \ 5, \ \ \ \ \ \ \ \ \ \, 5, \ \ldots)$.

  \vspace{.05in}
  \noindent  
 Then  
    $M_{m,\Phi,\Psi}=M_{m,\Psi,\Phi}=I$. The unconditional convergence follows as in Example \ref{nbnb9}(ii).
  
       \vspace{.1in} \hspace{.4in}
(iii) Let \
  $m=(\ \, \frac{1}{2},\ \ \ \ \, \frac{1}{2}, \ \ \ \ \  \, \frac{1}{2}, \ \ \ \ \ \, \ \ \frac{1}{2}, \ \ \ \  \frac{1}{2^2},  \ \ \ \  \frac{1}{3}, \ \ \ \ \ \ \, 3, \ \ \ \ \ \ \ \ \ \ 3, \ \ \ \,  \frac{1}{2^3}, \ \ \ \  \frac{1}{4}, \ \ \ \ \  \ \frac{1}{4}, \ \ \ \ \  \ \, \ \ \frac{1}{4}, \ \ \ \,  \frac{1}{2^4},  \ \ \ \ \, \frac{1}{5}, \ \ \ \ \ \ \ 5, \ \ \ \ \ \ \ \ \ \, 5, \ \ldots)$.

 \vspace{.05in}
  \noindent  
 Then  
    $M_{m,\Phi,\Psi}=M_{m,\Psi,\Phi}=G_1$ - non-invertible on $\h$ (see Lemma \ref{lemg}). The unconditional convergence follows as in Example \ref{nbnb9}(iii).
    \end{ex}

\begin{ex}
\label{nbnbnew2}
 Let 
  $\Phi=(\,e_1, \ \ \ \ \ \ e_2,  \ \ \ \ e_1, \ \ \ \ \ \, e_3, \ \ \ \ \, e_1, \ \ \ \ \ e_4, \ \ \ \ e_1, \ \ \ \ \ \ \, e_5, \ \ \ \  e_1, \ \ \ \ \ e_6, \ \ \ \ e_1, \ \ \ \ \ \ \, e_7,\ \ldots)$,
  
  \vspace{.05in}\hspace{.91in} 
 $\Psi=(\,e_1,\ \ \ \, 2 e_2, \ \ \ \ e_1,\ \ \ \ \frac{1}{3}e_3, \ \ \ \ e_1, \ \ \  4e_4, \ \ \ \ e_1,\ \ \ \ \frac{1}{5}e_5, 
 \ \ \ \ e_1, \ \ \ \, 6e_6, \ \ \ \ e_1,\ \ \ \, \frac{1}{7}e_7, \ \ldots)$,
  
  \vspace{.05in}\hspace{.9in}
  $m=(\ \, \frac{1}{2},\ \ \ \ \  \frac{1}{2^2},  \ \ \ \  \frac{1}{2^2}, \ \  \ \ \ \ \ \ 1,  \ \ \ \,  \frac{1}{2^3}, \ \ \ \, \ \frac{1}{4^2}, \ \ \ \,  \frac{1}{2^4}, \ \  \ \ \ \ \ \, 1, \ \ \ \  \frac{1}{2^5}, \ \ \ \ \ \frac{1}{6^2}, \ \ \ \,  \frac{1}{2^6}, \ \  \ \ \ \ \ \ 1,  \  \ldots)$.
  
  \vspace{.05in}
  \noindent  Then  
    $M_{m,\Phi,\Psi}=M_{m,\Psi,\Phi}=G_1$ - non-invertible on $\h$ (see Lemma \ref{lemg}). The unconditional convergence on $\h$ follows in the same way as in Example \ref{nbnb9}(i).
    \end{ex}

\begin{ex}\label{nbnbnew5}
 Let 
  $\Phi=(\frac{1}{2}e_1, \ \ \ \,  e_2, \ \ \ \ \ \ e_1, \ \ \ \  e_3, \ \ \ \frac{1}{2^3}  e_1, \ \ \ \  e_4, \ \ \ \ \ \ e_1, \ \ \ \ e_5, 
   \ \ \ \frac{1}{2^5} e_1, \ \ \ \, e_6, \ \ \ \ \ \ e_1, \ \ \ \ \,  e_7, \ \ldots)$,
  
   \vspace{.1in} \hspace{.49in}
(i) Let \  
 $\Psi=(\ \  e_1,\ \ \ \, e_2, \ \ \ \frac{1}{2^2} e_1,\ \ \ \ e_3, \ \ \ \ \ \ e_1, \ \ \ \  e_4, \ \ \ \frac{1}{2^4} e_1,\ \ \ \  e_5, 
 \ \ \ \ \ \ e_1, \ \ \ \, e_6, \ \ \ \frac{1}{2^6} e_1,\ \ \ \ \, e_7, \ \ldots)$.

        \vspace{.05in}
  \noindent  
 Then  
    $M_{(1),\Phi,\Psi}=M_{(1),\Psi,\Phi}=I$. The convergence is unconditional on $\h$, because 
    $M_{(1),\Phi,\Psi}\neweq M_{(1),\Psi,\Phi} \neweq M_{(1),\Theta,\Theta}$, where $\Theta$ is the same as in Example \ref{nbnb21}(i).

    \vspace{.1in} \hspace{.45in}
(ii) Let \ 
 $\Psi=(\ \  e_1,\ \,  \frac{1}{2} e_2, \ \ \ \frac{1}{2^2} e_1,\ \,   \frac{1}{3} e_3, \ \ \ \ \ \ e_1, \ \   \frac{1}{4}  e_4, \ \  \frac{1}{2^4} e_1,\ \,   \frac{1}{5}  e_5, 
 \ \ \ \ \ \, e_1, \ \, \frac{1}{6} e_6, \ \ \, \frac{1}{2^6} e_1,\ \ \, \frac{1}{7} e_7, \ \ldots)$.

\vspace{.05in}
  \noindent  
 Then  
    $M_{(1),\Phi,\Psi}=M_{(1),\Psi,\Phi}=G_1$ - non-invertible on $\h$ (see Lemma \ref{lemg}). The convergence is unconditional on $\h$, because 
    $M_{(1),\Phi,\Psi}\neweq M_{(1),\Psi,\Phi} \neweq M_{(1),\Theta,\Theta}$, where $\Theta$ is the same as in Example \ref{nbnb21}(ii) (apply Prop. \ref{bmh}).   
    
    \end{ex}

\begin{ex}\label{nbnbnew6}
 Let 
  $\Phi=(e_1, \ \ \ \,  e_2, \ \ \ \ \frac{1}{2} e_1, \ \ \ \  e_3, \ \ \ \  e_1, \ \ \ \,  e_4, \ \ \  \ \frac{1}{2^2} e_1, \ \ \ \, e_5, 
   \ \ \  e_1, \ \ \ \, e_6, \ \ \ \ \frac{1}{2^3} e_1, \ \ \ \,  e_7, \ \ldots)$,

  \vspace{.05in}\hspace{.96in} 
 $m =(\ \frac{1}{2},\ \ \ \ \ 1, \ \ \ \ \ \ \ \ \, 1,\ \ \ \ \ \, 1, \ \ \ \,  \frac{1}{2^3}, \ \ \ \ \ 1, \ \ \ \ \ \ \ \ \ \, 1,\ \ \ \ \ \  1, 
\ \ \,   \frac{1}{2^5}, \ \ \ \ \ 1, \ \ \ \ \ \ \ \ \ \, 1,\ \ \ \ \ \, 1, \ \ldots)$.
  
   \vspace{.05in}
  \noindent  
 Then  
    $M_{(1),\Phi,\Phi}=I$. The convergence is unconditional on $\h$, because 
    $M_{(1),\Phi,\Phi}\neweq  M_{(1),(\sqrt{m_n}\phi_n),\sqrt{m_n}\phi_n)}$ and $(\sqrt{m_n}\phi_n)$ 
    is Bessel for $\h$ (apply Prop. \ref{bmh}).
   
\end{ex}

\begin{ex}\label{nbnbnew7}
Let
  $\Phi=(\,e_1,\ \ \ \frac{1}{2} e_2, \ \ \ \ e_1,\ \ \ \ \frac{1}{3}e_3, \ \ \ \ 
e_1, \ \ \  \frac{1}{4}e_4, \ \ldots)$.
  
     \vspace{.1in} \hspace{.54in}
(i) Let \  $m=(\ \, \frac{1}{2},\ \ \ \ \ \ \ 1,  \ \ \ \ \,  \frac{1}{2^2}, \ \  \ \ \ \ \ \, 1,  \ \ \ \  \frac{1}{2^3}, \ \ \ \ \ \ \, 1,\  \ldots)$.

   \vspace{.05in}
  \noindent  Then  
    $M_{m,\Phi,\Phi}=G_2$ - non-invertible on $\h$ (see Lemma \ref{lemg}). The unconditional convergence on $\h$ follows in the same way as in Example \ref{nbnbnew6}.

       \vspace{.1in} \hspace{.5in}
(ii) Let \  $m=(\ \, \frac{1}{2},\ \ \ \ \ 2^2,  \ \ \ \ \, \frac{1}{2^2}, \ \  \ \ \ \ 3^2,  \ \ \ \  \frac{1}{2^3}, \ \ \ \ \, 4^2,\  \ldots)$.

  \vspace{.05in}
  \noindent  Then  
    $M_{m,\Phi,\Phi}=I$. 
    The unconditional convergence follows as in Example \ref{nbnbnew6}.

     \vspace{.1in} \hspace{.46in}
(iii) Let \  $m=(\ \, \frac{1}{2},\ \ \ \ \ \ \, 2,  \ \ \ \ \, \frac{1}{2^2}, \ \ \ \ \ \ \ \ 3,  \ \ \ \   \frac{1}{2^3}, \ \ \ \ \ \ \, 4,\  \ldots)$.

    \vspace{.05in}
  \noindent  Then  
    $M_{m,\Phi,\Phi}=G_1$ - non-invertible on $\h$ (see Lemma \ref{lemg}). The unconditional convergence on $\h$ follows in the same way as in Example \ref{nbnbnew6}.

\end{ex}

\begin{ex}\label{nbnbnew8}
 Let 
  $\Phi=(\,e_1, \ \ \ \frac{1}{2}  e_2, \ \  \,e_1, \ \ \ \, \frac{1}{3}  e_3,\ \  \,e_1, \ \ \ \frac{1}{4}  e_4, \ \ \,e_1, \ \ \ \frac{1}{5}  e_5,\ \ldots)$,
  
  \vspace{.05in}\hspace{.99in}
 $\Psi=(\,e_1,\ \ \ 2 \,e_2, \ \ \,e_1,\ \ \ 3 \,e_3, \ \ \,e_1,\ \ \ 4 \,e_4,\ \  \,e_1,\ \ \ 5 \,e_5,\ \ldots)$.

   \vspace{.1in} \hspace{.55in}
(i) Let \ 
     $m=(\ \, \frac{1}{2},\ \ \ \ \ \ \ 1,  \ \ \ \frac{1}{2^2}, \ \ \ \ \ \ \  1, \ \ \,  \frac{1}{2^3}, \ \ \ \ \ \ \  1, \ \ \,  \frac{1}{2^4}, \ \ \ \ \ \ \ 1, \  \ldots)$.

   \vspace{.05in}
  \noindent  Then  
    $M_{m,\Phi,\Psi}=M_{m,\Psi,\Phi}=I$. The unconditional convergence follows as in Example \ref{nbnb21}(i).

     \vspace{.1in} \hspace{.51in}
(ii) Let \ 
$m=(\ \, \frac{1}{2},\ \ \ \ \ \ \, \frac{1}{2},  \ \ \ \frac{1}{2^2}, \ \ \ \ \ \ \,  \frac{1}{3}, \ \ \,  \frac{1}{2^3}, \ \ \ \ \ \ \,  \frac{1}{4}, \ \ \,  \frac{1}{2^4}, \ \ \ \ \ \ \, \frac{1}{5}, \  \ldots)$.

  \vspace{.05in}
  \noindent  Then  
    $M_{m,\Phi,\Psi}=M_{m,\Psi,\Phi}=G_1$ - non-invertible on $\h$ (see Lemma \ref{lemg}). The unconditional convergence follows as in Example \ref{nbnb21}(ii).

\end{ex}

\begin{ex}\label{nbnb1718}
Let
 $\Psi=(e_1, \ \ \ \, 2\,e_2,\ \ \ \ e_1,\ \ \ \ \ \, e_3,\ \ \ \ \ e_1,\ \ \ 4\,e_4,\ \ \ \ \ e_1,\ \ \ \ \ \, e_5,\ \ldots)$,

 \vspace{.05in}\hspace{.98in}
 $m=(\,\frac{1}{2},\ \ \ \ \ \ \ \, 1,\ \  \ \, \frac{1}{2^2},\ \ \ \ \ \ \ \,  3,\ \ \ \ \, \frac{1}{2^3},\ \ \ \ \ \ \ 1, \ \ \ \ \ \frac{1}{2^4}, \ \ \ \ \ \ \, 5, \ \ldots)$.
 
\vspace{.1in} \hspace{.55in}
(i) Let \  $\Phi=(e_1,\ \ \ \, \frac{1}{2}\,e_2,\ \ \ \ e_1, \ \ \ \frac{1}{3}\,e_3,\ \ \ \ \, e_1, \ \ \, \frac{1}{4}\,e_4, \ \ \ \ \  e_1, \ \ \ \frac{1}{5} e_5, \ \ldots)$.

   \vspace{.05in}
\noindent 
Then  $M_{m,\Phi,\Psi}=M_{m,\Psi,\Phi}=I$. 
 The unconditional convergence follows as in Example \ref{nbnb21}(i).

 \vspace{.1in} \hspace{.51in}
(ii) Let \  $\Phi=(\,e_1,\ \ \frac{1}{2^2}\,e_2, \ \ \ e_1, \ \ \, \frac{1}{3^2}\,e_3, \ \ \ \,e_1, \  \, \frac{1}{4^2}\,e_4,\ \ \ \ \, e_1, \ \ \frac{1}{5^2} e_5, \ \ldots)$.

\vspace{.05in}
\noindent Then  $M_{m,\Phi,\Psi}=M_{m,\Psi,\Phi}=G_1$ - non-invertible on $\h$ (see Lemma \ref{lemg}).  
 The unconditional convergence follows as in Example \ref{nbnb21}(ii).

\end{ex}

\begin{ex}\label{nbnb13}
Let
$\Phi=(\ \ \, e_1,\ \  \frac{1}{2}\,e_2, \ \ \ \ \ \ \ e_1,\ \ \frac{1}{3}\,e_3,\ \ \ \ \ \ \ \,e_1,\ \ \frac{1}{4}\,e_4,\ \ldots)$,

\vspace{.05in} \hspace{.98in}
 $\Psi=(\frac{1}{2}\,e_1, \ \ \ 2e_2, \ \ \frac{1}{2^2}\,e_1, \ \ \ 3e_3, \ \ \ \frac{1}{2^3}\,e_1,\ \ \ 4e_4,\ \ldots)$.

\vspace{.1in} \hspace{.52in}
(i) 
Then  $M_{(1),\Phi,\Psi}=M_{(1),\Psi,\Phi}=I$.  
The convergence  is unconditional on $\h$, because  $M_{(1),\Phi,\Psi}\neweq M_{(1),\Psi,\Phi} \neweq M_{(1),\Theta,\Theta}$, where
$\Theta$ is the same as in Example \ref{nbnb21}(i) \,(apply Prop. \ref{bmh}).

   \vspace{.1in} \hspace{.51in}
(ii) Let \  $m=(\ \ \ \ \, 1,\ \ \ \ \ \ \frac{1}{2}, \ \ \ \ \ \ \ \  \ 1, \ \ \ \ \ \ \frac{1}{3}, \ \ \ \  \ \ \ \ \ 1, \ \ \ \ \ \ \, \frac{1}{4},\  \ldots)$.

 \vspace{.05in} \noindent
Then  $M_{m,\Phi,\Psi}=M_{m,\Psi,\Phi}=G_1$ - non-invertible on $\h$ (see Lemma \ref{lemg}).  
The unconditional convergence follows as in Example \ref{nbnb21}(ii).

\end{ex}

\begin{ex}\label{nbnb14}
Let
 $\Phi=(\ \ \ \, e_1, \ \, \frac{1}{2^2}\,e_2,\ \ \ \ \ \ e_1, \ \ \frac{1}{3^2}\,e_3, \ \ \ \ \ \ e_1, \ \ \frac{1}{4^2}\,e_4,\ \ldots)$,

\vspace{.05in} \hspace{.99in}
$\Psi=(\frac{1}{2}\,e_1,\ \ \ \ 2e_2, \ \, \frac{1}{2^2}\,e_1,\ \ \ \ \, 3e_3,\ \ \frac{1}{2^3}\,e_1,\ \ \ \ 4e_4,\ \ldots)$.

 \vspace{.1in} \hspace{.54in}
(i) 
Then  $M_{(1),\Phi,\Psi}=M_{(1),\Psi,\Phi}=G_1$ - non-invertible on $\h$ (see Lemma \ref{lemg}).   
The convergence is unconditional on $\h$, because  $M_{(1),\Phi,\Psi}\neweq M_{(1),\Psi,\Phi} \neweq M_{(1),\Theta,\Theta}$, where
$\Theta$ is the same as in Example \ref{nbnb21}(ii) \,(apply Prop. \ref{bmh}).

    \vspace{.1in} \hspace{.51in}
(ii) Let \  $m=(\ \ \ \, \ 1,\ \ \ \ \ \ \ \ 2, \ \  \ \ \ \ \ \ \, 1, \ \ \ \ \ \ \ \, 3,  \ \ \ \ \ \ \ \ \, 1, \ \ \ \ \ \ \ \ 4,\  \ldots)$.

 \vspace{.05in} \noindent
Then  $M_{m,\Phi,\Psi}=M_{m,\Psi,\Phi}=I$.  
The unconditional convergence follows as in Example \ref{nbnb21}(i).

\end{ex}

\begin{ex}\label{nbnbnew9} 
Let
$\Phi=(\ \ \ \ \ \,e_1,\ \ \, \frac{1}{2^2}\,e_2, \ \ \ \ \ \ \ \ \ e_1,\ \ \frac{1}{3^2}\,e_3,\ \ \ \ \ \ \ \ \ \ \, e_1,\ \ \, \frac{1}{4^2}\,e_4,\ \ldots)$,

\vspace{.05in} \hspace{.99in}
 $\Psi=(\frac{1}{\sqrt{2}}\,e_1, \ \ \ 2^2e_2, \ \ \frac{1}{\sqrt{2^2}}\,e_1, \ \ \ 3^2e_3, \ \ \ \frac{1}{\sqrt{2^3}}\,e_1,\ \ \ \ 4^2e_4,\ \ldots)$,

\vspace{.05in} \hspace{.99in}
  $m=(\ \ \ \ \frac{1}{\sqrt{2}},\ \ \ \ \ \ \ \ \, 1, \ \  \ \ \ \ \frac{1}{\sqrt{2^2}}, \ \ \ \ \ \ \ \, \ 1, \ \ \ \  \ \ \  \frac{1}{\sqrt{2^3}}, \ \  \ \ \ \ \ \ \ \, 1,\  \ldots)$.

\vspace{.05in} \noindent
Then  $M_{m,\Phi,\Psi}=M_{m,\Psi,\Phi}=I$.  
The unconditional convergence follows as in Example \ref{nbnb21}(i).

\end{ex}

\begin{ex}
\label{nbnbnew9new} 
Let 
$\Phi=(\ \ \ \ \ \ e_1,\ \  \frac{1}{2^3}\,e_2, \ \ \ \ \ \ \ \ \ e_1,\ \ \frac{1}{3^3}\,e_3,\ \ \ \ \ \ \ \ \ \ \,e_1,\ \ \, \frac{1}{4^3}\,e_4,\ \ldots)$,

\vspace{.05in} \hspace{.99in}
 $\Psi=(\frac{1}{\sqrt{2}}\,e_1, \ \ \ \ 2e_2, \ \ \, \frac{1}{\sqrt{2^2}}\,e_1, \ \ \ \ 3e_3, \ \ \ \, \frac{1}{\sqrt{2^3}}\,e_1,\ \ \ \ \ 4e_4,\ \ldots)$,

\vspace{.05in} \hspace{.99in}
  $m=(\ \ \ \ \frac{1}{\sqrt{2}},\ \ \ \ \ \ \, \ 2, \ \  \ \ \ \, \ \frac{1}{\sqrt{2^2}}, \ \ \ \ \ \ \ \ 3,  \ \ \  \ \ \ \ \ \frac{1}{\sqrt{2^3}}, \ \ \ \ \ \ \ \, \ 4,\  \ldots)$.

\vspace{.05in} \noindent
Then  $M_{m,\Phi,\Psi}=M_{m,\Psi,\Phi}=G_1$ - non-invertible on $\h$ (see Lemma \ref{lemg}).  
The unconditional convergence follows as in Example \ref{nbnb21}(ii).
\end{ex}

\begin{ex}\label{nbnb15i} 
 Let $\Phi=(ne_n)$ and $m=(\frac{1}{n^2})$. Then  $M_{m,\Phi,\Phi}\neweq M_{(1),(e_n),(e_n)}=I$. 
\end{ex}

\begin{ex} \label{nbnb15ii} 
 Let $\Phi=(ne_n)$ and $m=(\frac{1}{n^3})$. Then $M_{m,\Phi,\Phi}\neweq M_{(\frac{1}{n}),(e_n),(e_n)}=G_1$ 
- unconditionally convergent and non-invertible on $\h$ (see Lemma \ref{lemg}).
  \end{ex}

\begin{ex} \label{nbnbnew10}
Let
$\Phi=(e_1, \ \ \ \, 2e_2, \ \ \ \ \, e_1, \ \ \ \, 3e_3, \ \ \ \ \ e_1,\ \ \ \, 4e_4,  \ \ \ \ \  \ e_1,\ \ \ 5e_5,\ \ldots)$,

\vspace{.05in} \hspace{.99in}
 $\Psi=(  e_1,\ \  \frac{1}{2}\,e_2, \ \ \, 2 e_1,\ \ \frac{1}{3}\,e_3,\ \  2^2e_1,\ \ \frac{1}{4}\,e_4,\ \ \, 2^3e_1,\ \ \frac{1}{5}\,e_5,\ \ldots)$.

 \vspace{.1in} \hspace{.54in}
(i) Let \  $m=(\, \ \frac{1}{2},\ \ \ \ \ \ \, 1, \ \ \ \ \ \frac{1}{2^3}, \ \ \ \ \ \ \, 1,  \ \ \  \ \  \frac{1}{2^5}, \ \ \, \ \ \ \ 1, \ \ \ \ \ \, \ \frac{1}{2^7}, \ \  \ \ \ \ 1,\  \ldots)$.

\vspace{.05in} \noindent
Then  $M_{m,\Phi,\Psi}=M_{m,\Psi,\Phi}=I$.  
The unconditional convergence follows as in Example \ref{nbnb21}(i).

 \vspace{.1in} \hspace{.51in}
(ii) Let \ $m=(\, \ \frac{1}{2},\ \ \ \ \ \, \frac{1}{2}, \ \ \ \ \ \frac{1}{2^3}, \ \ \ \ \ \ \, \frac{1}{3},  \ \ \  \ \  \frac{1}{2^5}, \ \  \ \ \ \, \frac{1}{4}, \ \ \ \ \ \ \,  \frac{1}{2^7}, \ \  \ \ \ \, \frac{1}{5},\  \ldots)$.

 \vspace{.05in} \noindent
Then  $M_{m,\Phi,\Psi}=M_{m,\Psi,\Phi}=G_1$ - non-invertible on $\h$ (see Lemma \ref{lemg}).   
The unconditional convergence follows as in Example \ref{nbnb21}(ii).

\end{ex}

\begin{ex} \label{nbnbnew11}
Let
$\Phi=(e_1, \ \ \ \ \, 2e_2, \ \ \ \ \, e_1, \ \ \ \ \, 3e_3, \ \ \ \ \ \ e_1,\ \ \ \ \, 4e_4,  \ \ \ \ \ \, \ e_1,\ \ \ \ \ 5e_5,\ \ldots)$,

\vspace{.05in} \hspace{.96in}
$m=(\,  \frac{1}{2},\ \ \ \ \ \ \ \ 2, \ \  \ \ \ \frac{1}{2^3}, \ \ \ \ \ \ \ \ 3,  \ \ \  \ \, \ \frac{1}{2^5}, \ \ \ \ \ \ \ \ \, 4, \ \ \ \ \ \  \frac{1}{2^7}, \ \ \ \ \ \ \ \ \, 5,\  \ldots)$.

 \vspace{.1in} \hspace{.56in}
(i) Let \ $\Psi=(  e_1,\ \  \frac{1}{2^2}\,e_2, \ \ \, 2 e_1,\ \ \frac{1}{3^2}\,e_3,\ \  2^2e_1,\ \ \frac{1}{4^2}\,e_4,\ \ \, 2^3e_1,\ \ \, \frac{1}{5^2}\,e_5,\ \ldots)$.

\vspace{.05in} \noindent
Then  $M_{m,\Phi,\Psi}=M_{m,\Psi,\Phi}=I$.  
The unconditional convergence follows as in Example \ref{nbnb21}(i).

    \vspace{.1in} \hspace{.52in}
(ii) Let \ $\Psi=(  e_1,\ \  \frac{1}{2^3}\,e_2, \ \ \, 2 e_1,\ \ \frac{1}{3^3}\,e_3,\ \  2^2e_1,\ \ \frac{1}{4^3}\,e_4,\ \ \, 2^3e_1,\ \ \, \frac{1}{5^3}\,e_5,\ \ldots)$.
 
\vspace{.05in} \noindent
Then  $M_{m,\Phi,\Psi}=M_{m,\Psi,\Phi}=G_1$ - non-invertible on $\h$ (see Lemma \ref{lemg}).  
The unconditional convergence follows as in Example \ref{nbnb21}(ii).
 
 \end{ex}

\begin{ex} \label{nbnbnew12}
Let
$\Phi=(\frac{1}{\sqrt{2}}\, e_1, \ \ \ \,  2\,e_2, \ \ \  \frac{1}{\sqrt{2^2}}\, e_1, \ \ \ \frac{1}{3}\,e_3, \ \ \  \frac{1}{\sqrt{2^3}}\, e_1,\ \ \ \, 4\,e_4,  \ \ \  \frac{1}{\sqrt{2^4}}\, e_1,\ \ \, \frac{1}{5}\,e_5,\ \ldots)$,

\vspace{.05in} \hspace{.99in}
$\Psi=( \frac{1}{\sqrt{2}}\, e_1,\ \ \ \frac{1}{2}\,e_2, \ \ \ \frac{1}{\sqrt{2^2}}\, e_1,\ \ \ \,  3\,e_3,\ \ \  \frac{1}{\sqrt{2^3}}\, e_1,\ \ \ \frac{1}{4}\,e_4, \ \ \  \frac{1}{\sqrt{2^4}}\, e_1,\ \ \ 5\,e_5,\ \ldots)$.
 
   \vspace{.1in} 
\noindent
 Then  $M_{(1),\Phi,\Psi}=M_{(1),\Psi,\Phi}=I$.  
The convergence  is unconditional on $\h$, because  $M_{(1),\Phi,\Psi}\neweq M_{(1),\Psi,\Phi} \neweq M_{(1),\Theta,\Theta}$, where
$\Theta$ is the same as in Example \ref{nbnb21}(i).

\end{ex}

\begin{ex} \label{nbnbnew13}
Let
$\Phi=(\frac{1}{\sqrt{2}}\, e_1, \ \ \ \  2\,e_2, \ \ \  \frac{1}{\sqrt{2^2}}\, e_1, \ \  \frac{1}{3^2}\,e_3, \ \ \  \frac{1}{\sqrt{2^3}}\, e_1,\ \ \ \, 4\,e_4,  \ \ \  \frac{1}{\sqrt{2^4}}\, e_1,\  \, \frac{1}{5^2}\,e_5,\ \ldots)$,

\vspace{.05in} \hspace{.99in}
$\Psi=( \frac{1}{\sqrt{2}}\, e_1,\ \  \frac{1}{2^2}\,e_2, \ \ \ \frac{1}{\sqrt{2^2}}\, e_1,\ \ \ \,  3\,e_3,\ \ \  \frac{1}{\sqrt{2^3}}\, e_1,\ \  \frac{1}{4^2}\,e_4, \ \ \  \frac{1}{\sqrt{2^4}}\, e_1,\ \ \ \, 5\,e_5,\ \ldots)$.
 
    \vspace{.1in} 
\noindent
 Then  $M_{(1),\Phi,\Psi}=M_{(1),\Psi,\Phi}=G_1$ - non-invertible on $\h$ (see Lemma \ref{lemg}).  
The convergence  is unconditional on $\h$, because  $M_{(1),\Phi,\Psi}\neweq M_{(1),\Psi,\Phi} \neweq M_{(1),\Theta,\Theta}$, where
$\Theta$ is the same as in Example \ref{nbnb21}(ii) \,(apply Prop. \ref{bmh}).
 
\end{ex}

\begin{ex} \label{nbnbnew14}
Let
$\Phi=(\frac{1}{\sqrt{2}}\, e_1, \ \ \ \,  2\,e_2, \ \ \  \frac{1}{\sqrt{2^2}}\, e_1, \ \ \ \, 3 \,e_3, \ \ \  \frac{1}{\sqrt{2^3}}\, e_1,\ \ \ \, 4\,e_4,  \ \ \  \frac{1}{\sqrt{2^4}}\, e_1,\ \ \ \, 5\,e_5,\ \ldots)$.

\vspace{.1in} \hspace{.54in}
(i) Let \ $m=( \ \ \ \ \ \ \ \, 1,\ \ \ \ \ \ \frac{1}{2^2},\ \ \ \ \ \ \ \ \ \ \ \ 1,\ \ \ \ \ \, \frac{1}{3^2}, \ \ \ \ \ \ \ \ \ \ \ \  1,\ \ \ \ \ \ \frac{1}{4^2}, \ \ \ \ \ \ \ \ \ \, \ \ \, 1,\ \ \ \ \ \, \frac{1}{5^2},\ \ldots)$.

\vspace{.05in} \noindent
 Then  $M_{m,\Phi,\Phi}=I$.  
The convergence  is unconditional on $\h$, because  $M_{m,\Phi,\Phi} \neweq M_{(1),\Theta,\Theta}$, where
$\Theta$ is the same as in Example \ref{nbnb21}(i) \,(apply Prop. \ref{bmh}).

\vspace{.1in} \hspace{.51in}
(ii) Let \ $m=( \ \ \ \ \ \ \ \, 1,\ \ \ \ \ \ \frac{1}{2^3},\ \ \ \ \ \ \ \ \ \ \ \ 1,\ \ \ \ \ \, \frac{1}{3^3}, \ \ \ \ \ \, \ \ \ \ \ \ \, 1,\ \ \ \ \ \, \frac{1}{4^3}, \ \ \ \ \ \ \ \ \ \, \ \ \, 1,\ \ \ \ \ \, \frac{1}{5^3},\ \ldots)$.

\vspace{.05in} \noindent
 Then  $M_{m,\Phi,\Phi}=G_1$ - non-invertible on $\h$ (see Lemma \ref{lemg}).  
The convergence  is unconditional on $\h$, because  $M_{m,\Phi,\Phi} \neweq M_{(1),\Theta,\Theta}$, where
$\Theta$ is the same as in Example \ref{nbnb21}(ii).

\end{ex}

 \begin{ex} \label{nbnbnew15}
 Let
$\Phi=(\frac{1}{\sqrt{2}}\, e_1, \ \ \ \,  2\,e_2, \ \ \  \frac{1}{\sqrt{2^2}}\, e_1, \ \ \ \frac{1}{3} \,e_3, \ \ \  \frac{1}{\sqrt{2^3}}\, e_1,\ \ \ \, 4\,e_4,  \ \ \  \frac{1}{\sqrt{2^4}}\, e_1,\ \ \ \frac{1}{5}\,e_5,\ \ldots)$.

\vspace{.1in} \hspace{.54in}
(i) Let \ $m=( \ \ \ \ \ \ \ \, 1,\ \ \ \  \, \ \frac{1}{2^2}, \ \ \ \ \ \ \ \ \ \ \ \ 1,\ \ \ \ \ \, \, 3^2, \ \ \ \ \ \ \ \ \ \ \ \ 1,\ \ \ \ \ \ \frac{1}{4^2}, \ \ \ \ \ \ \ \ \ \ \ \, 1,\ \ \ \ \ \ \, 5^2,\ \ldots)$.

\vspace{.05in} \noindent
 Then  $M_{m,\Phi,\Phi}=I$.  
The convergence  is unconditional on $\h$, because  $M_{m,\Phi,\Phi} \neweq M_{(1),\Theta,\Theta}$, where
$\Theta$ is the same as in Example \ref{nbnb21}(i) \,(apply Prop. \ref{bmh}).

\vspace{.1in} \hspace{.51in}
(ii) Let \ $m=( \ \ \ \ \ \ \ \, 1, \ \ \ \ \, \ \frac{1}{2^3}, \ \ \ \ \ \ \ \ \ \ \ \ 1,\ \ \ \ \ \ \, \ 3, \ \ \ \ \ \ \ \ \ \ \ \ \, 1,\ \ \ \ \ \, \frac{1}{4^3}, \ \ \ \ \ \ \ \ \ \ \ \, 1,\ \ \ \ \ \ \ \ 5,\ \ldots)$. 
  
\vspace{.05in} \noindent
 Then  $M_{m,\Phi,\Phi}=G_1$ - non-invertible on $\h$ (see Lemma \ref{lemg}).   
The convergence  is unconditional on $\h$, because  $M_{m,\Phi,\Phi} \neweq M_{(1),\Theta,\Theta}$, where
$\Theta$ is the same as in Example \ref{nbnb21}(ii).
 
 \end{ex}

\subsection{Examples for one Bessel non-frame and one non-Bessel sequence; TABLE 2 
 on page \pageref{table2}}

\begin{ex}\label{bnb4}
Let
 $\Phi= (e_2, \ \ \ e_3, \ \  \ \, e_4,\ \ \ \ e_5,\ \ \ \ \, e_6,\ \ \ e_7, \ \ldots)$,

 \vspace{.05in} \hspace{.885in}
 $m=(\, \, \frac{1}{2}, \ \ \ \ \, 1, \ \ \ \frac{1}{2^2},\ \ \ \ \ \ 1, \ \ \ \ \frac{1}{2^3}, \ \ \ \ \, 1, \ \ldots)$.

 \vspace{.1in} \hspace{.505in}
(i) Let
 $\Psi=(\,e_1,\ \ \ \, e_2,\ \ \ e_1,\ \ \ \ e_3,\ \ \ \ \, e_1,\ \ \ \, e_4,\ \ldots)$.

  \vspace{.05in}
 \noindent 
 Then $M_{m,\Phi,\Psi}$ and  $M_{m,\Psi,\Phi}$ are unconditionally convergent on $\h$, because $\Phi$ and $m\Psi$ are Bessel for $\h$ (apply Prop. \ref{bmh}).  
 The multiplier $M_{m,\Phi,\Psi}$ is not surjective on $\h$, because $\Phi$ is not complete in $\h$. The multiplier $M_{m,\Psi,\Phi}$ is not injective, for example $M_{m,\Psi,\Phi}(2e_2)=e_1=M_{m,\Psi,\Phi}(4e_4)$.
 
 \vspace{.1in} \hspace{.465in}
(ii) Let 
 $\Psi=(\, e_2, \  \frac{1}{2}e_2,\  \ \ e_2,\ \ \, \frac{1}{3}e_3, \ \ \ \, e_2, \ \, \frac{1}{4} e_4, \ \ldots)$.

  \vspace{.05in}
 \noindent  
  Then $M_{m,\Phi,\Psi}$ and  $M_{m,\Psi,\Phi}$ are unconditionally convergent on $\h$, because $\Phi$ and $m\Psi$ are Bessel for $\h$ (apply Prop. \ref{bmh}). The multiplier $M_{m,\Phi,\Psi}$ (resp. $M_{m,\Psi,\Phi}$) is not invertible on $\h$, because $\Phi$ (resp. $\Psi$) is not complete in $\h$.

\end{ex}

\begin{ex}\label{exadd2}
Let
 $\Phi= ( \, e_2,\ \ \ \ \ \ e_3, \ \ \ \, e_4, \ \ \, \ \ \ \, \, e_5, \ \ \ \ e_6, \ \ \ \ \ \ \, e_7, \ \, \ldots)$, 
 
 \vspace{.05in} \hspace{.91in}
 $\Psi=(\, e_2, \ \  \frac{1}{2^2}e_2,\  \ \ \,e_2,\ \ \ \frac{1}{3^2}e_3, \ \ \ \ e_2, \ \ \ \frac{1}{4^2} e_4, \ \ldots)$, 
 
 \vspace{.05in} \hspace{.91in}
 $m=(\ \frac{1}{2},\ \ \ \ \ \ \ \ 2, \ \ \ \,  \frac{1}{2^2}, \ \ \ \ \ \ \ \ \, 3, \ \ \ \  \frac{1}{2^3}, \ \ \ \ \ \ \ \, \, 4, \   \ldots)$.
 
  \vspace{.05in}
 \noindent The conclusion is the same as in Example \ref{bnb4}(ii). 
 \end{ex}

\begin{ex}\label{bnb10}
Let
 $\Phi= ( \ \ e_2,\ \ \ \ e_3, \ \ \ \, e_4, \ \ \ \, \, e_5,\ \ldots)$, 
 
 \vspace{.05in} \hspace{.91in}
 $\Psi=(2e_2, \ \ 3e_3,\ \ 4e_4,\ \ 5e_5,\ \ldots)$, 
 
 \vspace{.05in} \hspace{.91in}
 $m=(\ \ \ \frac{1}{2},\ \ \ \ \, \frac{1}{3}, \ \ \ \ \ \, \frac{1}{4}, \ \ \ \ \, \frac{1}{5}, \  \ldots)$.
 
  \vspace{.05in}
 \noindent 
The conclusion is the same as in Example \ref{bnb4}(ii).
 \end{ex}

\begin{ex}\label{bnbnew1}
Let
 $\Phi= ( \, e_2,\ \ \ \ \, \, e_3, \ \ \ \, \, e_4, \ \ \, \, \ \ e_5, \ \ \ \ \, \,  e_6, \ \ \ \ \ \, e_7, \ \ \ \ \ \, e_8, \, \ \ \ \ \ e_9,\ \ldots)$, 
 
 \vspace{.05in}\hspace{.91in}
 $\Psi=(\, e_2, \ \ \, \frac{1}{2}e_2,\  \ 2 e_2,\ \ \, \frac{1}{3}e_3, \ \  2^2e_2, \ \ \ \frac{1}{4} e_4, \ \  2^3e_2, \ \ \ \frac{1}{5} e_5,\ \ldots)$, 
 
 \vspace{.05in} \hspace{.91in}
 $m=(\ \frac{1}{2},\ \ \ \ \  \ \ 1, \ \ \ \, \, \frac{1}{2^3}, \ \ \ \ \ \, \, 1, \ \ \ \ \, \, \frac{1}{2^5},\, \, \ \ \ \ \ \ 1, \, \, \ \ \ \  \frac{1}{2^7}, \ \ \ \ \ \ \, \, 1,\   \ldots)$.
 
   \vspace{.05in}
 \noindent  
The conclusion is the same as in Example \ref{bnb4}(ii). 
\end{ex}

\begin{ex}\label{bnb12}
Let
  $\Phi= ( \ \, \, e_2,\ \ \ \  \, e_3, \ \ \ \, e_4, \ \ \ \ \, e_5, \ \ldots)$, 
 
 \vspace{.05in}\hspace{.91in}
 $\Psi=(2e_2, \ \ \frac{1}{3}e_3, \ \ 4e_4, \ \, \frac{1}{5}e_5, \ \ldots)$, 
 
 \vspace{.05in} \hspace{.91in}
 $m=(\ \ \, \frac{1}{2},\ \ \ \ \ \ 3,\ \ \ \ \, \, \frac{1}{4},\ \ \ \ \, \ \, 5, \  \ldots)$.
 
 \noindent 
The conclusion is the same as in Example \ref{bnb4}(ii).
\end{ex}

\begin{ex}\label{bnb13}
Let $\Phi=(e_1, \ \ \frac{1}{2}e_2, \ \ \frac{1}{2}e_2,\ \ \frac{1}{3}e_3, \ \ \frac{1}{3}e_3,\ \ \frac{1}{3}e_3,\ \ldots)$,

  \vspace{.05in} \hspace{.91in}
   $\Psi=(e_1,\  \ \ \ e_2,\ \ \ \ \ e_2,\ \ \ \ \, e_3, \ \ \ \ \, e_3, \ \ \ \, \, e_3, \ \ldots)$.

    \vspace{.1in} \hspace{.48in}
 (i) Then $M_{(1),\Phi,\Psi}=M_{(1),\Psi,\Phi}=I$. 
The convergence is unconditional on $\h$, because $M_{(1),\Phi,\Psi}\neweq M_{(1),\Psi,\Phi}\ \neweq M_{(1),\Theta,\Theta}$, where 

\vspace{.05in}
$\Theta=(e_1, \frac{1}{\sqrt{2}}e_2, \frac{1}{\sqrt{2}}e_2, \frac{1}{\sqrt{3}}e_3,\frac{1}{\sqrt{3}}e_3,\frac{1}{\sqrt{3}}e_3, \ldots)$  is Bessel for $\h$ (apply Prop. \ref{bmh}). 

    \vspace{.1in} \hspace{.47in}
 (ii)   Let   
$m=(\ \, 1, \ \ \ \ \, \ \frac{1}{2},\ \ \ \ \ \, \frac{1}{2},\ \ \ \ \ \, \frac{1}{3},\ \ \ \ \ \, \, \frac{1}{3},\ \ \ \ \, \, \frac{1}{3}, \ \ldots)$.

 \vspace{.05in} \noindent
 Then $M_{m,\Phi,\Psi}=M_{m,\Psi,\Phi}=G_1$ - non-invertible on $\h$ (see Lemma \ref{lemg}). 
  The convergence is unconditional on $\h$, because $M_{m,\Phi,\Psi}\neweq M_{m,\Psi,\Phi} \neweq M_{m,\Theta,\Theta}$, where 
  $\Theta$ is the same as in (i) \,(apply Prop. \ref{bmh}). 

 \end{ex}

\begin{ex}\label{bnb14}
Let
$\Phi=(e_1, \ \ \frac{1}{2^2}e_2, \ \ \frac{1}{2^2}e_2,\ \  \frac{1}{3^2}e_3,\ \ \frac{1}{3^2}e_3,\ \ \frac{1}{3^2}e_3, \ \ \frac{1}{4^2}e_4,\ \ \frac{1}{4^2}e_4,\ \ \frac{1}{4^2}e_4,\ \ \frac{1}{4^2}e_4,\ \ldots)$, 
 
  \vspace{.05in} \hspace{.91in}
  $\Psi=(e_1,\ \ \ \ \ \ e_2,\ \ \ \ \ \ e_2,\ \ \ \ \ \ e_3,\ \ \ \ \ \, \ e_3,\ \ \ \ \ \  e_3,\ \ \ \ \ \ e_4,\ \ \ \ \ \, \, e_4,\ \ \ \ \ \ e_4,\ \ \ \ \ \ e_4,\ \ldots)$. 
  
   \vspace{.05in}
  \noindent  
  Then $M_{(1),\Phi,\Psi}=M_{(1),\Psi,\Phi}=G_1$ - non-invertible on $\h$ (see Lemma \ref{lemg}). 
 The convergence is unconditional on $\h$, because $M_{(1),\Phi,\Psi}\neweq M_{(1),\Psi,\Phi} \neweq M_{(1),\Theta,\Theta}$, where 
  
   \vspace{.05in} 
  $\Theta=(e_1, \ \frac{1}{2}e_2,\ \frac{1}{2}e_2,\ \frac{1}{3}e_3, \ \frac{1}{3}e_3, \ \frac{1}{3}e_3, \ \frac{1}{4}e_4,\ \frac{1}{4}e_4,\ \frac{1}{4}e_4,\ \frac{1}{4}e_4,\ \ldots)$  is Bessel for $\h$ (apply Prop. \ref{bmh}). 
 \end{ex}

\begin{ex}\label{bnb15}
Let
$\Phi=(\frac{1}{\sqrt{2}}e_1, \ \ \, \frac{1}{2}\,e_2, \ \frac{1}{2}\,e_2, \ \ \frac{1}{\sqrt{2^2}}e_1, \ \ \frac{1}{3}\,e_3, \ \frac{1}{3}\,e_3, \ \frac{1}{3}\,e_3,  \ \ \frac{1}{\sqrt{2^3}}e_1, \ \ \frac{1}{4}\,e_4, \ \frac{1}{4}\,e_4, \ \frac{1}{4}\,e_4, \ \frac{1}{4}\,e_4, \ \ldots)$, 

\vspace{.05in}\hspace{.91in}
  $\Psi=(\ \ \ \ \, \, e_1, \, \ \ \ \ \, e_2, \ \ \  \, \, e_2, \ \ \ \ \ \ \ \ \ e_1, \ \ \, \ \ \   e_3, \ \ \ \  e_3, \ \ \ \, \, e_3, \, \, \ \ \ \ \ \  \ e_1, \ \, \ \ \ \,  e_4, \ \ \, \ \, e_4, \ \, \ \ \, e_4, \, \ \ \ \, e_4, \ \ldots)$, 
  
  \vspace{.05in} \hspace{.905in}
  $m=(\ \ \, \, \frac{1}{\sqrt{2}}, \ \, \ \ \ \ \ \, 1, \ \, \ \ \ \ 1, \ \ \ \ \, \ \ \frac{1}{\sqrt{2^2}},  \, \ \ \ \ \ \ 1, \, \ \ \ \ \  1,\, \ \ \ \ \ 1, \ \ \ \ \, \, \frac{1}{\sqrt{2^3}}, \ \ \ \ \ \ \, \, 1, \, \ \ \ \ \, \, 1,\, \ \ \ \ \ \, 1,\, \ \ \ \ \, 1, \ \ldots)$.
  
   \vspace{.05in}
 \noindent  
Then $M_{m,\Phi,\Psi}=M_{m,\Psi,\Phi}=I$. 
  The convergence is unconditional on $\h$, because $M_{m,\Phi,\Psi} \neweq M_{m,\Psi,\Phi} \neweq M_{(1),\Theta,\Theta}$, where
  
   \vspace{.05in} 
  $\Theta=(\frac{1}{\sqrt{2}}e_1, \ \frac{1}{\sqrt{2}}e_2,\ \frac{1}{\sqrt{2}}e_2,\frac{1}{\sqrt{2^2}}e_1, \ \frac{1}{\sqrt{3}}e_3, \ \frac{1}{\sqrt{3}}e_3, \ \frac{1}{\sqrt{3}}e_3, \frac{1}{\sqrt{2^3}}e_1, \ \frac{1}{\sqrt{4}}e_4,\ \frac{1}{\sqrt{4}}e_4,\ \frac{1}{\sqrt{4}}e_4,\ \frac{1}{\sqrt{4}}e_4,\ \ldots)$  is Bessel for $\h$ (apply Prop. \ref{bmh}). 
  
  \end{ex}

\begin{ex}\label{bnb17}
Let
  $\Psi=(\ \ \, e_1, \ \ \ \ \, \, \, \, e_2,\ \ \ \ \ \, \, e_1,\ \ \ \ \, \, e_3, \ \ \ \ \ \, \, \, e_1,\ \ \ \ \ \,  e_4,\ \ldots)$, 
  
  \vspace{.05in}\hspace{.91in}
  $m=(\ \ \ \, \, 1, \ \ \ \ \ \, \, \, \ 2, \ \ \ \ \ \ \ \,\, 1,\ \ \ \ \ \ \, 3,\ \ \ \ \ \ \, \, \, \, 1, \ \ \ \ \ \ \, \, 4,\ \ldots)$.

 \vspace{.1in} \hspace{.51in}
(i) Let \ $\Phi=(\frac{1}{2}\,e_1,  \ \ \, \frac{1}{2}\,e_2, \ \ \frac{1}{2^2}\,e_1, \ \ \frac{1}{3}\,e_3, \ \ \frac{1}{2^3}\,e_1, \ \ \, \frac{1}{4}\,e_4,\ \ldots)$.

   \vspace{.05in}
\noindent 
 Then $M_{m,\Phi,\Psi}=M_{m,\Psi,\Phi}=I$. 
 The unconditional convergence follows as in Example \ref{nbnb21}(i).
 
  \vspace{.1in} \hspace{.48in}
(ii) Let \ $\Phi=(\frac{1}{2}\,e_1,  \ \frac{1}{2^2}\,e_2, \ \ \frac{1}{2^2}\,e_1, \  \frac{1}{3^2}\,e_3, \ \ \frac{1}{2^3}\,e_1,  \ \frac{1}{4^2}\,e_4, \ \ldots)$.
 
 \vspace{.05in}
\noindent  
Then $M_{m,\Phi,\Psi}=M_{m,\Psi,\Phi}=G_1$ - non-invertible on $\h$ (see Lemma \ref{lemg}).  
The unconditional convergence follows as in Example \ref{nbnb21}(ii).
  \end{ex}

\begin{ex}\label{bnbnew2}
Let $\Phi=(e_1, \ \ \frac{1}{2}e_2, \ \ \frac{1}{\sqrt{2}}e_2,\ \ \frac{1}{3}e_3, \ \ \frac{1}{3}e_3,\ \ \frac{1}{\sqrt{3}}e_3, 
\ \ \frac{1}{4}e_4, \ \ \frac{1}{4}e_4, \ \ \frac{1}{4}e_4, \ \ \frac{1}{\sqrt{4}}e_4,\ \ldots)$,

  \vspace{.05in} \hspace{.99in}
   $\Psi=(e_1,\  \ \ \ e_2,\ \ \frac{1}{\sqrt{2}}e_2,\ \ \ \, \ e_3, \ \ \ \ \, e_3, \ \ \frac{1}{\sqrt{3}}e_3, \,  \ \ \ \ e_4, \ \, \ \ \ e_4, \,\ \ \ \ e_4, \ \ \frac{1}{\sqrt{4}}e_4,\ \ldots)$.

\vspace{.05in} \noindent
  Then $M_{(1),\Phi,\Psi}=M_{(1),\Psi,\Phi}=I$. The unconditional convergence on $\h$ follows in the same way as in Example \ref{bnb13}(i).

\end{ex}

\begin{ex}\label{bnbnew3}
Let 
$\Phi=(\frac{1}{2}\,e_1, \ \ \frac{1}{\sqrt{2}}\,e_2, \ \ \frac{1}{2^2}\,e_1, \ \ \frac{1}{\sqrt{3}}\,e_3, \ \ \frac{1}{2^3}\,e_1, 
\ \ \frac{1}{\sqrt{4}}\,e_4, \ \ldots)$, 

\vspace{.05in}\hspace{.99in}
  $\Psi=(\ \ \, \, e_1, \ \ \frac{1}{\sqrt{2}}\,e_2, \, \ \ \ \ \ \, e_1, \ \ \frac{1}{\sqrt{3}}\,e_3, \, \ \ \ \ \ \, e_1, \ \ \frac{1}{\sqrt{4}}\,e_4, \ \ldots)$. 
  
    \vspace{.05in}
\noindent  
Then $M_{(1),\Phi,\Psi}=M_{(1),\Psi,\Phi}=G_1$ - non-invertible on $\h$ (see Lemma \ref{lemg}). 
The convergence is unconditional on $\h$, because $M_{(1),\Phi,\Psi}\neweq M_{(1),\Psi,\Phi} \neweq M_{(1),\Theta,\Theta}$, where 
$\Theta$ is the same as in Example \ref{nbnb21}(ii) \,(apply Prop. \ref{bmh}).

  \end{ex}

  \begin{ex}\label{bnbnew2new} 
Let $\Phi=(e_1, \ \ \frac{1}{2}e_2, \ \ \, \frac{1}{\sqrt{2}}e_2,\ \ \frac{1}{3}e_3, \ \ \frac{1}{3}e_3,\ \ \frac{1}{\sqrt{3}}e_3, 
\ \ \frac{1}{4}e_4, \ \ \frac{1}{4}e_4, \ \ \frac{1}{4}e_4, \ \ \frac{1}{\sqrt{4}}e_4,\ \ldots)$,

  \vspace{.05in} \hspace{.99in}
   $\Psi=(e_1,\  \ \ \ e_2,\ \ \frac{1}{\sqrt[4]{2}}e_2,\ \ \ \, \, \ e_3, \ \ \ \  e_3, \ \   \frac{1}{\sqrt[4]{3}}e_3,   \ \ \ \ e_4, \ \ \ \ e_4, \ \ \ \ e_4, \ \ \frac{1}{\sqrt[4]{4}}e_4,\ \ldots)$, 
  
   \vspace{.05in} \hspace{.98in}
$m=(\ \, 1, \ \ \ \ \,  \ 1,\ \ \ \ \ \, \frac{1}{\sqrt[4]{2}}, \ \ \ \ \ \ \, 1,\ \ \ \ \ \, 1,\ \ \ \ \ \ \frac{1}{\sqrt[4]{3}}, 
\ \ \ \ \, \ 1,\ \ \ \ \ \, 1, \ \ \ \ \ \, 1, \ \ \ \ \ \, \frac{1}{\sqrt[4]{4}}, \ \ldots)$.

 \vspace{.05in}
\noindent  Then $M_{m,\Phi,\Psi}=M_{m,\Psi,\Phi}=I$. The unconditional convergence on $\h$ follows in the same way as in Example \ref{bnb13}(i).

\end{ex}

\begin{ex}\label{bnbnew4} 
Let $\Phi=(e_1, \ \ \frac{1}{2}e_2, \ \ \frac{1}{2}e_2,\ \ \frac{1}{3}e_3, \ \ \frac{1}{3}e_3,\ \ \frac{1}{3}e_3, \ \ \frac{1}{4}e_4, \ \ \frac{1}{4}e_4,\ \ \frac{1}{4}e_4, \ \ \frac{1}{4}e_4, \ \ldots)$,

  \vspace{.05in} \hspace{.99in}
   $\Psi=(e_1,\  \ \ \ e_2, \ \ \frac{1}{2}e_2,\ \ \ \, \, e_3, \ \ \ \ \, e_3, \ \ \frac{1}{3}e_3, \ \ \ \ \, e_4, \ \ \ \ \, e_4, \ \ \ \ \, e_4, \ \ \frac{1}{4}e_4, \ \ldots)$.

   \vspace{.1in} \hspace{.54in}
(i) Let \   
$m=(\ \ 1, \ \ \ \ \, \, \frac{1}{2},\ \ \ \ \ \  1,\ \ \ \ \, \, \frac{1}{3},\ \ \ \ \, \ \frac{1}{3},\ \ \ \ \ \ 1, \ \ \ \ \ \,  \frac{1}{4},\ \ \ \ \ \ \frac{1}{4},\ \ \ \ \, \ \frac{1}{4}, \ \ \ \ \  \, \ 1, \ \ldots)$.

\vspace{.05in}
\noindent Then $M_{m,\Phi,\Psi}=M_{m,\Psi,\Phi}=G_1$ - non-invertible on $\h$ (see Lemma \ref{lemg}). 
  The convergence is unconditional on $\h$, because $M_{m,\Phi,\Psi}\neweq M_{m,\Psi,\Phi} \neweq M_{(1),\Theta,\Theta}$, where 
$\Theta$ is the same as in Example \ref{bnb14} (apply Prop. \ref{bmh}).

    \vspace{.1in} \hspace{.51in}
(ii) Let \
 $m=(\ \ 1, \ \ \ \ \  \, 1,\ \ \ \ \, \ \, 2,\ \ \ \ \ \, 1,\ \ \ \ \ \ 1,\ \ \ \ \ \ \, 3, \ \ \ \ \  \, 1,\ \ \ \ \ \, \, 1,\ \ \ \  \, \ \ 1, \ \ \ \ \  \, \ 4, \ \ldots)$.  
  
    \vspace{.05in}
 \noindent  
  Then  $M_{m,\Phi,\Psi}=M_{m,\Psi,\Phi}=I$. 
The convergence  is unconditional on $\h$, because
    $M_{m,\Phi,\Psi}\neweq M_{m,\Psi,\Phi} \neweq M_{(1),\Theta,\Theta}$, where 
$\Theta$ is the same as in Example \ref{bnb13}(i).
  
  \end{ex}

\begin{ex}\label{bnbnew5} 
Let $\Phi=(e_1, \ \ \frac{1}{2}e_2, \ \ \ \, \frac{1}{2}e_2,\ \ \frac{1}{3}e_3, \ \ \frac{1}{3}e_3,\ \ \, \, \frac{1}{3}e_3, \ \ \frac{1}{4}e_4, \ \ \frac{1}{4}e_4,\ \ \frac{1}{4}e_4, \ \ \ \, \frac{1}{4}e_4, \ \ldots)$,

  \vspace{.05in} \hspace{.99in}
   $\Psi=(e_1,\  \ \ \ e_2, \ \ \frac{1}{2^2}e_2,\ \ \ \, \ e_3, \ \ \ \ \,  e_3, \ \ \frac{1}{3^2}e_3, \ \ \ \ e_4, \ \ \ \ \ e_4, \ \ \ \ \, e_4, \ \ \frac{1}{4^2}e_4, \ \ldots)$. 
  
   \vspace{.05in} \hspace{.98in}
$m=(\  1,\, \ \ \ \ \,  \frac{1}{2},\ \ \ \ \ \ \ \, \, 2,\ \ \ \ \, \ \frac{1}{3},\ \ \ \ \, \ \frac{1}{3},\ \ \ \ \ \ \, \, 3, \ \ \ \ \ \, \frac{1}{4},\ \ \ \ \ \, \frac{1}{4},\ \ \ \ \, \, \, \frac{1}{4}, \ \ \ \ \ \ \, \ 4, \ \ldots)$.

\vspace{.05in}
\noindent Then $M_{m,\Phi,\Psi}=M_{m,\Psi,\Phi}=G_1$ - non-invertible on $\h$ (see Lemma \ref{lemg}). 
  The convergence is unconditional on $\h$, because $M_{m,\Phi,\Psi}\neweq M_{m,\Psi,\Phi} \neweq M_{(1),\Theta,\Theta}$, where 
$\Theta$ is the same as in Example \ref{bnb14} (apply Prop. \ref{bmh}).

  \end{ex}

\begin{ex}\label{bnb19}
Let $\Phi=(\frac{1}{n}e_n)$ and $\Psi=(ne_n)$. 
Then  $M_{(1),\Phi,\Psi}\neweq M_{(1),\Psi,\Phi} \neweq M_{(1),(e_n),(e_n)}=I$. 
\end{ex}

\begin{ex}\label{bnb20}
 Let $\Phi=(\frac{1}{n^2}e_n)$ and $\Psi=(ne_n)$. 
Then  $M_{(1),\Phi,\Psi}\neweq M_{(1),\Psi,\Phi}\neweq M_{(\frac{1}{n}),(e_n),(e_n)} = G_1$
- unconditionally convergent and non-invertible on $\h$ (see Lemma \ref{lemg}).
\end{ex}

\begin{ex}\label{bnb21i} Let $\Phi=(\frac{1}{n}e_n)$, $\Psi=(n^2e_n)$, and $m=(\frac{1}{n})$. 
Then  $M_{m,\Phi, \Psi}\neweq M_{m,\Psi,\Phi} \neweq M_{(1),(e_n),(e_n)}=I$.
\end{ex}

\begin{ex}\label{bnb21ii} 
Let $\Phi=(\frac{1}{n}e_n)$, $\Psi=(n^2e_n)$, and $m=(\frac{1}{n^2})$. 
Then $M_{m,\Phi, \Psi}\neweq M_{m,\Psi,\Phi}  \neweq M_{(\frac{1}{n}),(e_n),(e_n)}=G_1$
- unconditionally convergent and non-invertible on $\h$ (see Lemma \ref{lemg}).
\end{ex}

\begin{ex}\label{bnb23i}
Let $\Phi=(\frac{1}{n^3}e_n)$, $\Psi=(ne_n)$, and $m=(n^2)$. 
Then  $M_{m,\Phi,\Psi}\neweq M_{m,\Psi,\Phi}\neweq M_{(1),(e_n),(e_n)} =I$. 
\end{ex}

\begin{ex}\label{bnb23ii} 
Let $\Phi=(\frac{1}{n^3}e_n)$, $\Psi=(ne_n)$, and $m=(n)$. 
Then  $M_{m,\Phi,\Psi}\neweq M_{m,\Psi,\Phi} \neweq M_{(\frac{1}{n}),(e_n),(e_n)}=G_1$
- unconditionally convergent and non-invertible on $\h$ (see Lemma \ref{lemg}).
\end{ex}

\begin{ex}\label{bnbnew6}
Let 
$\Phi=(\frac{1}{\sqrt{2}}e_1, \ \ \  \frac{1}{2}e_2, \ \  \frac{1}{\sqrt{2^2}}e_1, \ \ \ \frac{1}{3}e_3, \ \ \, \frac{1}{\sqrt{2^3}}e_1, \ \ \ \frac{1}{4}e_4,\ \ldots)$, 

\vspace{.05in} \hspace{.99in}
$\Psi=(\frac{1}{\sqrt{2}} e_1,\ \ \ 2 e_2, \ \ \, \frac{1}{\sqrt{2^2}} e_1, \ \ \ 3 e_3, \ \ \, \frac{1}{\sqrt{2^3}} e_1,\ \ \ \, 4 e_4,\ \ldots)$.

\vspace{.05in}
\noindent Then $M_{(1),\Phi,\Psi}=M_{(1),\Psi,\Phi}=I$. The convergence is unconditional on $\h$, because $M_{(1),\Phi,\Psi} \neweq M_{(1),\Psi,\Phi} \neweq M_{(1),\Theta,\Theta}$, where
$\Theta$ is the same as in Example \ref{nbnb21}(i).

  \end{ex}

\begin{ex}\label{bnbnew7}
Let 
$\Phi=(\frac{1}{\sqrt{2}}e_1, \ \   \frac{1}{2^2}e_2, \ \  \frac{1}{\sqrt{2^2}}e_1, \  \ \frac{1}{3^2}e_3, \ \  \frac{1}{\sqrt{2^3}}e_1,  \ \ \frac{1}{4^2}e_4,\ \ldots)$, 

\vspace{.05in} \hspace{.99in}
$\Psi=(\frac{1}{\sqrt{2}} e_1,\ \ \ \, 2 e_2, \ \ \, \frac{1}{\sqrt{2^2}} e_1, \ \ \ \, 3 e_3, \ \ \, \frac{1}{\sqrt{2^3}} e_1,\ \ \ \, 4 e_4,\ \ldots)$.

\vspace{.05in}
\noindent Then $M_{(1),\Phi,\Psi}=M_{(1),\Psi,\Phi}=G_1$ - non-invertible on $\h$ (see Lemma \ref{lemg}). The convergence is unconditional on $\h$, because $M_{(1),\Phi,\Psi} \neweq M_{(1),\Psi,\Phi} \neweq M_{(1),\Theta,\Theta}$, where
$\Theta$ is the same as in Example \ref{nbnb21}(ii) (apply Prop. \ref{bmh}).
  
  \end{ex}

  \begin{ex}\label{bnbnew8}
  Let 
$\Phi=(\frac{1}{\sqrt{2}}e_1, \ \ \ \, \frac{1}{2}e_2, \ \ \, \frac{1}{\sqrt{2^2}}e_1, \ \ \ \, \frac{1}{3}e_3, \ \ \, \, \frac{1}{\sqrt{2^3}}e_1, \ \ \ \, \frac{1}{4}e_4,\ \ldots)$, 

\vspace{.05in} \hspace{.99in}
$\Psi=(\frac{1}{\sqrt{2}} e_1,\ \  2^2 e_2, \ \ \, \frac{1}{\sqrt{2^2}} e_1, \ \ \, 3^2 e_3, \ \ \, \frac{1}{\sqrt{2^3}} e_1,\ \ \, 4^2 e_4,\ \ldots)$.

\vspace{.1in} \hspace{.59in}
(i) Let
$m=(\ \ \  \ \ \ \, 1,\ \ \ \ \ \, \, \ \frac{1}{2}, \ \ \ \ \ \ \ \ \ \  \, 1, \ \ \ \ \ \ \ \, \frac{1}{3}, \ \ \ \ \ \ \ \ \ \ \  1, \ \ \ \ \ \ \ \,   \frac{1}{4},\ \ldots)$.

\vspace{.05in}
\noindent Then $M_{m,\Phi,\Psi}=M_{m,\Psi,\Phi}=I$. The unconditional convergence follows as in Example \ref{nbnb21}(i).

 \vspace{.1in} \hspace{.555in}
(ii) Let $m=(\ \ \  \ \ \ \, 1,\ \ \ \  \, \ \frac{1}{2^2}, \ \ \ \ \ \ \ \ \ \  \, 1, \ \ \ \ \ \ \frac{1}{3^2}, \ \ \ \ \ \ \ \ \ \ \ 1, \ \ \ \ \ \  \, \frac{1}{4^2},\ \ldots)$. 

 \vspace{.05in}
\noindent Then $M_{m,\Phi,\Psi}=M_{m,\Psi,\Phi}=G_1$ - non-invertible on $\h$ (see Lemma \ref{lemg}). 
The convergence is unconditional on $\h$, because $M_{m,\Phi,\Psi}\neweq M_{m,\Psi,\Phi}\neweq M_{(1),m\Psi,\Phi}$ and the sequences $\Phi$ and $m\Psi$  are Bessel for $\h$ (apply Prop. \ref{bmh}).
  \end{ex}

    \begin{ex}\label{bnbnew10}
      Let
  $\Phi= ( e_1, \ \ \, \frac{1}{2} e_2,\ \ \frac{1}{3}e_3,\ \ \frac{1}{4}  e_4, \ \ \frac{1}{5}e_5, \ \ldots)$, 
 
 \vspace{.05in}\hspace{.99in}
 $\Psi=(e_1, \ \ 2\,e_2, \ \ \frac{1}{3}e_3, \ \ 4\, e_4, \ \ \frac{1}{5}e_5, \ \ldots)$,

   \vspace{.1in} \hspace{.55in}
(i) Let \ 
  $m=(\ \, 1, \ \ \ \ \ \ \, 1,\ \ \ \, \, 3^2,\ \ \ \ \ \ 1,\ \ \ \  \, 5^2, \  \ldots)$.

  \vspace{.05in} \noindent
Then  $M_{m,\Phi,\Psi}\neweq M_{m,\Psi,\Phi}\neweq M_{(1),(e_n), (e_n)}=I$.

     \vspace{.1in} \hspace{.52in}
(ii) Let \
   $m=(\ \, 1, \ \ \ \ \ \   \frac{1}{2},\ \ \ \ \ \  3,\ \ \ \ \ \, \frac{1}{4},\ \ \ \ \ \, \, 5, \  \ldots)$.
      
    \vspace{.05in} \noindent
 Then $M_{m,\Phi,\Psi}\neweq M_{m,\Psi,\Phi}\neweq M_{(\frac{1}{n}),(e_n), (e_n)}=G_1$  
- unconditionally convergent and non-invertible on $\h$ (see Lemma \ref{lemg}).
   \end{ex}

\subsection{Examples for two Bessel non-frame sequences; TABLE 3 
on page \pageref{table3}}

\begin{ex}\label{bb12} Let $\Phi=(e_2, e_3, e_4, e_5, \ldots)$.

 \vspace{.05in}
 (i) Let $m=(1,1,1,1,\ldots)$. 
  Then $M_{m,\Phi,\Phi}$ is clearly unconditionally convergent on $\h$ and not surjective.
 
  \vspace{.05in}
 (ii) Let $m=(\frac{1}{2}, \frac{1}{3}, \frac{1}{4}, \frac{1}{5}, \ldots)$. 
  Then $M_{m,\Phi,\Phi}$ is clearly unconditionally convergent on $\h$ and not surjective.
 
\end{ex}

\begin{ex}\label{bb34}
Let $\Phi=(e_2, e_3, e_4, e_5, \ldots)$ and $\Psi=(\frac{1}{2}e_2,\frac{1}{3}e_3,\frac{1}{4}e_4,\frac{1}{5}e_5,\ldots)$.

 \vspace{.05in}
 (i) Let $m=(1,1,1,1,\ldots)$. 
  Then $M_{m,\Phi,\Psi}$ and $M_{m,\Psi,\Phi}$ are clearly unconditionally convergent on $\h$ and not surjective.
 
  \vspace{.05in}
 (ii) Let $m=(\frac{1}{2}, \frac{1}{3}, \frac{1}{4}, \frac{1}{5}, \ldots)$. 
   Then $M_{m,\Phi,\Psi}$ and $M_{m,\Psi,\Phi}$ are clearly unconditionally convergent on $\h$ and not surjective.
 
 \vspace{.05in}
 (iii) Let $m=(2,3,4,5,\ldots)$. 
  Then $M_{m,\Phi,\Psi}$ and $M_{m,\Psi,\Phi}$ are clearly unconditionally convergent on $\h$ and not surjective.
 \end{ex}

\begin{ex}\label{bb7}
Let $\Phi=(\frac{1}{n} e_n)$. 

\vspace{.05in}
(i) Let $m=(1)$. 
Then $M_{m,\Phi,\Phi}\neweq M_{(\frac{1}{n^2}),(e_n),(e_n)} =G_2$
- unconditionally convergent and non-invertible on $\h$ (see Lemma \ref{lemg}).

\vspace{.05in}
(ii) Let $m=(\frac{1}{n})$. 
Then $M_{m,\Phi,\Phi}\neweq M_{(\frac{1}{n^3}),(e_n),(e_n)}=G_3$
- unconditionally convergent and non-invertible on $\h$ (see Lemma \ref{lemg}).

\vspace{.05in}
(iii) Let $m=(n^2)$. 
Then $M_{m,\Phi,\Phi}\neweq M_{(1),(e_n),(e_n)}=I$.

\vspace{.05in}
(iv) Let $m=(n)$. 
Then $M_{m,\Phi,\Phi}\neweq M_{(\frac{1}{n}),(e_n),(e_n)}=G_1$
- unconditionally convergent and non-invertible on $\h$ (see Lemma \ref{lemg}).
\end{ex}

\subsection{Examples for one overcomplete frame and one non-Bessel sequence; TABLE 4 
on page \pageref{table4}}

\begin{ex} \label{ex11table4}
Let
$\Phi=(e_1, \ \ e_2, \ \ \ \,\, e_2, \ \ e_2, \ \  e_3, \ \ \ \, \, e_3, \  \ e_3, \  \ e_4, \ \ \ \, \, e_4, \ \ e_4,\ \ldots)$,

\vspace{.05in} \hspace{.91in}
$\Psi=(e_1, \  \ e_1, \ - e_1,  \ \ e_2,  \  \ e_1,   \ - e_1,  \ \ \, e_3, \  \ e_1,    \ - e_1,  \ \ e_4, \ \ldots)$,

\vspace{.05in} \hspace{.9in}
$m=(\ \, 1, \ \ \ \frac{1}{2}, \ \ \ \  \ \frac{1}{2}, \ \  \, \ 1, \ \ \, \frac{1}{2^2}, \ \  \, \frac{1}{2^2}, \ \ \ \,\, 1, \ \ \, \frac{1}{2^3}, \ \ \ \frac{1}{2^3}, \ \ \ \, 1, \ \ldots)$.

\vspace{.05in}
\noindent 
Then $M_{m,\Phi,\Psi}=M_{m,\Psi,\Phi}=I$. 
The convergence of $M_{m,\Phi,\Psi}$ and $M_{m,\Psi,\Phi}$ is unconditional on $\h$,  because $\Phi$ and $m\Psi$ are Bessel for $\h$ (apply Prop. \ref{bmh}).  
\end{ex}

\begin{ex}\label{fnb2} Let
$\Phi=(e_1, \ \ \ e_1, \ \ \ e_2, \ \ \ e_3, \ \ \ e_4, \ \ \ e_5, \ \ldots)$,

\vspace{.05in} \hspace{.91in}
$\Psi=(e_2, \ \ \ e_2,  \ \ \ e_3, \ \ \ e_2, \ \ \  e_4, \ \ \ e_2, \ \ldots)$,

\vspace{.05in} \hspace{.91in}
$m=(\ \, 1, \ \ \ \ \frac{1}{2}, \ \ \ \, \, 1, \ \ \  \frac{1}{2^2}, \ \ \ \,\, 1, \ \ \, \frac{1}{2^3}, \ \ldots)$.

\vspace{.05in}
\noindent 
 The unconditional convergence of $M_{m,\Phi,\Psi}$ and $M_{m,\Psi,\Phi}$ follows as in Example \ref{ex11table4}. The multiplier $M_{m,\Phi,\Psi}$ is not injective, for example $M_{m,\Phi,\Psi}e_1 =0$.
 The multiplier $M_{m,\Psi,\Phi}$ is not surjective, 
 because $\Psi$ is not complete in $\h$.

\end{ex}

\begin{ex} \label{fnb51} 
Let 
$\Phi=( e_1, \ \ \ \  e_1, \ \ \  \, \, e_1, \ \ \ \ \, \ e_1, \ \ \ \ \ \ \ \, e_1,  
\ \ \ \ \  e_2, \ \  \ \ \ e_2, \ \ \ \ \, \, e_2,  \ \ \ \ \, \ e_2, \ \ \ \ \ \ \ \, e_2, \ \ \ \ \  e_3, \ \   \ \ \ e_3, \ \ \ \ \, \, e_3, \ \ \ \ \, \, e_3, \ \ \ \ \ \ \, \ \, e_3, \ \ldots)$,

\vspace{.05in} \hspace{.91in}
$\Psi=( e_1, \ \ \ \ e_1,  \ - e_1, \ \, \frac{1}{2^2}  e_1, \,  - \frac{1}{2^2}e_1,  
\ \ \ \ \  e_2, \ \  \ \ \ e_1, \ \ - e_1, \ \, \frac{1}{2^4}  e_1, \, -\frac{1}{2^4} e_1, \ \ \ \ \  e_3, \ \  \ \ \ e_1, \ \ - e_1,  \ \, \frac{1}{2^6} e_1, \  - \frac{1}{2^6} e_1,  \ \ldots)$,

\vspace{.05in} \hspace{.9in}
$m=(\ \, 1,\ \ \ \ \ \frac{1}{2}, \ \ \ \ \ \frac{1}{2},\ \ \ \ \ \ \ \, 1,\ \ \ \ \ \ \ \ \ \, 1,\ \ \ \  \ \ \  1, \ \ \ \, \,  \frac{1}{2^3}, 
\ \ \ \ \, \frac{1}{2^3}, \ \ \ \ \ \, \ \, 1, \ \ \ \ \ \ \ \ \ \, 1, \ \ \ \ \ \  \, 1, \ \  \ \, \,  \frac{1}{2^5}, 
 \ \ \ \ \ \frac{1}{2^5},\ \ \ \ \ \ \,\, 1,\ \ \ \ \ \ \ \ \ \, 1,   \ \ldots)$.
 
 \vspace{.05in} \noindent
Then $M_{m,\Phi,\Psi}=M_{m,\Psi,\Phi}=I$. 
The convergence is unconditional on $\h$, because $\Phi$ and $m\Psi$ are Bessel for $\h$ (apply Prop. \ref{bmh}).
\end{ex}

\begin{ex} \label{fnb52}
Let 
$\Phi=( e_1,  \ \ \ e_1, \ \  \ e_2, \ \ \  e_3, \ \ \ \,  e_4,  
\ \ \ \, \ e_5,  \ \ \ e_6,  \ \ \ \ e_7, \ \ \ \, e_8, \ \ \ \ e_9, \ \ \,  e_{10}, \ \ldots)$,

\vspace{.05in} \hspace{.91in}
$\Psi=( e_2,  \ \ \ e_3, \ \ \ e_1, \ \ \  e_1, \ \ \ \,  e_4,  
 \ \ \frac{1}{5}  e_5, \ \ \ e_4, \ \ \frac{1}{6} e_6, \ \ \ e_4, \ \ \frac{1}{7} e_7, \ \ \ 
 \, e_4,  \ \ldots)$.
 
\vspace{.1in} \hspace{.47in}
(i) Let \ $m=(\ \, 1,\ \ \ \ \, 1, \ \ \ \ \, 1, \ \ \ \ \, 1,  \ \ \ \ \, \frac{1}{2}, 
\ \ \ \ \ \ 1, \ \ \   \frac{1}{2^2}, \ \ \ \ \ \ 1, \ \ \,  \frac{1}{2^3}, \ \ \ \ \ \, \ 1, \ \ \   \frac{1}{2^4}, \ \ldots)$.

 \vspace{.05in} \noindent
The unconditional convergence of $M_{m,\Phi,\Psi}$ and $M_{m,\Psi,\Phi}$ follows as in Example \ref{ex11table4}. The multiplier $M_{m,\Phi,\Psi}$ (resp. $M_{m,\Psi,\Phi}$) is not injective, because $M_{m,\Phi,\Psi}e_2=e_1=M_{m,\Phi,\Psi}e_3$ (resp. $M_{m,\Psi,\Phi}e_2 =e_1=M_{m,\Psi,\Phi}e_3$).

\vspace{.1in} \hspace{.44in}
(ii)  Let \ $m=(\ \, 1,\ \ \ \ \, 1, \ \ \ \ \, 1, \ \ \ \ \, 1,  \ \ \ \ \, \frac{1}{2}, 
\ \ \ \ \ \ 5, \ \ \   \frac{1}{2^2}, \ \ \ \ \ \ 6, \ \ \,  \frac{1}{2^3}, \ \ \ \ \ \, \ 7, \ \ \   \frac{1}{2^4}, \ \ldots)$.

 \vspace{.05in} \noindent
Then the same conclusion as in (i) holds.

\end{ex}

\begin{ex} \label{fnb5b}
Let
$\Phi=(e_1, \ \ e_2, \ \ \ \ e_2, \ \ \ \ \, e_2, \ \  e_3, \ \ \ \ e_3, \ \, \ \ \ e_3, \ \, \ e_4, \ \ \ \, e_4, \ \ \ \ \, e_4,\ \ldots)$,

\vspace{.05in} \hspace{.91in}
$\Psi=(e_1, \  \ e_1, \ - e_1,  \ \ \frac{1}{2} e_2,  \  \ e_1,   \ - e_1,  \ \ \frac{1}{3} e_3, \ \, \ e_1,    \ - e_1,  \ \ \frac{1}{4} e_4, \ \ldots)$,

\vspace{.05in} \hspace{.9in}
$m=(\ \, 1, \ \ \  \frac{1}{2}, \ \ \ \ \,\, \frac{1}{2}, \ \ \ \ \ \ 2, \ \ \, \frac{1}{2^2}, \ \  \ \frac{1}{2^2}, \ \ \ \ \ \ 3, \ \  \frac{1}{2^3}, \ \ \ \, \frac{1}{2^3}, \ \ \ \ \ \ 4, \ \ldots)$.

\vspace{.05in}
\noindent 
Then $M_{m,\Phi,\Psi}=M_{m,\Psi,\Phi}=I$. 
The convergence is unconditional on $\h$, because $\Phi$ and $m\Psi$ are Bessel for $\h$ (apply Prop. \ref{bmh}).

\end{ex}

\begin{ex}\label{fnb3}
Let $\Phi=(\ \ \,\, e_1, \ \ \ \ \,\, e_1, \ \ \ \ \, e_2, \ \ \ \ \, e_3, \ \ \ \ \, e_4, \ \ \ \ \,e_5, \ \ldots)$, 

\vspace{.05in} \hspace{.91in}
$\Psi=(\frac{1}{2}\,e_1, \ \ \frac{1}{2}\,e_1, \ \ \, 2e_2, \ \ \,  3e_3, \ \ \, 4e_4,\ \ \, 5e_5, \ \ldots)$.

\vspace{.1in} \hspace{.48in}
(i) Let \ $m=(\ \ \ \,\,\, 1, \ \ \ \ \ \ \, 1, \, \ \ \ \ \ \frac{1}{2}, \,\ \ \ \ \ \frac{1}{3}, \ \ \ \ \ \ \frac{1}{4}, \ \ \ \ \ \, \frac{1}{5},\  \ldots)$.

\vspace{.05in} \noindent 
Then $M_{m,\Phi,\Psi}=M_{m,\Psi,\Phi}=I$. 
The convergence is unconditional on $\h$, because $\Phi$ and $m\Psi$ are Bessel for $\h$ (apply Prop. \ref{bmh}).

\vspace{.1in} \hspace{.44in}
(ii) Let \ $m=(\ \ \ \,\,\, 1, \ \ \ \ \ \ \, 1,\, \ \ \ \  \frac{1}{2^2}, \ \ \ \, \frac{1}{3^2}, \ \ \  \ \, \frac{1}{4^2}, \ \ \ \  \frac{1}{5^2},\  \ldots)$.

\vspace{.05in} \noindent 
Then $M_{m,\Phi,\Psi}=M_{m,\Psi,\Phi}=G_1$ 
- non-invertible on $\h$ (see Lemma \ref{lemg}). 
The convergence is unconditional on $\h$, because $\Phi$ and $m\Psi$ are Bessel for $\h$ (apply Prop. \ref{bmh}).
\end{ex}

\begin{ex}\label{fnb31}
Let 
$\Phi=( e_1, \ \ \ \  e_1, \ \ \  \, e_1, \ \ \ \ \  e_2, \ \  \ \ \ \,\, e_2, \ \ \ \ \,\,  e_2,  \ \ \ \ \ \, e_3, \ \ \,  \ \ \ \, e_3, \ \ \ \ \, e_3, \ \ \ \ \, e_4, \ \ \ \ \ \ \, e_4, \ \ \ \ \ \ \, e_4, \ \ \ \ \, e_5, \ \ \ \ \ \ \, e_5, \ \ \ \ \,  e_5, \ \ldots)$,

\vspace{.05in} \hspace{.91in}
$\Psi=( e_1, \ \ \ \ e_1,  \ \ \ \, e_1, 
\ \ \ \ \  e_2, \ \  \ \ \frac{1}{2}e_2, \ \ \, \frac{1}{2}e_2,  \ \ \ \ \ \, e_3, \ \  \ \   3e_3, \ \ \ 3 e_3,   \ \ \ \ \, e_4, \ \  \ \   \frac{1}{4}e_4, \ \ \ \  \frac{1}{4} e_4,  \ \ \ \  \, e_5, \ \  \ \  \, 5e_5, \ \ \ 5 e_5,\ \ldots)$,

\vspace{.05in} \hspace{.9in}
$m=(\ \, 1,\ \ \ \ \ \, 1, \  \  -1, \ \ \ \ \ \ \, 1, \ \ \  \ \ \ \ \,\, 1, \ \ \ \, -1, \ \ \ \ \ \ \  1, 
\ \ \ \ \ \, \frac{1}{3^2}, \ \   -\frac{1}{3^2}, \ \ \ \ \ \,  1, \ \ \ \ \ \ \ \ \, 1, \ \  \ \ \ -1, \ \ \ \ \ \,\,  1, 
\ \ \ \ \ \  \frac{1}{5^2}, \ \   -\frac{1}{5^2},
   \ \ldots)$.
 
 \vspace{.05in} \noindent
Then 
$M_{m,\Phi,\Psi}=M_{m,\Psi,\Phi}=I$. 
The convergence is unconditional on $\h$, because $\Phi$ and $m\Psi$ are Bessel for $\h$ (apply Prop. \ref{bmh}).
\end{ex}

\begin{ex}\label{fnb32}
Let 
$\Phi=( e_1, \ \ \ \  e_1, \ \ \ \ \ \, e_2, \ \  \ \ \ \ \ e_2,   \ \ \ \ \ \ e_3, \ \   \ \ \ e_3, \ \ \ \ \ \ \, e_4, \ \ \ \ \ \, e_4,  \ \ \ \ \ \, e_5, \ \   \ \ \ \ \ \, e_5, \ \ \ \ \ \, e_6, \ \ \ \ \ \ \ e_6,  \ \ldots)$,

\vspace{.05in} \hspace{.91in}
$\Psi=( e_2, \ \ \ \ e_2,  \ \ \ \ 2e_2, 
\ \ \ \ \frac{1}{2}  e_2, \ \  \ \ 3e_3, \ \ \, \frac{1}{3}e_3,   \ \  \ \  \, 4e_4, \ \ \ \frac{1}{4} e_4,   \ \ \ \ 5e_5, \ \  \ \  \, \frac{1}{5}e_5, \ \ \ \  6e_6,  \ \ \ \ \, \frac{1}{6}e_6, \ \ldots)$,

\vspace{.05in} \hspace{.9in}
$m=(\ \, 1,\ \ \ \ \ \, 1, \ \ \ \ \ \ \  \frac{1}{2}, \ \ \ \ \ \ \ \ 1, \ \ \  \ \ \ \, \, \frac{1}{3}, \ \ \ \ \ \ \,  1, \ \ \ \ \ \ \ \   \frac{1}{4}, \ \ \ \ \ \ \ 1, \ \ \ \ \ \ \ \,  \frac{1}{5},\, \ \ \ \ \ \ \ \ 1,\ \ \ \ \ \ \ \frac{1}{6}, \ \ \  \ \ \ \ \ \ 1, \ \ldots)$.
   
 \vspace{.05in} \noindent 
  Then  $M_{m,\Phi,\Psi}$ and $M_{m,\Psi,\Phi}$ are unconditionally convergent on $\h$, because $\Phi$ and $m\Psi$ are Bessel for $\h$ (apply Prop. \ref{bmh}). They are not injective - for example, $M_{m,\Phi,\Psi}e_1=0$ and $M_{m,\Psi,\Phi}(\frac{1}{2}e_1)=e_2=M_{m,\Psi,\Phi}(\frac{2}{3}e_2)$.
\end{ex}

\begin{ex}\label{fnb5}
Let 
$\Phi=(\ \ e_1, \ \ \ \ \, e_1, \ \ \  \ e_2, \ \ \ \ \ \, e_3, \ \ \ \ e_4, \ \ \ \ \ \, e_5, \ \ldots)$,

\vspace{.05in} \hspace{.9in}
$\Psi=(\frac{1}{2}e_1, \ \ \frac{1}{2}e_1, \ \ 2e_2, \ \ \ \frac{1}{3}e_3,\ \  4e_4, \ \ \  \frac{1}{5}e_5,\ \ldots)$,

\vspace{.05in} \hspace{.9in}
$m=(\ \ \ \ 1,\ \ \ \ \ \ 1, \ \ \ \ \ \frac{1}{2}, \ \ \ \ \ \ \ \, 3, \ \ \ \ \,  \frac{1}{4},\ \ \  \ \ \ \ \, 5,\  \ldots)$. 

 \vspace{.05in} \noindent 
Then 
$M_{m,\Phi,\Psi}=M_{m,\Psi,\Phi}=I$. 
The convergence is unconditional on $\h$, because $\Phi$ and $m\Psi$ are Bessel  for $\h$ (apply Prop. \ref{bmh}). 
\end{ex}

\begin{ex}\label{fnb6}
Let 
$\Phi=(\ \, \ e_1, \ \ \ \ \, e_1, \ \ \  \ e_2, \ \ \ \ \ \ e_3, \ \ \ \ e_4, \ \ \ \ \ \ \, e_5, \ \ldots)$,

\vspace{.05in} \hspace{.99in}
$\Psi=(\frac{1}{2}e_1, \ \ \frac{1}{2}e_1, \ \ 2e_2, \ \ \frac{1}{3^2}e_3,\ \  4e_4, \ \ \frac{1}{5^2}e_5,\ \ldots)$,

\vspace{.05in} \hspace{.98in}
$m=(\ \ \ \ 1,\ \ \ \ \ \  1, \ \ \ \ \frac{1}{2^2}, \ \ \ \ \ \ \ \ 3, \ \ \ \, \frac{1}{4^2},\ \ \ \ \ \ \ \, 5,\  \ldots)$.

\vspace{.05in} \noindent 
Then 
$M_{m,\Phi,\Psi}=M_{m,\Psi,\Phi}=G_1$ - non-invertible on $\h$ (see Lemma \ref{lemg}). 
 The convergence is unconditional on $\h$, because $\Phi$ and $m\Psi$ are Bessel for $\h$ (apply Prop. \ref{bmh}).

\end{ex}

\begin{ex}\label{fnb7}
Let $\Phi=(\frac{1}{2}\,e_1, \ \ \ e_2, \ \  \frac{1}{2^2}\,e_1, \ \ \ \, e_3, \ \  \frac{1}{2^3}\,e_1, \ \ \ e_4, \ \ldots)$, 

\vspace{.05in} \hspace{.99in}
$\Psi=(\ \ \, e_1, \ \ \ e_2, \ \ \ \ \ \ \ e_1, \ \ \ \, e_3, \ \ \ \ \ \ \, e_1, \ \ \ e_4, \ \ldots)$.

\vspace{.1in} \hspace{.56in}
(i) Then $M_{(1),\Phi,\Psi}=M_{((1),\Psi,\Phi}=I$.
The convergence is unconditional on $\h$, because $M_{(1),\Phi,\Psi}\neweq M_{(1),\Psi,\Phi} \neweq M_{(1),\Theta,\Theta}$, where 
$\Theta$ is the same as in Example \ref{nbnb21}(i) (apply Prop. \ref{bmh}).

\vspace{.1in} \hspace{.55in}
(ii) Let 
$m=(\ \ \ \ 1, \ \ \ \ \frac{1}{2},\ \ \ \ \ \ \ \ \ 1,\ \ \ \ \, \frac{1}{3}, \ \ \ \ \  \ \ \ \, 1, \ \ \ \  \frac{1}{4}, \ \ldots)$. 
 
 \vspace{.05in} \noindent
Then $M_{m,\Phi,\Psi}=M_{m,\Psi,\Phi}=G_1$ - non-invertible on $\h$ (see Lemma \ref{lemg}).
 The convergence is unconditional on $\h$, because
$M_{m,\Phi,\Psi}\neweq M_{m,\Psi,\Phi} \neweq M_{m,\Theta,\Theta}$, where $\Theta$ is the same as in Example \ref{nbnb21}(i) (apply Prop. \ref{bmh}).

\end{ex}

\begin{ex}\label{fnb8}
Let 
$\Phi=(\frac{1}{2}e_1, \ \ \ e_2, \ \ \  \frac{1}{2^2}e_1, \ \ \ e_3, \ \   \frac{1}{2^3}e_1, \ \ \ e_4, \ \   \frac{1}{2^4}e_1, \ \ \ e_5, \ \ldots)$, 

\vspace{.05in} \hspace{.99in}
$\Psi=(\ \ \, e_3, \ \ \  e_3, \ \ \ \ \ \ \, e_1, \ \ \  e_3, \ \ \ \ \ \ e_1, \ \ \ e_4, \ \ \ \ \ \ e_1, \ \ \ e_5, \ \ldots)$.

\vspace{.05in} 
\noindent
 Then $M_{(1),\Phi,\Psi}$ and $M_{(1),\Psi,\Phi}$ are unconditionally convergent on $\h$, because 
$M_{(1),\Phi,\Psi}\neweq  M_{(1),\Theta,\Xi}$ and $M_{m,\Psi,\Phi} \neweq M_{(1),\Xi,\Theta}$, where

\vspace{.05in} 
$\Theta=(\frac{1}{2}e_1, \ \  e_2, \  \frac{1}{\sqrt{2^2}}e_1, \  \, e_3, \ \,   \frac{1}{\sqrt{2^3}}e_1, \ \ e_4, \    \frac{1}{\sqrt{2^4}}e_1, \ \, e_5,  \ \ldots)$ and

\vspace{.08in} 
$\Xi\, =(\ \  e_3, \ \ \, e_3, \  \frac{1}{\sqrt{2^2}}e_1, \  \, e_3, \ \,   \frac{1}{\sqrt{2^3}}e_1, \  \ e_4, \    \frac{1}{\sqrt{2^4}}e_1, \ \, e_5, \ \ldots)$ are Bessel for $\h$ (apply Prop. \ref{bmh}).

\vspace{.05in} \noindent 
The multiplier $M_{(1),\Phi,\Psi}$ is not injective, for example $M_{(1),\Phi,\Psi}e_2=0$.
The multiplier $M_{(1),\Psi,\Phi}$ is not surjective, 
because $\Psi$ is not complete in $\h$.
\end{ex}

\begin{ex}\label{fnb9}
Let 
$\Phi=(e_1, \ \frac{1}{\sqrt{2}}e_2, \ \frac{1}{\sqrt{2}}e_2, \ \frac{1}{\sqrt{3}}e_3, \ \, \frac{1}{\sqrt{3}}e_3, \ \frac{1}{\sqrt{3}}e_3, \ \ldots)$, 

\vspace{.05in} \hspace{.99in}
$\Psi=(e_1, \ \ \ \ \ e_2, \ \ \ \ \ \ \ e_2, \ \ \ \ \ \ e_3, \ \ \ \ \  \ \ e_3, \ \ \ \ \ \ e_3, \ \ldots)$, 

\vspace{.05in} \hspace{.98in}
$m=(\ \, 1, \ \ \ \, \frac{1}{\sqrt{2}},\ \ \ \, \ \ \frac{1}{\sqrt{2}},\ \ \ \ \frac{1}{\sqrt{3}}, \ \ \ \  \ \ \frac{1}{\sqrt{3}}, \ \ \ \  \frac{1}{\sqrt{3}}, \ \ldots)$.

 \vspace{.05in} \noindent
  Then $M_{m,\Phi,\Psi}=M_{m,\Psi,\Phi}=I$. The convergence is unconditional on $\h$, because
$\Phi$ and $m\Psi$ are Bessel for $\h$ (apply Prop. \ref{bmh}).
\end{ex}

\begin{ex}\label{fnb11}
Let 
$\Phi=(\frac{1}{2^2}e_1, \ \ \ e_2, \ \  \frac{1}{2^4}e_1, \ \ \ e_3, \ \  \frac{1}{2^6}e_1, \ \ \ e_4,\ \ldots)$, 

\vspace{.05in} \hspace{.99in}
$\Psi=(\ \ \ \, e_1,\ \  \ e_2, \ \ \ \ \ \ e_1, \ \  \ e_3, \ \  \ \ \ \ e_1,\ \   \ e_4,\ \ldots)$.

\vspace{.1in}
\hspace{.583in} (i) Let 
$m=(\ \ \ \,\ \ 2, \ \ \  \ \, 1,  \ \ \ \ \,\, 2^2, \ \ \ \  \ 1, \ \  \ \ \ \ 2^3,  \ \ \ \  \, 1,\  \ldots)$.

\vspace{.05in}\noindent 
 Then $M_{m,\Phi,\Psi}=M_{m,\Psi,\Phi}=I$. 
The unconditional convergence follows as in Example \ref{nbnb21}(i).

\vspace{.1in}
\hspace{.56in} (ii) Let 
$m=(\ \ \ \ \ \, 2, \ \ \  \  \frac{1}{2},  \ \ \ \ \, 2^2, \ \ \ \ \ \frac{1}{3}, \ \ \ \ \ \ 2^3,  \ \ \ \   \frac{1}{4},\ \ldots)$.

\vspace{.05in}\noindent 
 Then $M_{m,\Phi,\Psi}=M_{m,\Psi,\Phi}=G_1$ - non-invertible on $\h$ (see Lemma \ref{lemg}).  The unconditional convergence follows as in Example \ref{nbnb21}(ii).

\end{ex}

\begin{ex}\label{fnb13}
Let 
$\Phi=(\frac{1}{2}\,e_1, \ \ e_2, \ \ \frac{1}{2}\,e_2, \ \ \ \ \, \frac{1}{2}\,e_2, \ \ \frac{1}{2^2}\,e_1, \ \ \, e_3, \ \ \frac{1}{3}\,e_3, \ \ \  \ \  \frac{1}{3}\,e_3,  \ \ \frac{1}{2^3}\,e_1, \ \ \, e_4, \ \ \frac{1}{4}\,e_4, \ \ \  \ \,  \frac{1}{4}\,e_4,  \ \ldots)$, 

\vspace{.05in} \hspace{.99in}
$\Psi=(\ \ \, e_1, \ \ e_2, \ \ \frac{1}{2}\,e_2, \ \ - \frac{1}{2}\,e_2, \ \ \ \ \ \ \, e_1, \ \ \, e_3, \ \ \frac{1}{3}\,e_3, \ \ - \frac{1}{3}\,e_3,
\ \ \ \ \ \, \ e_1, \ \ \, e_4, \ \ \frac{1}{4}\,e_4, \ \ - \frac{1}{4}\,e_4, \ \ldots)$.

\vspace{.05in} 
\noindent Then $M_{(1),\Phi,\Psi}=M_{(1),\Psi,\Phi}=I$. The convergence is unconditional on $\h$, because 
 $M_{(1),\Phi,\Psi}\neweq M_{(1),\Psi,\Phi} \neweq M_{\nu,\Theta,\Theta}$, where 

\vspace{.05in} 
 $\Theta=(\frac{1}{\sqrt{2}}\,e_1, \ e_2, \ \frac{1}{2}\,e_2, \ \frac{1}{2}\,e_2, \ \, \frac{1}{\sqrt{2^2}}\,e_1, \ e_3, \ \frac{1}{3}\,e_3, \ \frac{1}{3}\,e_3, 
 \ \, \frac{1}{\sqrt{2^3}}\,e_1, \ e_4, \ \frac{1}{4}\,e_4, \ \frac{1}{4}\,e_4,\ \ldots)$ is Bessel for $\h$ and 
 
 \vspace{.05in} 
 $\nu\, =(\ \ \ \ \ \ \ \  1, \ \ \, 1, \ \ \ \ \ \, 1, \ \ \, -1, \ \ \ \ \ \  \ \ \ \ \ 1, \ \ \  1, \ \  \ \ \ \  1,\ \ \, -1, \ \ \ \ \ \ \ \ \ \  1, \ \ \ 1,\ \ \ \ \ \, 1, \ \ \, -1, \ \ldots)$ is $SN$ (apply Prop. \ref{bmh}).
\end{ex} 

\begin{ex}\label{fnb14}
Let 
$\Phi=(\frac{1}{2}e_1, \ \ \ \ \  e_2, \ \   \frac{1}{2^2}e_1, \ \ \ \ \ e_3, \ \ \ \frac{1}{2^3}e_1, \ \ \ \ \ \, e_4,\ \ldots)$, 

\vspace{.05in} \hspace{.99in}
$\Psi=(\ \ e_1,\ \ \, \frac{1}{2} e_2, \ \ \ \ \ \  e_1, \ \ \ \frac{1}{3} e_3, \ \ \ \ \ \ \, e_1,\ \ \  \frac{1}{4} e_4,\ \ldots)$.

\vspace{.1in}
\hspace{.58in} 
(i) Then $M_{(1),\Phi,\Psi}=M_{(1),\Psi,\Phi}=G_1$ - non-invertible on $\h$ (see Lemma \ref{lemg}).   The convergence is unconditional on $\h$, because $M_{(1),\Phi,\Psi} \neweq M_{(1),\Psi,\Phi} \neweq M_{(1),\Theta,\Theta}$, where 
$\Theta$ is the same as in Example \ref{nbnb21}(ii).

\vspace{.1in}
\hspace{.555in}  
(ii) Let $m=(\ \ \,  \ 1, \ \ \ \ \ \ \ 2, \ \ \ \ \ \ \ \, 1, \ \ \ \ \ \ \ 3, \ \ \ \ \ \ \ \, \ 1,  \ \ \ \ \ \ \ 4,\ \ldots)$.

\vspace{.05in} \noindent
Then $M_{m,\Phi,\Psi}=M_{m,\Psi,\Phi}=I$. The unconditional convergence follows as in Example \ref{nbnb21}(i).

\end{ex}

\begin{ex}\label{fnb33}
Let 
$\Phi=(\frac{1}{2}e_1, \ \ \ \ \,  e_2, \ \ \ \ \ \, \frac{1}{2}e_1, \ \ \ \ \ e_3, \ \ \ \, \frac{1}{2^2}e_1, \ \ \ \ \ \, e_4, \ \ \ \ \ \frac{1}{2^2}e_1, \ \ \ \ \  e_5, \ \ \ \frac{1}{2^3}e_1, \ \ \ \ \,  e_6, \ \ \ \ \ \, \frac{1}{2^3}e_1, \ \ \ \ \, e_7,\ \ldots)$, 

\vspace{.05in} \hspace{.99in}
$\Psi=(\ \, e_1,\ \ \ \ \, e_2, \ \ \ \  \frac{1}{\sqrt{2}} e_1, \ \ \ \ \ e_3, \ \ \ \ \ \ \, e_1,\ \ \ \ \ \,  e_4, \ \ \  \frac{1}{\sqrt{2^2}} e_1, \ \ \ \ \ e_5, \ \ \ \ \ \ \ e_1,\ \ \ \ \,  e_6, \ \ \  \frac{1}{\sqrt{2^3}} e_1, \ \ \ \ \ e_7,\ \ldots)$, 

\vspace{.05in} \hspace{.98in}
$m=(\ \ \, \ 1, \ \ \ \ \ \, 1,  \ \ \ \ \ \ \ \, \frac{1}{\sqrt{2}}, \ \ \ \ \ \ \, 1, \ \ \ \ \ \ \ \ \frac{1}{2},  \ \ \ \ \ \ \  1, \, \ \ \ \ \ \ \frac{1}{\sqrt{2^2}}, \ \ \ \ \ \ \ 1, \ \ \ \ \ \ \, \frac{1}{2^2}, \ \ \ \ \ \   1,  \ \ \ \ \ \ \ \frac{1}{\sqrt{2^3}}, \ \ \ \ \ \ \, 1, \ \ldots)$.

\vspace{.051in} \noindent
 Then $M_{m,\Phi,\Psi}=M_{m,\Psi,\Phi}=I$. 
The convergence is unconditional on $\h$, because $\Phi$ and $m\Psi$ are Bessel for $\h$ (apply Prop. \ref{bmh}).

\end{ex}

\begin{ex}\label{fnb16}
Let 
$\Phi=(\frac{1}{2}e_1, \ \ \ \ \ \, \ e_2, \ \  \, \frac{1}{2^2}e_1, \ \ \ \ \ \ \ \ \, e_3, \ \ \ \frac{1}{2^3}e_1, \ \ \ \ \ \ \ \, e_4,\ \ldots)$, 

\vspace{.05in} \hspace{.99in}
$\Psi=(\ \ e_1,\ \ \frac{1}{\sqrt{2}} e_2, \ \ \ \ \ \, \ e_1, \ \ \ \frac{1}{\sqrt{3}} e_3, \ \ \ \ \ \ \ e_1,\ \ \, \frac{1}{\sqrt{4}} e_4,\ \ldots)$, 

\vspace{.05in} \hspace{.98in}
$m=(\ \ \, \ 1, \ \ \ \ \, \ \frac{1}{\sqrt{2}},  \ \ \ \ \ \ \ \ 1, \ \ \ \ \ \ \ \frac{1}{\sqrt{3}}, \ \ \ \ \ \ \ \, \ 1, \ \ \ \ \ \ \frac{1}{\sqrt{4}},\  \ldots)$.

\vspace{.05in}
\noindent Then $M_{m,\Phi,\Psi}=M_{m,\Psi,\Phi}=G_1$ - non-invertible on $\h$ (see Lemma \ref{lemg}).  The unconditional convergence follows as in Example \ref{nbnb21}(ii).

\end{ex}

\begin{ex}\label{fnb18}
Let 
$\Phi=(\frac{1}{2}e_1, \ \ \ \ \ \ \ e_2, \ \  \, \frac{1}{2^2}e_1, \ \ \ \ \ \ \ e_3, \ \ \ \frac{1}{2^3}e_1, \ \ \ \ \ \ \, \, e_4,\ \ldots)$, 

\vspace{.05in} \hspace{.99in}
$\Psi=(\ \ e_1,\ \ \ \frac{1}{2^2} e_2, \ \ \ \ \ \, \ e_1, \ \ \ \frac{1}{3^2} e_3, \ \ \ \ \ \ \ e_1,\ \ \ \, \frac{1}{4^2} e_4,\ \ldots)$, 

\vspace{.05in} \hspace{.98in}
$m=(\ \ \ \, 1, \ \ \ \ \ \ \ \ \, 2,  \ \ \ \ \ \ \ \ 1, \ \ \ \ \ \ \ \  \ 3, \ \ \  \ \ \ \ \, \ 1,  \ \ \ \ \ \ \ \ \ 4,\ \ldots)$.

\vspace{.05in} 
\noindent Then $M_{m,\Phi,\Psi}=M_{m,\Psi,\Phi}=G_1$ - non-invertible on $\h$ (see Lemma \ref{lemg}).  The unconditional convergence follows as in Example \ref{nbnb21}(ii).

\end{ex}

\begin{ex}\label{fnb34}
Let 
$\Phi=(\frac{1}{2}\,e_1, \ \ e_2, \ \ \frac{1}{2}\,e_2, \ \ \ \ \, \frac{1}{2}\,e_2, \ \ \frac{1}{2^2}\,e_1, \ \  e_3, \ \ \frac{1}{3}\,e_3, \ \ \  \ \,  \frac{1}{3}\,e_3,  \ \ \frac{1}{2^3}\,e_1, \ \ \, e_4, \ \ \frac{1}{4}\,e_4, \ \ \  \ \,  \frac{1}{4}\,e_4,  \ \ldots)$, 

\vspace{.05in} \hspace{.99in}
$\Psi=(\ \ \ e_1, \ \ e_2, \ \  2\,e_2, \ \, - \ 2\,e_2, \ \ \ \ \ \ \, e_1, \ \  e_3, \ \ \, 3\,e_3, \ \,  -\ 3  \,e_3,
\ \ \ \ \ \ e_1, \ \ \, e_4, \ \ \, 4\,e_4, \ \, - \ 4\,e_4, \ \ldots)$.

\vspace{.05in} 
\noindent Then $M_{(1),\Phi,\Psi}=M_{(1),\Psi,\Phi}=I$. The convergence is unconditional on $\h$, because 
 $M_{(1),\Phi,\Psi}\neweq M_{(1),\Psi,\Phi} \neweq M_{\nu,\Theta,\Theta}$, where 

\vspace{.05in} 
 $\Theta=(\frac{1}{\sqrt{2}}\,e_1, \ \, e_2, \ \ \ e_2, \ \ \ e_2, \ \, \frac{1}{\sqrt{2^2}}\,e_1, \ \ e_3, \ \ \ \,e_3, \ \ \ \,e_3, 
 \ \, \frac{1}{\sqrt{2^3}}\,e_1,\ \ e_4, \ \ \ \,e_4, \ \ \ \,e_4,\ \ldots)$ is Bessel for $\h$ and 
 
 \vspace{.05in} 
 $\nu\, =(\ \ \ \ \ \  \ \, 1, \ \ \, \ 1, \ \ \ \ \, 1, \ \,  -1, \ \ \ \ \  \ \ \ \ \ \, 1, \ \ \ \, 1, \ \  \ \ \  1,\ \ \, -1,\ \ \ \ \ \ \ \ \ \ \, 1, \ \ \ \, 1,\ \ \ \ \ 1, \ \  -1, \ \ldots)$ is $SN$ (apply Prop. \ref{bmh}).

\end{ex}

\begin{ex}\label{fnb35}
Let 
$\Phi=(\frac{1}{2}\,e_1, \ \ e_2, \ \ \frac{1}{2}\,e_2, \ \ \ \ \, \, \frac{1}{2}\,e_2, \ \ \frac{1}{2^2}\,e_1, \ \  e_3, \ \ \frac{1}{3}\,e_3, \ \ \  \ \,  \frac{1}{3}\,e_3,  \ \ \frac{1}{2^3}\,e_1, \ \ \, e_4, \ \ \frac{1}{4}\,e_4, \ \ \  \ \,  \frac{1}{4}\,e_4,  \ \ldots)$, 

\vspace{.05in} \hspace{.99in}
$\Psi=(\ \ \ e_2, \ \ e_2, \ \  \, 2\,e_2, \ \, - \ 2\,e_2, \ \ \ \ \ \ \, e_2, \ \  e_3, \ \ \, 3\,e_3, \ \,  -\ 3  \,e_3,
\ \ \ \ \ \, \, e_2, \ \ \, e_4, \ \ \, 4\,e_4, \ \, - \ 4\,e_4, \ \ldots)$.

\vspace{.05in} 
\noindent Then $M_{(1),\Phi,\Psi}$ is not injective and $M_{(1),\Psi,\Phi}$ is not surjective. They are unconditionally convergent on $\h$, because 
 $M_{(1),\Phi,\Psi}\neweq  M_{(1),\Theta,\Xi}$ and $M_{(1),\Psi,\Phi} \neweq  M_{(1),\Xi,\Theta}$, where 

\vspace{.05in} 
 $\Theta=(\frac{1}{\sqrt{2}}\,e_1, \ \, e_2, \ \ \ e_2, \ \ \ \ \, e_2, \ \, \frac{1}{\sqrt{2^2}}\,e_1, \ \ e_3, \ \ \ \,e_3, \ \ \ \ \, e_3, 
 \ \, \frac{1}{\sqrt{2^3}}\,e_1,\ \ e_4, \ \ \ \,e_4, \ \ \ \ \, e_4,\ \ldots)$ and
 
 \vspace{.05in} 
 $\Xi\,=(\frac{1}{\sqrt{2}}\,e_2, \ \, e_2, \ \ \ e_2, \  - \,e_2, \ \, \frac{1}{\sqrt{2^2}}\,e_2, \ \ e_3, \ \ \ \,e_3, \  -\,e_3, 
 \ \, \frac{1}{\sqrt{2^3}}\,e_2,\ \ e_4, \ \ \ \,e_4, \ - \,e_4,\ \ldots)$ are Bessel for $\h$ (apply Prop. \ref{bmh}).
 \end{ex}

\begin{ex}\label{fnb15}
Let 
$\Phi=(\frac{1}{2}e_1, \ \ \ \ \, e_2, \ \  \frac{1}{2^2}e_1, \ \ \ \, e_3, \ \  \frac{1}{2^3}e_1, \ \ \ \ \, e_4,\ \ldots)$, 

\vspace{.05in} \hspace{.99in}
$\Psi=(\ \ \, e_1,\ \  2 e_2, \ \ \ \ \ \ e_1, \ \ 3 e_3, \ \  \ \ \ \ e_1,\ \ \, 4 e_4,\ \ldots)$.

\vspace{.1in}
\hspace{.58in} (i) Let 
$m=(\ \ \ \ \ 1, \ \ \ \ \ \frac{1}{2}, \ \ \ \ \ \ \ \, 1, \ \ \ \ \ \frac{1}{3},\ \ \ \ \, \ \ \  1,  \ \ \ \ \,\, \frac{1}{4},\ \ldots)$.

\vspace{.05in}
\noindent Then $M_{m,\Phi,\Psi}=M_{m,\Psi,\Phi}=I$. The unconditional convergence follows as in Example \ref{nbnb21}(i).

\vspace{.1in}
\hspace{.54in} (ii) Let 
$m=(\ \ \ \ \ 1, \ \ \ \  \frac{1}{2^2}, \ \ \ \ \ \ \ \, 1, \ \ \  \, \frac{1}{3^2}, \ \ \ \ \ \ \  1,  \ \ \  \,\, \frac{1}{4^2},\ \ldots)$.

\vspace{.05in}
\noindent Then $M_{m,\Phi,\Psi}=M_{m,\Psi,\Phi}=G_1$ - non-invertible on $\h$ (see Lemma \ref{lemg}). The unconditional convergence follows as in Example \ref{nbnb21}(ii).

\end{ex}

\begin{ex}\label{fnb36}
Let 
$\Phi=(\frac{1}{2^2}e_1, \ \ \ \, e_2, \ \  \frac{1}{2^4}e_1, \ \ \ \,e_3, \ \  \frac{1}{2^6}e_1,\, \ \ \ \ e_4,\ \ldots)$, 

\vspace{.05in} \hspace{.99in}
$\Psi=(\ \ \ \, e_1,\ \  2 e_2, \ \ \ \ \ \, e_1, \ \ 3 e_3, \ \ \ \ \ \, e_1,\ \ \, 4 e_4,\ \ldots)$.

\vspace{.1in}
\hspace{.58in} (i) Let 
$m=(\ \ \ \ \ \, 2, \ \ \ \ \, \frac{1}{2}, \ \ \ \ \ \, 2^2, \ \ \ \ \, \frac{1}{3}, \ \ \ \ \ \ 2^3,  \ \ \ \ \ \ \frac{1}{4},\ \ldots)$.

\vspace{.05in}
\noindent Then $M_{m,\Phi,\Psi}=M_{m,\Psi,\Phi}=I$. The unconditional convergence follows as in Example \ref{nbnb21}(i).

\vspace{.1in}
\hspace{.54in} (ii) Let 
$m=(\ \ \ \ \ \, 2, \ \ \  \, \frac{1}{2^2}, \ \ \ \ \ 2^2, \ \ \ \,  \frac{1}{3^2}, \ \ \ \ \ \, 2^3,  \ \ \  \,\, \frac{1}{4^2},\ \ldots)$.

\vspace{.05in}
\noindent Then $M_{m,\Phi,\Psi}=M_{m,\Psi,\Phi}=G_1$ - non-invertible on $\h$ (see Lemma \ref{lemg}). The unconditional convergence follows as in Example \ref{nbnb21}(ii).

\end{ex}

\begin{ex}\label{fnb37}
Let 
$\Phi=(\frac{1}{2}\,e_1, \ \ e_2, \ \ \frac{1}{2}\,e_2, \ \ \ \ \, \frac{1}{2}\,e_2, \ \ \frac{1}{2^2}\,e_1, \ \  e_3, \ \ \frac{1}{3}\,e_3, \ \ \  \ \,  \frac{1}{3}\,e_3,  \ \ \frac{1}{2^3}\,e_1, \ \ \, e_4, \ \ \frac{1}{4}\,e_4, \ \ \  \ \,  \frac{1}{4}\,e_4,  \ \ \frac{1}{2^4}\,e_1, \ \  e_5, \ \ \frac{1}{5}\,e_5, \ \ \  \ \,  \frac{1}{5}\,e_5,  \ \ldots)$, 

\vspace{.05in} \hspace{.99in}
$\Psi=(\ \ \ e_1, \ \ e_2, \ \  2\,e_2, \ \, - \ 2\,e_2, \ \ \ \ \ \ \, e_1, \ \  e_3, \ \ \frac{1}{3}\,e_3, \   -\ \frac{1}{3} \,e_3,
\ \ \ \ \ \ e_1, \ \ \, e_4, \ \ \, 4\,e_4, \ \, - \ 4\,e_4, \ \ \ \ \ \ \, e_1, \ \  e_5, \ \ \frac{1}{5}\,e_5, \   -\ \frac{1}{5} \,e_5, \ \ldots)$.

\vspace{.05in} 
\noindent Then $M_{(1),\Phi,\Psi}=M_{(1),\Psi,\Phi}=I$. The convergence is unconditional on $\h$, because 
 $M_{(1),\Phi,\Psi}\neweq M_{(1),\Psi,\Phi} \neweq M_{\nu,\Theta,\Theta}$, where 

\vspace{.05in} 
 $\Theta=(\frac{1}{\sqrt{2}}\,e_1, \ \, e_2, \ \ \ e_2, \ \ \ e_2, \ \, \frac{1}{\sqrt{2^2}}\,e_1, \ \ e_3, \ \ \frac{1}{3}\,e_3, \ \ \frac{1}{3}\,e_3, 
 \ \, \frac{1}{\sqrt{2^3}}\,e_1,\ \ e_4, \ \ \ \,e_4, \ \ \ \,e_4, \ \, \frac{1}{\sqrt{2^4}}\,e_1, \ \ e_5, \ \ \frac{1}{5}\,e_5, \ \ \frac{1}{5}\,e_5,  \ \ldots)$ is Bessel for $\h$ and 
 
 \vspace{.05in} 
 $\nu\, =(\ \ \ \ \ \, \ \ 1, \ \ \ \, 1, \ \ \ \ \, 1, \  \, -1, \ \ \ \ \  \ \ \ \ \ 1, \ \ \ \  1, \ \ \ \ \ \ \ 1,\ \ \ \, -1, \ \ \, \ \ \ \ \ \ \ \ \, 1, \ \ \ 1,\ \ \ \ \, \ 1, \ \  -1, \ \ \ \ \  \ \ \ \ \ \, 1, \ \ \ \, 1, \ \  \ \ \ \ \ 1,\ \ \ \, -1, \ \ldots)$ is $SN$ (apply Prop. \ref{bmh}).
 \end{ex}

\begin{ex}\label{fnb38}
Let 
$\Phi=(\frac{1}{2}\,e_1, \ \ e_2, \ \ \frac{1}{2}\,e_2, \ \ \ \ \, \frac{1}{2}\,e_2, \ \ \frac{1}{2^2}\,e_1, \ \  e_3, \ \ \frac{1}{3}\,e_3, \ \ \  \ \,  \frac{1}{3}\,e_3,  \ \ \frac{1}{2^3}\,e_1, \ \ \, e_4, \ \ \frac{1}{4}\,e_4, \ \ \  \ \,  \frac{1}{4}\,e_4,  \ \ \frac{1}{2^4}\,e_1, \ \  e_5, \ \ \frac{1}{5}\,e_5, \ \ \  \ \,  \frac{1}{5}\,e_5,  \ \ldots)$, 

\vspace{.05in} \hspace{.99in}
$\Psi=(\ \ \ e_2, \ \ e_2, \ \  2\,e_2, \ \, - \ 2\,e_2, \ \ \ \ \ \ \, e_2, \ \  e_3, \ \ \frac{1}{3}\,e_3, \   -\ \frac{1}{3} \,e_3,
\ \ \ \ \ \, \, e_2, \ \ \, e_4, \ \ \, 4\,e_4, \ \, - \ 4\,e_4, \ \ \ \ \ \ \, e_2, \ \  e_5, \ \ \frac{1}{5}\,e_5, \   -\ \frac{1}{5} \,e_5, \ \ldots)$.

\vspace{.05in} 
\noindent Then $M_{(1),\Phi,\Psi}$ is not injective, because $M_{(1),\Phi,\Psi}e_1=0$. The multiplier $M_{(1),\Psi,\Phi}$ is not surjective, because $\Psi$ is not complete in $\h$. Both $M_{(1),\Phi,\Psi}$ and $M_{(1),\Psi,\Phi}$
 are unconditionally convergent on $\h$, because 
 $M_{(1),\Phi,\Psi}\neweq  M_{(1),\Theta,\Xi}$ and $M_{(1),\Psi,\Phi} \neweq  M_{(1),\Xi,\Theta}$, where 

\vspace{.05in} 
 $\Theta=(\frac{1}{\sqrt{2}}\,e_1, \ \, e_2, \ \ \ e_2, \ \ \ \ \, e_2, \ \, \frac{1}{\sqrt{2^2}}\,e_1, \ \ e_3, \ \ \ \frac{1}{3}\,e_3, \ \ \ \ \, \frac{1}{3}\,e_3, 
 \ \ \frac{1}{\sqrt{2^3}}\,e_1,\ \ \, e_4, \ \ \ \, e_4, \ \ \ \ e_4, \ \ \frac{1}{\sqrt{2^4}}\,e_1, \ \ \, e_5, \ \ \frac{1}{5}\,e_5, \ \ \ \ \, \frac{1}{5}\,e_5,  \ \ldots)$ and 

 \vspace{.05in} 
 $\Xi\, =(\frac{1}{\sqrt{2}}\,e_2, \ \, e_2, \ \ \ e_2, \  - \,e_2, \ \, \frac{1}{\sqrt{2^2}}\,e_2, \ \ e_3, \ \ \ \frac{1}{3}\,e_3, \  -\,\frac{1}{3}\,e_3, 
 \ \, \frac{1}{\sqrt{2^3}}\,e_2,\ \ \  e_4, \  \ \ \,e_4, \ - \,e_4,\ \, \frac{1}{\sqrt{2^4}}\,e_2, \ \ e_5, \ \ \frac{1}{5}\,e_5, \ -\, \frac{1}{5}\,e_5, \ \ldots)$ are Bessel for $\h$ (apply Prop. \ref{bmh}).

\end{ex}

\begin{ex}\label{fnb39}
Let 
$\Phi=(\frac{1}{2}\,e_1, \ \ \  e_2, \  \ \ \ \ e_2,\ \ \ \ \ \ \ \, e_2, \ \ \, \frac{1}{2^2}\,e_1, \ \ \,  e_3, \ \ \ \ \ e_3, \ \ \ \ \ \ \ \ e_3,  \ \ \frac{1}{2^3}\,e_1, \ \ \ e_4, \ \ \ \ \, e_4, \ \ \ \ \ \ \ \ \ e_4,\ \ \, \frac{1}{2^4}\,e_1, \ \ \  e_5, \ \ \ \ \ e_5, \ \ \ \ \ \ \ \, e_5,\ \ldots)$, 

\vspace{.05in} \hspace{.99in}
$\Psi=(\ \ \, e_1, \ \ \ e_2, \ \,  2\,e_2, \ \, - \ 2\,e_2, \ \ \ \ \ \ \ e_1, \ \ \ e_3, \ \ \frac{1}{3}\,e_3, \  \, - \frac{1}{3} \,e_3,
\ \ \ \ \ \, \ e_1, \ \ \ e_4, \ \ \, 4\,e_4, \ \, - \ 4\,e_4, \ \, \ \ \ \ \ \ e_1, \ \ \ e_5, \ \ \frac{1}{5}\,e_5, \  \, - \frac{1}{5} \,e_5,\ \ldots)$. 

\vspace{.1in}
\hspace{.58in} (i) Let 
$m=(\ \ \ \ \ 1, \ \ \ \ \, 1,\ \ \ \ \, \ \frac{1}{2}, \ \ \ \ \ \ \ \ \ \frac{1}{2}, \ \ \ \ \ \ \ \ \, 1, \ \ \ \ \, 1, \ \ \ \ \ \ \frac{1}{3}, \ \ \ \ \ \ \ \ \ \, \frac{1}{3},\ \ \ \ \ \ \  \, 1,  \ \ \ \ \, 1,\ \ \ \ \,\ \ \frac{1}{4}, \ \ \ \ \ \ \ \ \ \, \frac{1}{4}, \ \ \ \ \ \ \ \ \ 1, \ \ \ \ \  1, \ \ \ \ \ \ \frac{1}{5},\ \ \ \ \ \ \ \ \ \,\frac{1}{5},\ \ldots)$.

\vspace{.05in} 
\noindent Then $M_{m,\Phi,\Psi}=M_{m,\Psi,\Phi}=I$. The convergence is unconditional on $\h$, because 
 $M_{m,\Phi,\Psi}\neweq M_{m,\Psi,\Phi} \neweq M_{\nu,\Theta,\Theta}$, where $\Theta$ and $\nu$ are the same as in Example \ref{fnb37}.

 \vspace{.1in}
\hspace{.54in} (ii) Let 
$m=(\ \ \ \ \ 1, \ \ \ \ \frac{1}{2},\ \ \ \ \, \ \frac{1}{2}, \ \ \ \ \ \ \ \ \ \frac{1}{2}, \ \ \ \ \ \ \ \ \ 1, \ \ \ \ \frac{1}{3}, \ \ \ \ \ \  \frac{1}{3}, \ \ \ \ \ \ \ \ \ \,  \frac{1}{3},\ \ \ \ \ \ \  \, 1,  \ \ \ \ \, \frac{1}{4},\ \ \ \ \,\, \ \frac{1}{4}, \ \ \ \ \ \ \ \ \ \, \frac{1}{4}, \ \ \ \ \ \ \  \ \ 1, \ \ \ \  \frac{1}{5}, \ \ \ \ \ \ \frac{1}{5},\ \ \ \ \ \ \ \ \ \, \frac{1}{5},\ \ldots)$.
 
 \vspace{.05in} 
\noindent Then $M_{m,\Phi,\Psi}=M_{m,\Psi,\Phi}=G_1$ - non-invertible on $\h$ (see Lemma \ref{lemg}). The convergence is unconditional on $\h$, because 
 $M_{m,\Phi,\Psi}\neweq M_{m,\Psi,\Phi} \neweq M_{\nu,\Theta,\Theta}$, where 
 
\vspace{.05in} 
 $\Theta=(\frac{1}{\sqrt{2}}\,e_1, \ \, \frac{1}{\sqrt{2}} \,e_2, \ \ \ \, e_2, \ \ \ e_2, \ \, \frac{1}{\sqrt{2^2}}\,e_1, \ \ \frac{1}{\sqrt{3}} \, e_3, \ \, \frac{1}{3} e_3, \ \ \frac{1}{3} e_3, 
 \ \, \frac{1}{\sqrt{2^3}}\,e_1,\ \ \frac{1}{\sqrt{4}} \, e_4, \ \ \ \,e_4, \ \ \ \,e_4,
 \ \, \frac{1}{\sqrt{2^4}}\,e_1,\ \ \frac{1}{\sqrt{5}} \, e_5, \  \frac{1}{5} e_5, \ \ \frac{1}{5} e_5,\ \ldots)$ is Bessel for $\h$ and 
 
 \vspace{.05in} 
 $\nu\, =(\ \ \ \ \ \, \ \ 1, \ \ \ \, \ \ \ \ \ \, 1, \ \ \ \ \ 1, \ \  -1, \ \ \ \ \ \, \ \ \ \ \ 1, \ \ \ \ \ \ \ \  \ 1, \ \ \ \, \ \ \  1,\ \ \, -1, \ \ \ \ \ \ \ \ \ \  \ 1,\ \ \ \ \ \ \ \ \ \, 1,\ \ \ \ \  1, \ \  -1, 
 \ \ \ \ \ \ \ \ \ \ \  1, \ \ \ \ \ \ \ \ \ \, 1,\ \ \ \ \ 1, \ \ \ -1,\ \ldots)$ is $SN$ (apply Prop. \ref{bmh}).
 
 \end{ex}

\begin{ex}\label{fnb40}
Let 
$\Phi=(\frac{1}{2}e_1, \ \ \ \ \  e_2, \ \  \, \frac{1}{2^2}e_1, \ \ \ \ \ e_3, \ \ \ \frac{1}{2^3}e_1, \ \ \ \ \ \, e_4, \ \ \ \frac{1}{2^4}e_1, \ \ \ \ \  e_5,\ \ldots)$, 

\vspace{.05in} \hspace{.99in}
$\Psi=(\ \ e_1,\ \ \ \frac{1}{2} e_2, \ \ \ \ \ \ e_1, \ \ \ \, 3e_3, \ \ \ \ \ \ e_1,\ \ \ \, \frac{1}{4} e_4, \ \ \ \ \ \ \ e_1,\ \ \ \, 5 e_5,\ \ldots)$, 

\vspace{.05in} \hspace{.98in}
$m=(\ \ \, \ 1, \ \ \ \ \ \ \, 2,  \ \ \ \ \ \ \ \  1, \ \ \ \ \ \ \, \frac{1}{3}, \ \ \ \ \ \ \ \, 1,  \ \ \ \ \ \ \ \ \,  4, \ \ \ \ \ \ \ \  1,  \ \ \ \ \ \ \  \frac{1}{5},\ \ldots)$.

\vspace{.05in} \noindent
Then $M_{m,\Phi,\Psi}=M_{m,\Psi,\Phi}=I$. The unconditional convergence follows as in Example \ref{nbnb21}(i).

\end{ex}

\begin{ex}\label{fnb41}
Let 
$\Phi=(\frac{1}{2}e_1, \ \ \ \ \ \ e_2, \ \  \, \frac{1}{2^2}e_1, \ \ \ \ \ \, e_3, \ \ \ \frac{1}{2^3}e_1, \ \ \ \ \ \ e_4, \ \ \ \frac{1}{2^4}e_1, \ \ \ \ \ \, e_5,\ \ldots)$, 

\vspace{.05in} \hspace{.99in}
$\Psi=(\ \ e_1,\ \  \frac{1}{2^2} e_2, \ \ \ \ \ \, \ e_1, \ \ \ \, 3e_3, \ \ \ \ \ \ \ e_1,\ \ \, \frac{1}{4^2} e_4, \ \ \ \ \ \ \, e_1,\ \ \ \, 5 e_5,\ \ldots)$, 

\vspace{.05in} \hspace{.98in}
$m=(\ \ \ \,  1, \ \ \ \ \ \ \ \, 2,  \ \ \ \ \ \ \ \ 1, \ \ \ \ \ \,  \frac{1}{3^2}, \ \ \ \ \ \ \ \  1,  \ \ \ \ \ \ \ \ \, 4, \ \ \ \ \ \ \ \ \, 1,  \ \ \ \ \    \frac{1}{5^2},\ \ldots)$.

\vspace{.05in} \noindent
Then $M_{(1),\Phi,\Psi}=M_{(1),\Psi,\Phi}=G_1$ - non-invertible on $\h$ (see Lemma \ref{lemg}).   The convergence is unconditional on $\h$, because $M_{(1),\Phi,\Psi} \neweq M_{(1),\Psi,\Phi} \neweq M_{(1),\Theta,\Theta}$, where 
$\Theta$ is the same as in Example \ref{nbnb21}(ii).

\end{ex}

\subsection{Examples for one overcomplete frame and one Bessel non-frame sequence; TABLE 5 
on page \pageref{table5}}

\begin{ex}\label{fb1}
Let $\Phi=(e_1, e_1, e_2, e_3, e_4, e_5, \ldots)$ and $\Psi=(e_2,e_3,e_4,e_5,e_6,\ldots)$. 

\vspace{.05in}
(i) Then $M_{(1),\Phi,\Psi}$ and $M_{(1),\Psi,\Phi}$ are unconditionally convergent on $\h$ due to Proposition \ref{bmh} and non-invertible on $\h$ due to Proposition \ref{c1}.

\vspace{.05in}
(ii) If  $m=(\frac{1}{n})$, then $M_{m,\Phi,\Psi}$ and $M_{m,\Psi,\Phi}$ are unconditionally convergent on $\h$ due to Proposition \ref{bmh} and non-invertible on $\h$ due to Proposition \ref{c1}. \end{ex}

\begin{ex}\label{fb3}
Let $\Phi=(e_1, e_1, e_2, e_3, e_4, e_5,  \ldots)$ and $\Psi=(\frac{1}{2}\,e_1,\frac{1}{2}\,e_1,\frac{1}{2}\,e_2,\frac{1}{3}\,e_3,\frac{1}{4}\,e_4,\frac{1}{5}\,e_5,\ldots)$. 

\vspace{.05in}
(i) Then $M_{(1),\Phi,\Psi}=M_{(1),\Psi,\Phi}=G_1$ - non-invertible on $\h$ (see Lemma \ref{lemg}). The unconditional convergence on $\h$ follows from Proposition \ref{bmh}.

\vspace{.05in}
(ii) If $m=(1,1,\frac{1}{2},\frac{1}{3},\frac{1}{4},\frac{1}{5},\ldots)$, then $M_{m,\Phi,\Psi}=M_{m,\Psi,\Phi}=G_2$ - non-invertible on $\h$ (see Lemma \ref{lemg}). The unconditional convergence follows from Proposition \ref{bmh}. 

\vspace{.05in}
(iii) If $m=(1,1,2,3,4,5,\ldots)$, then $M_{m,\Phi,\Psi}=M_{m,\Psi,\Phi}=I$. The convergence is unconditional on $\h$, because $\Phi$ and $m\Psi$ are Bessel for $\h$ (apply Prop. \ref{bmh}).
\end{ex}

\begin{ex}\label{fb6}
Let 
$\Phi=( \ \ \, \, e_1,\ \ \ \ \ e_1,\ \ \ \ \ \, \ e_2, \ \ \ \ \ \ \ e_3, \ \ \ \ \ \, \ e_4, \ \ \ \ \ \, \ e_5,\ \ldots)$,

\vspace{.05in} \hspace{.91in}
 $\Psi=(\frac{1}{2}\,e_1, \ \ \frac{1}{2}\,e_1, \ \ \frac{1}{2^2}\,e_2, \  \ \frac{1}{3^2}\,e_3, \ \ \frac{1}{4^2}\,e_4, \ \ \frac{1}{5^2}\,e_5, \ \ldots)$,
 
 \vspace{.05in} \hspace{.91in}
  $m=(\ \ \ \,\,\, 1, \ \ \ \ \ \ 1, \ \ \ \ \ \ \ \ \ 2, \ \ \ \ \ \ \ \ \ 3, \ \ \ \ \ \ \ \ 4, \ \ \ \ \ \ \ \ \, 5, \ \ldots)$. 
  
  \vspace{.05in}
  \noindent Then $M_{m,\Phi,\Psi}=M_{m,\Psi,\Phi}=G_1$ - non-invertible on $\h$  (see Lemma \ref{lemg}). The convergence is unconditional on $\h$, because $\Phi$ and $m\Psi$ are Bessel for $\h$ (apply Prop. \ref{bmh}).

\end{ex}

\begin{ex}\label{fb7}
Let $\Phi=(\frac{1}{2}\,e_1, e_2, \frac{1}{2^2}\,e_1, e_3, \frac{1}{2^3}\,e_1, e_4, \ldots)$ and $\Psi=(e_2,e_3,e_4,e_5,\ldots)$.

\vspace{.05in}
(i) Then $M_{(1),\Phi,\Psi}$ and $M_{(1),\Psi,\Phi}$ are unconditionally convergent on $\h$ due to Proposition \ref{bmh} and non-invertible on $\h$ due to Proposition \ref{c1}.

\vspace{.05in}
(ii) If  $m=(\frac{1}{n})$, then $M_{m,\Phi,\Psi}$ and $M_{m,\Psi,\Phi}$ are unconditionally convergent on $\h$ due to Proposition \ref{bmh} and non-invertible on $\h$ due to Proposition \ref{c1}.
\end{ex}

\begin{ex}\label{fb10}
Let 
$\Phi=(\frac{1}{2}\,e_1, \ \ e_2, \ \ \frac{1}{2^3}\,e_1, \ \  e_3, \ \  \frac{1}{2^5}\,e_1, \ \ \, e_4, \ \ \frac{1}{2^7}\,e_1, \ \ \, e_5,\ \ldots)$,

\vspace{.05in} \hspace{.91in}
 $\Psi=(\ \ \, e_2, \ \ \, e_3,\ \ \ \ \ \ e_4, \ \ \, e_5,\ \ \ \ \ \, \ e_6, \ \ \, e_7,\ \ \ \ \, \ \ e_8, \ \ \, e_9,\ \ldots)$, 
 
 \vspace{.05in} \hspace{.905in}
 $m=( \ \ \  \ 1,\ \ \ \,\, 1, \ \ \ \ \ \ \, \ 2,\ \ \ \,\, 1,\ \ \ \ \ \, \ 2^2,\ \ \ \ 1,\ \ \ \ \, \ \ 2^3,\ \ \ \, \ 1,\ \ldots)$.
  
  \vspace{.05in}
\noindent  
Since $m\Phi$ is a frame for $\h$ and $\Psi$ is Bessel for $\h$, which is not a frame for $\h$, it follows that $M_{(1),m\Phi,\Psi}$ and $M_{(1),\Psi,m\Phi}$ are unconditionally convergent 
(see Prop. \ref{bmh}), but not invertible on $\h$ (see Prop. \ref{c1}). 

\end{ex}

\begin{ex}\label{fb11}
Let $\Phi=(\frac{1}{2}\,e_1, e_2, \frac{1}{2^2}\,e_1, e_3, \frac{1}{2^3}\,e_1, e_4, \ldots)$, and $\Psi=(\frac{1}{n}e_n)$.

\vspace{.05in}
(i) Then $M_{(1),\Phi,\Psi}$ and $M_{(1),\Psi,\Phi}$ are unconditionally convergent on $\h$ due to Proposition \ref{bmh} and non-invertible on $\h$ due to Proposition \ref{c1}.

\vspace{.05in}
(ii) If  $m=(\frac{1}{n})$, then $M_{m,\Phi,\Psi}$ and $M_{m,\Psi,\Phi}$ are unconditionally convergent on $\h$ due to Proposition \ref{bmh} and non-invertible on $\h$ due to Proposition \ref{c1}.
\end{ex}

\begin{ex}\label{fb13}
Let
 $\Phi=(\frac{1}{2}\,e_1,\ \ \ \ \, \, e_2,\ \  \frac{1}{2^2}\,e_1, \ \ \ \ \,\, e_3, \ \  \frac{1}{2^3}\,e_1, \ \ \ \ \,  e_4, \ \ldots)$, 
 
 \vspace{.05in} \hspace{.91in}
 $\Psi=(\frac{1}{2}\,e_1,\ \ \frac{1}{2}\,e_2, \ \ \frac{1}{2^2}\,e_1, \ \ \frac{1}{3}\,e_3, \ \ \frac{1}{2^3}\,e_1, \ \ \frac{1}{4}\,e_4, \ \ldots)$,

 \vspace{.1in}
\hspace{.51in} (i) Let 
 $m=(\ \ \ \, \ 2,  \ \ \ \ \ \, \ \ 2, \ \ \ \ \ \, 2^2, \ \ \ \ \ \ \ 3, \ \ \ \ \ \ \   2^3, \ \ \ \ \ \ \ 4,\ \ldots)$.

  \vspace{.05in} \noindent
   Then  $M_{m,\Phi,\Psi}=M_{m,\Psi,\Phi}=I$. The unconditional convergence follows as in Example \ref{nbnb21}(i).

  \vspace{.1in}
\hspace{.48in} (ii) Let 
 $m=(\ \ \ \, \ 2,  \ \ \ \ \ \, \ \ 1, \ \ \ \ \ \, 2^2, \ \ \ \ \ \ \ 1, \ \ \ \ \ \ \  2^3, \ \ \ \ \ \ \ 1,\ \ldots)$.
   
  \vspace{.05in} \noindent
   Then  $M_{m,\Phi,\Psi}=M_{m,\Psi,\Phi}=G_1$. The unconditional convergence follows as in Example \ref{nbnb21}(ii).

\end{ex}

\subsection{Examples for two overcomplete frames; TABLE 6 
on page \pageref{table6}} 

\begin{ex}\label{ff1}
Let $\Phi=(e_1, e_1, e_2, e_2, e_3, e_3, \ldots)$.

\vspace{.05in}
(i) If $m=(\frac{1}{2})$, then $M_{m,\Phi,\Phi}=I$.

\vspace{.05in}
(ii) If $m=(\frac{1}{2},1- \frac{1}{2},\frac{1}{3},1- \frac{1}{3}, \frac{1}{4},1- \frac{1}{4},\ldots)$, then $M_{m,\Phi,\Phi}=I$.

 \vspace{.05in}
(iii) If  $m=(\frac{1}{2}, \frac{1}{2}, \frac{1}{4}, \frac{1}{4},\frac{1}{6}, \frac{1}{6},\ldots)$, then $M_{m,\Phi,\Phi}=G_1$ - non-invertible on $\h$ (see Lemma \ref{lemg}).

 \vspace{.05in} \noindent The convergence of $M_{m,\Phi,\Phi}$  in (i)-(iii) is unconditional on $\h$ due to Prop. \ref{bmh}.
\end{ex}

\begin{ex}\label{ff2}
Let $\Phi=(e_1, e_1, e_2, e_2, e_3, e_3, \ldots)$ and $\Psi=(e_1, e_1, e_2, e_3, e_4, e_5, \ldots)$. 
By Proposition \ref{bmh}, $M_{(1),\Phi,\Psi}$ and  $M_{(1),\Psi,\Phi}$ are unconditionally convergent on $\h$. However, 
$M_{(1),\Phi,\Psi}$ is  not injective, for example $M_{(1),\Phi,\Psi}e_2=e_2=M_{(1),\Phi,\Psi}e_3$. 
Furthermore, $M_{(1),\Psi,\Phi}$ is injective, but not surjective. Indeed, observe that $M_{(1),\Psi,\Phi}= M_{(1),(e_n),\Theta}$, where $\Theta=(2e_1,e_2,e_2,e_3,e_3,\ldots)$ is an overcomplete frame for $\h$, and apply Proposition \ref{rbis}(a).
\end{ex}

\begin{ex}\label{ff9}
Let 
$\Phi=(e_1, \ e_1,\, -\,e_1, \ \ \ e_2,\ \ \ \ \, e_2,\ \, - \,e_2, \ \ \ e_3, \ \ \ \,\, e_3, \ \ - \,e_3, \ \ldots)$,

\vspace{.05in} \hspace{.91in}
$\Psi=(e_1, \ e_1,\ \ \ \, \, e_1,\ \ \ \, e_2,\ \ \frac{1}{2}e_2, \ \ \, \frac{1}{2}e_2, \ \ \ e_3, \ \ \frac{1}{3}e_3, \ \ \ \frac{1}{3}e_3, \ \ldots)$.

\vspace{.05in} 
(i) If $m=(1)$, then $M_{m,\Phi,\Psi}=M_{m,\Psi,\Phi}=I$. The convergence is unconditional on $\h$ due to Proposition \ref{bmh}.

\vspace{.05in}
(ii) If $m=(1,1,1,1,\frac{1}{2},\frac{1}{2},1,\frac{1}{3},\frac{1}{3},\ldots)$, then $M_{m,\Phi,\Psi}=M_{m,\Psi,\Phi}=I$.  The convergence is unconditional on $\h$ due to Proposition \ref{bmh}.

\vspace{.05in}
(iii) If $m=(1,1,1,\frac{1}{2},\frac{1}{2},\frac{1}{2},\frac{1}{3},\frac{1}{3},\frac{1}{3},\ldots)$, then $M_{m,\Phi,\Psi}=M_{m,\Psi,\Phi}=G_1$ - non-invertible on $\h$ (see Lemma \ref{lemg}). The convergence is unconditional on $\h$ due to Proposition \ref{bmh}.

\vspace{.05in}
(iv) If $m=(1,1,1,1,2,2,1,3,3,\ldots)$, then $M_{m,\Phi,\Psi}=M_{m,\Psi,\Phi}=I$.  The convergence is unconditional on $\h$, because $\Phi$ and $m\Psi$ are Bessel for $\h$ (apply Prop. \ref{bmh}).

\vspace{.05in}
(v) If $m=(1,1,1,\frac{1}{2},2,2,\frac{1}{3},3,3,\ldots)$, then $M_{m,\Phi,\Psi}=M_{m,\Psi,\Phi}=G_1$ - non-invertible on $\h$ (see Lemma \ref{lemg}).  The convergence is unconditional on $\h$, because $\Phi$ and $m\Psi$ are Bessel for $\h$ (apply Prop. \ref{bmh}).
\end{ex}

\begin{ex}\label{ff6}
Let 
$\Phi=(e_1, \ \ e_2,  \ \ e_1, \ \  e_3, \ \  e_4, \ \ \ \ \,\,\, e_4, \ \ \  e_5, \ \ \ \ \ \, \ e_5, \ \  \ e_6, \, \,\ \ \ \ \ \, e_6, \ \  \ \ e_7, \ \ \ \ \ \,\, e_7, \ \ldots)$,

\vspace{.05in} \hspace{.91in}
$\Psi=(e_2, \ \ e_1,  \ \ e_3, \ \  e_1, \  \ e_4, \ \ \, \frac{1}{2}\,e_1, \ \ \, e_5, \ \ \, \frac{1}{2^2}\,e_1, \ \ \, e_6, \ \  \frac{1}{2^3}\,e_1, \ \ \ \ e_7,\ \  \frac{1}{2^4}\,e_1, \ \ldots)$.

\vspace{.05in}
\noindent 
Then $M_{(1),\Phi,\Psi}$ and $M_{(1),\Psi,\Phi}$ are unconditionally convergent on $\h$ due to Proposition \ref{bmh}. 
However,  $M_{(1),\Phi,\Psi}$ and $M_{(1),\Psi,\Phi}$  are not injective, for example $M_{(1),\Phi,\Psi}e_2=e_1=M_{(1),\Phi,\Psi}e_3$ and $M_{(1),\Psi,\Phi}e_2=e_1=M_{(1),\Psi,\Phi}e_3$.

\end{ex}

\begin{ex}\label{ff11}
Let 
$\Phi=(\frac{1}{\sqrt{2}}\,e_1, e_2, \frac{1}{\sqrt{2^2}}\,e_1, e_3, \frac{1}{\sqrt{2^3}}\,e_1, e_4, \ldots)$. Then $M_{(1),\Phi,\Phi}=I$.  The convergence is unconditional on $\h$ due to Proposition \ref{bmh}.
\end{ex}

\begin{ex}\label{ff12}
Let 
$\Phi=(e_1, \ \ e_2,  \ \ e_1, \ \  e_3, \ \  e_4, \ \ \  \frac{1}{2}\,e_1, \ \ \  e_5, \ \ \  \frac{1}{2^2}\,e_1, \ \  \ e_6,  \ \ \ \frac{1}{2^3}\,e_1, \ \  \, e_7,  \ \ldots)$,

\vspace{.05in} \hspace{.91in}
$\Psi=(e_2, \ \ e_1,  \ \ e_3, \ \  e_1, \  \ e_4, \ \  \ \frac{1}{2}\,e_1, \ \ \  e_5, \ \ \  \frac{1}{2^2}\,e_1, \ \  \ e_6,  \ \  \ \frac{1}{2^3}\,e_1, \ \  \, e_7,  \ \ldots)$.

\vspace{.05in}
\noindent Then $M_{(1),\Phi,\Psi}$ and $M_{(1),\Psi,\Phi}$ are unconditionally convergent on $\h$ due to Prop. \ref{bmh}. However, $M_{(1),\Phi,\Psi}$ and $M_{(1),\Psi,\Phi}$  are not injective, for example $M_{(1),\Phi,\Psi}e_2=e_1=M_{(1),\Phi,\Psi}e_3$ and $M_{(1),\Psi,\Phi}e_2=e_1=M_{(1),\Psi,\Phi}e_3$.
\end{ex}

\begin{ex}\label{ff13}
Let 
$\Phi=(\frac{1}{\sqrt{2}}\,e_1, \ \ e_2, \ \ \frac{1}{\sqrt{2^2}}\,e_1, \ \  \, e_3, \ \ \frac{1}{\sqrt{2^3}}\,e_1, \ \ e_4,\ \ldots)$,

\vspace{.05in} \hspace{.91in}
$\Psi=(\frac{1}{\sqrt[4]{2}}\,e_1, \ \ e_2,\ \  \frac{1}{\sqrt[4]{2^2}}\,e_1, \ \ e_3, \ \ \frac{1}{\sqrt[4]{2^3}}\,e_1, \ \ e_4, \ \ldots)$.

\vspace{.05in} 
 (i) If $m=(\frac{1}{\sqrt[4]{2}} \ , \ \, 1, \ \frac{1}{\sqrt[4]{2^2}} \ , \ \ \ 1, \ \frac{1}{\sqrt[4]{2^3}}\ , \ \ 1, \ldots)$, then $M_{m,\Phi,\Psi}=M_{(1),\Psi,\Phi}=I$.

\vspace{.05in}
 (ii) If $m=(\frac{1}{\sqrt[4]{2}} \ , \  \frac{1}{2}, \ \frac{1}{\sqrt[4]{2^2}} \ , \ \, \frac{1}{3}, \ \frac{1}{\sqrt[4]{2^3}}\ , \ \, \frac{1}{4}, \ldots)$, then $M_{m,\Phi,\Psi}=M_{m,\Psi,\Phi}=G_1$ - non-invertible on $\h$ (see Lemma \ref{lemg}). 

\vspace{.05in} \noindent
The convergence of $M_{m,\Phi,\Psi}$ and $M_{m,\Psi,\Phi}$ in (i)-(ii) is unconditional on $\h$ due to Proposition \ref{bmh}.
\end{ex}

\begin{ex}\label{ff15}
Let 
$\Phi=(\frac{1}{2}\,e_1, \ e_2, \  \frac{1}{2^2}\,e_1,  \ e_3, \ \frac{1}{2^3}\,e_1,  \ e_4,\ \ldots)$.

\vspace{.05in} 
 (i) \, If $m=(2,\ 1, \ 2^2, \ 1, \ 2^3, \ 1, \ \ldots)$, then $M_{m,\Phi,\Phi}=I$. 

\vspace{.05in} 
 (ii) If $m=(2,\ \frac{1}{2}, \ 2^2, \ \frac{1}{3}, \ 2^3, \ \frac{1}{4}, \ \ldots)$, then $M_{m,\Phi,\Phi}=G_1$ - non-invertible on $\h$ (see Lemma \ref{lemg}). 

\vspace{.05in} \noindent
The convergence of $M_{m,\Phi,\Phi}$ in (i)-(ii) is unconditional on $\h$, because $M_{m,\Phi,\Phi}\neweq M_{(1),(\sqrt{m_n}\phi_n),(\sqrt{m_n}\phi_n)}$ and $(\sqrt{m_n}\phi_n)$ is Bessel for $\h$ (apply Prop. \ref{bmh}).
\end{ex}

\subsection{Examples for one Riesz basis and one non-Bessel sequence; TABLE 7 
on page \pageref{table7}}

\begin{ex}\label{rnb1}
Let $\Phi=\seq[e]$, $\Psi=(e_1,e_2,e_1,e_3,e_1,e_4,\ldots)$, and $m=(\frac{1}{2},1,\frac{1}{2^2},1,\frac{1}{2^3},1,\ldots)$. 
The multipliers $M_{m,\Phi,\Psi}$ and $M_{m,\Psi,\Phi}$ are unconditionally convergent on $\h$, because $\Phi$ and $m\Psi$ are Bessel for $\h$ (apply Prop. \ref{bmh}).  
The non-invertibility is due to Proposition \ref{rc}.
\end{ex}

\begin{ex}\label{rnb31}
Let  $\Phi=\seq[e]$, $\Psi=(e_1,\frac{1}{2}e_2,e_1, \frac{1}{3}e_3,e_1, \frac{1}{4}e_4,\ldots)$, and $m=(\frac{1}{2},1,\frac{1}{2^2},1,\frac{1}{2^3},1,\ldots)$. 
Then the conclusion is the same as in Example \ref{rnb1}.
\end{ex}

\begin{ex}\label{rnb32}
Let $\Phi=\seq[e]$, $\Psi=(e_1,\frac{1}{2}e_2,e_1, \frac{1}{3}e_3,e_1, \frac{1}{4}e_4,\ldots)$, and $m=(\frac{1}{2},2,\frac{1}{2^2},3,\frac{1}{2^3},4,\ldots)$. 
Then the conclusion is the same as in Example \ref{rnb1}.

\end{ex}

\begin{ex}\label{rnb2}
Let $\Phi=\seq[e]$, $\Psi=(ne_n)$, and $m=(\frac{1}{n})$. Then $M_{m,\Phi,\Psi} \neweq M_{m,\Psi,\Phi} \neweq M_{(1),(e_n),(e_n)}=I$.
\end{ex}

\begin{ex}\label{rnb3}
Let $\Phi=\seq[e]$, $\Psi=(ne_n)$, and $m=(\frac{1}{n^2})$. 
Then $M_{m,\Phi,\Psi} \neweq M_{m,\Psi,\Phi}\neweq M_{(\frac{1}{n}),(e_n),(e_n)}=G_1$ 
- unconditionally convergent and non-invertible on $\h$ (see Lemma \ref{lemg}).
\end{ex}

\begin{ex}\label{rnb33}
Let  $\Phi=\seq[e]$, $\Psi=(e_1,\frac{1}{2}e_2,e_1, 3e_3,e_1, \frac{1}{4}e_4,e_1, 5e_5,\ldots)$, and $m=(\frac{1}{2},1,\frac{1}{2^2},\frac{1}{3},\frac{1}{2^3},1,\frac{1}{2^4},\frac{1}{5},\ldots)$. 
Then the conclusion is the same as in Example \ref{rnb1}.
\end{ex}

\begin{ex}\label{rnb4}
Let $\Phi=\seq[e]$, $\Psi=(e_1,\frac{1}{2}e_2,3e_3,\frac{1}{4}e_4,5e_5,\ldots)$, and $m=(1,2,\frac{1}{3},4,\frac{1}{5},\ldots)$. Then $M_{m,\Phi,\Psi} \neweq M_{m,\Psi,\Phi} \neweq M_{(1),(e_n),(e_n)}=I$.
\end{ex}

\begin{ex}\label{rnb5}
Let $\Phi=\seq[e]$, $\Psi=(e_1,\frac{1}{2^2}e_2,3e_3,\frac{1}{4^2}e_4,5e_5,\ldots)$, and $m=(1,2,\frac{1}{3^2},4,\frac{1}{5^2},\ldots)$. Then $M_{m,\Phi,\Psi} \neweq M_{m,\Psi,\Phi}\neweq M_{(\frac{1}{n}),(e_n),(e_n)}=G_1$ 
- unconditionally convergent and non-invertible on $\h$ (see Lemma \ref{lemg}).
\end{ex}

\subsection{Examples for one Riesz basis and one Bessel non-frame sequence; TABLE 8 
on page \pageref{table8}}

\begin{ex}\label{rb1}
Let $\Phi=\seq[e]$, $\Psi=(e_2, e_2, e_3, e_4, e_5, \ldots)$, and $m=(1)$. Then $M_{m,\Phi,\Psi}$ and $M_{m,\Psi,\Phi}$ are unconditionally convergent due to Proposition \ref{bmh} and non-invertible on $\h$ due to Proposition \ref{rc}.
\end{ex}

\begin{ex}\label{rb1b}
Let $\Phi=\seq[e]$, $\Psi=(e_2, e_2, e_3, e_4, e_5, \ldots)$, and  $m=(\frac{1}{n})$. Then the conclusion is the same as in Example \ref{rb1}. 
\end{ex}

\begin{ex}\label{rb2}
Let $\Phi=\seq[e]$, $\Psi=(\frac{1}{n}e_n)$, and $m=(1)$. Then $M_{m,\Phi,\Psi}\neweq M_{m,\Psi,\Phi}\neweq M_{(\frac{1}{n}),(e_n),(e_n)} =G_1$
- unconditionally convergent and non-invertible on $\h$ (see Lemma \ref{lemg}).
\end{ex}

\begin{ex}\label{rb3}
Let $\Phi=\seq[e]$, $\Psi=(\frac{1}{n}e_n)$, and  $m=(\frac{1}{n})$. Then $M_{m,\Phi,\Psi}\neweq M_{m,\Psi,\Phi}\neweq M_{(\frac{1}{n^2}),(e_n),(e_n)} 
 =G_2$
- unconditionally convergent and non-invertible on $\h$ (see Lemma \ref{lemg}).
\end{ex}

\begin{ex}\label{rb4}
Let $\Phi=\seq[e]$, $\Psi=(\frac{1}{n}e_n)$, and  $m=(n)$. Then $M_{m,\Phi,\Psi}\neweq M_{m,\Psi,\Phi} \neweq M_{(1),(e_n),(e_n)}=I$.
\end{ex}

\begin{ex}\label{rb5}
Let $\Phi=\seq[e]$, $\Psi=(\frac{1}{n^2}e_n)$, and $m=(n)$. Then $M_{m,\Phi,\Psi}\neweq M_{m,\Psi,\Phi}\neweq M_{(\frac{1}{n}),(e_n),(e_n)} = G_1$ 
- unconditionally convergent and non-invertible on $\h$ (see Lemma \ref{lemg}).
\end{ex}

\subsection{Examples for one Riesz basis and one overcomplete frame; TABLE 9 
on page \pageref{table9}}

\begin{ex}\label{rf1}
Let $\Phi=\seq[e]$, $\Psi=(e_1,e_1,e_2,e_3,e_4,\ldots)$, and $m=(1)$. Then $M_{m,\Phi,\Psi}$ and $M_{m,\Psi,\Phi}$ are unconditionally convergent due to Proposition \ref{bmh} and non-invertible on $\h$ due to Proposition \ref{rbis}.
\end{ex}

\begin{ex}\label{rf2}
Let  $\Phi=\seq[e]$,  $\Psi=(e_1,e_1,e_2,e_3,e_4,\ldots)$, and $m=(1,\frac{1}{2},1,\frac{1}{3},1,\frac{1}{4},\ldots)$.
Then the conclusion is the same as in Example \ref{rf1}.
\end{ex}

\begin{ex}\label{rf3}
Let $\Phi=\seq[e]$,    $\Psi=(\frac{1}{2}e_1,e_2,\frac{1}{2^2}e_1,e_3,\frac{1}{2^3}e_1,e_4,\ldots)$, and  $m=(1)$. 
Then the conclusion is the same as in Example \ref{rf1}.
\end{ex}

\begin{ex}\label{rf4}
Let $\Phi=\seq[e]$,    $\Psi=(\frac{1}{2}e_1,e_2,\frac{1}{2^2}e_1,e_3,\frac{1}{2^3}e_1,e_4,\ldots)$, and  $m=(1,\frac{1}{2},1,\frac{1}{3},1,\frac{1}{4},\ldots)$. 
Then the conclusion is the same as in Example \ref{rf1}.
\end{ex}

\begin{ex}\label{rf5} Let $\Phi=\seq[e]$, 
$\Psi=(\frac{1}{2}e_1,e_2,\frac{1}{2^3}e_1,e_3,\frac{1}{2^5}e_1,e_4,\frac{1}{2^7}e_1, e_5,\ldots)$, 
and $m=(1,1,2,1,2^2,1,2^3,1,\ldots)$. 
 Then $M_{m,\Phi,\Psi}$ and $M_{m,\Psi,\Phi}$ are 
unconditionally convergent on $\h$, because 
$\Phi$ and $m\Psi$ are Bessel for $\h$  (apply Prop. \ref{bmh}). Both $M_{m,\Phi,\Psi}$ and $M_{m,\Psi,\Phi}$ are non-invertible on $\h$ by Proposition \ref{rc}. 

\end{ex}

\subsection{Examples for two Riesz bases; TABLE 10 
on page \pageref{table10}}

\begin{ex} \label{rry} Consider $M_{(1),(e_n),(e_n)}=I$.
\end{ex}

\begin{ex} \label{rrn} Consider $M_{(\frac{1}{n}),(e_n),(e_n)}$  
and see Lemma \ref{lemg}.
\end{ex}

 {\bf Acknowledgments} 
The authors are grateful for the hospitality of the 
Institut de Physique Th\'eorique, Universit\'e Catholique de Louvain. 
The second author is also grateful for the hospitality of the Acoustics Research Institute and the support from the MULAC-project. She thanks 
 UACEG and the Department of Mathematics of UACEG for their support and openness  
  in order for the research on this paper to be done.

\bibliographystyle{plain}

\end{document}